\documentstyle[amscd]{amsart}
\numberwithin{equation}{section}
\theoremstyle{plain}
\newtheorem{thm}{Theorem}[section]
\newtheorem{cor}[thm]{Corollary}
\newtheorem{lem}[thm]{Lemma}
\newtheorem{prop}[thm]{Proposition}
\begin{document}
\title{$C^*$-algebras associated with textile  dynamical systems
}
\author{Kengo Matsumoto}
\address{ 
Department of Mathematics, 
Joetsu University of Education,
Joetsu 943-8512, Japan}
\email{kengo@@juen.ac.jp}
\maketitle
\begin{abstract}
A $C^*$-symbolic dynamical system
$({\cal A}, \rho,  \Sigma)$
is a finite family 
$\{ \rho_\alpha\}_{\alpha \in\Sigma}$
of endomorphisms of a $C^*$-algebra ${\cal A}$ 
with some conditions.
The endomorphisms yield a $C^*$-algebra ${\cal O}_\rho$
from the associated 
Hilbert $C^*$-bimodule.
In this paper, we will extend the notion of
$C^*$-symbolic dynamical system to 
$C^*$-textile  dynamical system
$({\cal A}, \rho, \eta, {\Sigma^\rho}, {\Sigma^\eta}, \kappa)$
which consists of two 
$C^*$-symbolic dynamical systems 
$({\cal A}, \rho,  {\Sigma^\rho})$
and 
$({\cal A}, \eta, {\Sigma^\eta})$
with certain commutation relations
$\kappa$
between their endomorphisms
$\{\rho_\alpha\}_{\alpha \in \Sigma^\rho}$
and
$\{\eta_a \}_{a \in \Sigma^\eta}$.
$C^*$-textile dynamical systems yield two-dimensional tilings and 
$C^*$-algebras ${\cal O}^{\kappa}_{\rho,\eta}$.
We will study the structure of the algebras ${\cal O}^\kappa_{\rho,\eta}$
and present its K-theory formulae .
 \end{abstract}


\def\Zp{{ {\Bbb Z}_+ }}
\def\CTDS{{ ({\cal A}, \rho, \eta, \Sigma^\rho, \Sigma^\eta, \kappa)}}
\def\OL{{{\cal O}_{{\frak L}}}}
\def\FL{{{\cal F}_{{\frak L}}}}
\def\FKL{{{\cal F}_k^{l}}}
\def\FHKL{{{\cal F}_{[k]}^{[l]}}}
\def\FKI{{{\cal F}_k^{\infty}}}
\def\FHKI{{{\cal F}_{[k]}^{\infty}}}
\def\FLI{{{\cal F}_{\frak L}}}
\def\FHLI{{{\cal F}_{[{\frak L}]} }}
\def\FN{{{\cal F}_{\theta_1,\dots,\theta_N}}}
\def\FNk{{{\cal F}^k_{\theta_1,\dots,\theta_N}}}
\def\FRHO{{ {\cal F}_\rho}}
\def\DRHO{{{{\cal D}_{\rho}}}}
\def\ORHO{{{\cal O}_\rho}}
\def\ORE{{{\cal O}^{\kappa}_{\rho,\eta}}}
\def\OHRE{{{\widehat{\cal O}}^{\kappa}_{\rho,\eta}}}
\def\FRE{{{\cal F}_{\rho,\eta}}}
\def\DRE{{{\cal D}_{\rho,\eta}}}
\def\M{{{\cal M}}}
\def\N{{{\cal N}}}
\def\H{{{\cal H}}}
\def\K{{{\cal K}}}
\def\P{{{\cal P}}}
\def\Q{{{\cal Q}}}
\def\A{{{\cal A}}}
\def\E{{ {\cal E} }}
\def\B{{{\cal B}}}
\def\R{{{\cal R}}}
\def\S{{{\cal S}}}
\def\sms{{{symbolic  matrix system }}}
\def\smss{{{symbolic  matrix systems }}}
\def\Aut{{{\operatorname{Aut}}}}
\def\End{{{\operatorname{End}}}}
\def\Ext{{{\operatorname{Ext}}}}
\def\Hom{{{\operatorname{Hom}}}}
\def\Ker{{{\operatorname{Ker}}}}
\def\ker{{{\operatorname{ker}}}}
\def\Coker{{{\operatorname{Coker}}}}
\def\id{{{\operatorname{id}}}}
\def\dim{{{\operatorname{dim}}}}
\def\min{{{\operatorname{min}}}}
\def\exp{{{\operatorname{exp}}}}
\def\Homeo{{{\operatorname{Homeo}}}}
\def\supp{{{\operatorname{supp}}}}
\def\Span{{{\operatorname{Span}}}}
\def\Proj{{{\operatorname{Proj}}}}
\def\Im{{{\operatorname{Im}}}}
\def\span{{{span}}}
\def\OA{{{\cal O}_A}}
\def\OB{{{\cal O}_B}}
\def\DA{{{\cal D}_A}}
\def\DB{{{\cal D}_B}}
\def\T{{  {\cal T}_{{\K}^{\M}_{\N}} }}


\section{Introduction}

In \cite{MaDoc99}, the author has introduced a notion of 
$\lambda$-graph system as presentations of subshifts.
The $\lambda$-graph systems are labeled Bratteli diagram with shift transformation.
They yield $C^*$-algebras so that its  K-theory groups  are
related to topological conjugacy invariants of the underlying symbolic dynamical systems.
The class of these $C^*$-algebras  include the Cuntz-Krieger algebras.
He has extended the notion of $\lambda$-graph system to  $C^*$-{\it symbolic dynamical system},
which 
is a generalization of both
a $\lambda$-graph system and an automorphism of a unital $C^*$-algebra.
It is a finite family $\{ \rho_{\alpha} \}_{\alpha \in \Sigma}$ 
of endomorphisms of a unital $C^*$-algebra 
${\cal A}$ 
such that the closed ideal generated by $\rho_\alpha(1), \alpha \in \Sigma$
coincides with $\A$.
A  finite labeled graph ${\cal G}$
gives rise to a $C^*$-symbolic dynamical system 
$({\cal A}_{\cal G}, \rho^{\cal G}, \Sigma)$ such that ${\cal A} = {\Bbb C}^N$ for some 
$N \in {\Bbb N}$.
A $\lambda$-graph system ${\frak L}$
is a generalization of a finite 
labeled graph and  yields a $C^*$-symbolic dynamical system 
$({\cal A}_{\frak L},\rho^{\frak L}, \Sigma)$ 
such that ${\cal A}_{\frak L}$ 
is $C(\Omega_{\frak L})$ for some compact Hausdorff space 
$\Omega_{\frak L}$ with 
$\dim \Omega_{\frak L} = 0$.
It also yields
a $C^*$-algebra $\OL$.
A $C^*$-symbolic dynamical system $({\cal A},\rho, \Sigma)$
provides a subshift 
denoted by $\Lambda_\rho$  
over $\Sigma$
and a Hilbert $C^*$-right $\cal A$-module 
$(\phi_\rho, {\cal H}_{\cal A}^{\rho}, \{u_\alpha\}_{\alpha \in \Sigma})$
with an orthogonal finite basis 
$\{u_{\alpha}\}_{\alpha\in \Sigma}$
and a unital faithful diagonal left action 
$\phi_{\rho} :{\cal A} \rightarrow L({\cal H}_{\cal A}^{\rho})$.
By using general construction of $C^*$-algebras 
 from Hilbert $C^*$-bimodules established by M. Pimsner \cite{Pim},
a $C^*$-algebra
denoted  by 
${\cal O}_\rho$
from 
$(\phi_{\rho},{\cal H}_{\cal A}^{\rho}, 
\{u_{\alpha}\}_{\alpha\in \Sigma})$ 
has been presented in \cite{MaCrelle}.
We call the algebra 
${\cal O}_\rho$ 
{\it the}\ $C^*$-{\it symbolic crossed product of}\ $\cal A$ {\it by the subshift}\ $\Lambda_\rho$. 
If ${\cal A} = C(X)$ with $\dim X=0$, 
there exists a $\lambda$-graph system ${\frak L}$ 
 such that the subshift $\Lambda_\rho$ 
is the subshift $\Lambda_{\frak L}$ 
presented by ${\frak L}$ 
and the $C^*$-algebra $\ORHO$ is the $C^*$-algebra ${\cal O}_{\frak L}$
associated with  $\frak L$.
If in particular, ${\cal A} = {\Bbb C}^N$, 
the subshift $\Lambda_\rho$ is a sofic shift
and 
$\ORHO$ is a Cuntz-Krieger algebra.
If $\Sigma = \{ \alpha \}$ an automorphism
$\alpha$ of a unital $C^*$-algebra $\A$,
the $C^*$-algebra $\ORHO$ is the ordinary crossed product
$\A\times_\alpha{\Bbb Z}$.

G. Robertson-T. Steger \cite{RoSte}
have initiated a certain study of higher dimensional analogue of Cuntz-Krieger algebras
from the view point of tiling systems of 2-dimensional plane.
After their work, 
A. Kumjian-D. Pask \cite{KP}
have generalized their construction to introduce the notion of 
higher rank graphs and its $C^*$-algebras.
The $C^*$-algebras constructed from 
higher rank graphs are called the higher rank graph $C^*$-algebras.
Since then, there have been 
many studies  on these $C^*$-algebras by many authors
(see for example \cite{Dea}, \cite{EGS}, \cite{KP}, \cite{PCHR}, \cite{PRW1}, \cite{RoSte}, etc.).

M. Nasu in \cite{NaMemoir} has introduced the notion of textile system 
which is useful in analyzing automorphisms and endomorphisms of topological Markov shifts.
A textile system also gives rise to a two-dimensional tiling called Wang tiling. 
Among textile systems,
 LR textile systems have specific properties 
that consist of two commuting symbolic matrices
$\M \P \cong \P \M$. 
In \cite{MaYMJ2008},
the author has extended the notion of textile systems to $\lambda$-graph systems
and has defined a notion of textile systems on $\lambda$-graph systems,
which are called textile $\lambda$-graph systems for short.
 $C^*$-algebras associated to  textile systems  
have been initiated by V. Deaconu (\cite{Dea}).

In this paper, we will extend the notion of
$C^*$-symbolic dynamical system to 
$C^*$-textile  dynamical system
which is a higher dimensional analogue of $C^*$-symbolic dynamical system.
The $C^*$-textile  dynamical system
$({\cal A}, \rho, \eta, \Sigma^\rho, \Sigma^\eta, \kappa)$
 consists of two 
$C^*$-symbolic dynamical systems 
$({\cal A}, \rho,  {\Sigma^\rho})$
and 
$({\cal A}, \eta, {\Sigma^\eta})$
with the following  commutation relations
between $\rho$ and $\eta$ through $\kappa$.
Set 
$$
\Sigma_{\rho \eta} = \{ (\alpha, b ) \in \Sigma^\rho \times \Sigma^\eta \mid 
\eta_b \circ \rho_\alpha  \ne 0 \},\quad
\Sigma_{\eta \rho} = \{ (a, \beta) \in \Sigma^\eta \times \Sigma^\rho \mid 
\rho_\beta \circ \eta_a  \ne 0 \}.
$$
We require that there exists a bijection
$
\kappa : \Sigma_{\rho \eta} \longrightarrow \Sigma_{\eta\rho},
$
which we fix and call a specification.
Then the required commutation relations are 
\begin{equation}
\eta_b \circ \rho_\alpha = \rho_\beta \circ \eta_a
\qquad
\text{ if } 
\quad 
\kappa(\alpha, b) = (a,\beta).
\end{equation}
$C^*$-textile dynamical systems provide  two-dimensional tilings and 
$C^*$-algebra $\ORE$.
The $C^*$-algbebra $\ORE$ 
is defined to be the universal $C^*$-algebra 
$C^*(x, S_\alpha, T_a ; x \in \A, \alpha \in \Sigma^\rho, a \in \Sigma^\eta)$  
generated by 
$x \in \A$ 
and two family of partial isometries 
$S_{\alpha}, \alpha \in \Sigma^\rho$,
$T_a, a \in \Sigma^\eta$
subject to the following relations called $(\rho,\eta;\kappa)$:
\begin{align}
\sum_{\beta \in \Sigma^\rho}S_{\beta}S_{\beta}^*  =1,\qquad
x S_\alpha S_\alpha^* & =  S_\alpha S_\alpha^*  x,\qquad
S_\alpha^* x S_\alpha  = \rho_\alpha(x),\\
\sum_{b \in \Sigma^\eta} T_{b} T_{b}^*  =1,\qquad
x T_a T_a^*  & =  T_a T_a^*  x,\qquad
T_a^* x T_a  = \eta_a(x),\\  
S_\alpha T_b  = T_a S_\beta &
\qquad
\text{ if } 
\quad
\kappa(\alpha, b) = (a,\beta)
\end{align}
for all $x \in {\cal A}$ and $\alpha \in \Sigma^\rho, a \in \Sigma^\eta.$

We will construct a tiling system in the plane from a
$C^*$-textle dynamical system.
The resulting tiling system is a two-dimensional subshift.

In this paper,
we will study the $C^*$-algebra
$\ORE$.
We will introduce a condition called (I) on $\CTDS$ 
which will be studied 
as a generalization of 
the condition (I)  on $C^*$-symbolic dynamical system
 \cite{Ma11}(cf. \cite{Ma7})
(and hence on a finite matrix of Cuntz-Krieger \cite{CK} ).
Under the assumption that $\CTDS$ satisfies condition (I),
the simplicity conditions of 
the algebra $\ORE$ will be  discussed in Section 4.
We will show the following 
 
\begin{thm}
Let $\CTDS$  be 
a $C^*$-textile  dynamical system satisfying condition (I).
Then the $C^*$-algebra $\ORE$ is the unique $C^*$-algebra subject to the relations
$(\rho,\eta;\kappa)$.
If $\CTDS$ is irreducible,  $\ORE$ is simple.
\end{thm}

We denote by 
$Z_\A$ the center of $\A$.
We next assume that 
$\rho_\alpha(Z_\A) \subset Z_\A, \alpha \in \Sigma^\rho$
and
$\eta_a(Z_\A) \subset Z_\A, a \in \Sigma^\eta$.
All examples of $C^*$-symbolic dynamical systems 
$(\A,\rho,\Sigma)$
appearing in the previous papers
\cite{MaCrelle}, \cite{MaMZ2010}
satisfy the conditions
$\rho_\alpha(Z_\A) \subset Z_\A, \alpha \in \Sigma$. 
A $C^*$-textile  dynamical system
$\CTDS$ is said to form squares
if 
the $C^*$-subalgebra of $\A$ generated by projections
$\rho_\alpha(1), \alpha \in \Sigma^\rho$
and
the $C^*$-subalgebra of $\A$ generated by projections
$\eta_a(1), a \in \Sigma^\eta$
coincide.
We will restrict our interest to the $C^*$-textile dynamical systems forming squares.
Then the  K-theory formulae hold as in the following way:
\begin{thm}
Suppose that $\CTDS$ froms squares.
There exists short exact sequences for $K_0(\ORE)$ and $K_1(\ORE)$ such as 
\begin{align*}
0 
& \longrightarrow 
K_0(\A) /
(1 - \lambda_\eta)K_0(\A) + (1 - \lambda_\rho)K_0(\A) \\
& \longrightarrow
K_0(\ORE) \\
& \longrightarrow
\Ker(1 - \lambda_\eta) \cap \Ker(1 - \lambda_\rho)
\text{ in }
K_0(\A)
\longrightarrow 0 \\
\intertext{ and}
0 
& \longrightarrow 
\Ker(1-\lambda_\eta) 
\text{ in }
K_0(\A) / 
(1 - \lambda_\rho) 
(\Ker(1-\lambda_\eta) 
\text{ in }
K_0(\A)) \\
& \longrightarrow
K_1(\ORE) \\
& \longrightarrow
\Ker(1 - \lambda_\rho)\text{ in }
K_0(\A) / (1 - \lambda_\eta)K_0(\A)
\longrightarrow 0
\end{align*}
where the endomorphisms
$
\lambda_\rho, \lambda_\eta: 
K_0(\A) \longrightarrow K_0(\A)
$
are defined by
\begin{align*}
\lambda_\rho([p]) & = \sum_{\alpha \in \Sigma^\rho} [\rho_\alpha(p)] \in K_0(\A)
\text{ for } [p] \in K_0(\A),\\
\lambda_\eta([p]) & = \sum_{a \in \Sigma^\eta} [\eta_a(p)] \in K_0(\A)
\text{ for } [p] \in K_0(\A).  
\end{align*}
\end{thm} 
Let
$A, B$ be  $N\times N$ matrices with entries in nonnegative integers such that
\begin{equation*}
AB = BA.
\end{equation*} 
Let $G_A, G_B$ be labeled directed graphs whose transition matrices are $A, B$ respectively.
Let $\M_A, \M_B$ denote symbolic marices for  $G_A, G_B$ 
whose conponents consist of formal sums of directed edges  respectively.
By the condition $AB = BA$,
one may take  
 a bijection
$\kappa: \Sigma^{AB} \longrightarrow \Sigma^{BA}$
which gives rise to
a specified equivalence
$\M_A \M_B \overset{\kappa}{\cong}\M_A \M_B$.
We then have a $C^*$-textile dynamical system
written as 
\begin{equation*}
(\A, \rho^{A}, \rho^{B},\Sigma^A, \Sigma^B, \kappa).
\end{equation*}
The associated $C^*$-algebra
is denoted by
${\cal O}_{A,B}^{\kappa}$.
The $C^*$-algebra
${\cal O}_{A,B}^{\kappa}$
is realized as a $2$-graph $C^*$-algebra constructed from Kumjian and Pask.
It is also seen in Deaconu's paper in \cite{Dea}.

\begin{prop}
Keep the above situations.
There exist short  exact sequences: 
\begin{enumerate}
\renewcommand{\labelenumi}{(\roman{enumi})}
\item
\begin{align*}
0 
& \longrightarrow 
{\Bbb Z}^N /
( (1 - A){\Bbb Z}^N + (1 - B) {\Bbb Z}^N ) \\
& \longrightarrow
K_0({\cal O}_{A,B}^{\kappa}) \\
& \longrightarrow
\Ker(1 - A) \cap \Ker(1 - B)
\text{ in }
{\Bbb Z}^N
\longrightarrow 0
\end{align*}
\item
\begin{align*}
0 
& \longrightarrow 
\Ker(1-B) \text{ in } {\Bbb Z}^N /
(1 - A)(\Ker(1 - B) \text{ in } {\Bbb Z}^N )\\
& \longrightarrow
K_1({\cal O}_{A,B}^{\kappa}) \\
& \longrightarrow
\Ker(1 - \Bar{A}) 
\text{ in }
{\Bbb Z}^N /
(1 - B) {\Bbb Z}^N 
\longrightarrow 0,
\end{align*}
\end{enumerate}
where $\Bar{A}$ is an endomorphism on the abelian group
$
{\Bbb Z}^N /
(1 - B) {\Bbb Z}^N 
$
induced by 
the matrix $A$.
\end{prop}

Throughout the paper,
we will denote by $\Zp$ 
and
by ${\Bbb N}$
the sets of nonnegative integers 
and the set of positive integers respectively.


\section{ $\lambda$-graph systems, $C^*$-symbolic dynamical systems
and their $C^*$-algebras}
In this section, we will briefly review 
 $\lambda$-graph systems and $C^*$-symbolic dynamical systems.
Throughout the section,
$\Sigma$ denotes a finite set with its discrete topology, 
that is called an alphabet.
Each element of $\Sigma$ is called a symbol.
Let $\Sigma^{\Bbb Z}$ 
be the infinite product space 
$\prod_{i\in {\Bbb Z}}\Sigma_{i},$ 
where 
$\Sigma_{i} = \Sigma$,
 endowed with the product topology.
 The transformation $\sigma$ on $\Sigma^{\Bbb Z}$ 
given by 
$
\sigma( (x_i)_{i \in {\Bbb Z}})
 = (x_{i+1})_{i \in {\Bbb Z}}
$ 
is called the full shift over $\Sigma$.
 Let $\Lambda$ be a shift invariant closed subset of $\Sigma^{\Bbb Z}$ i.e. 
 $\sigma(\Lambda) = \Lambda$.
  The topological dynamical system 
  $(\Lambda, \sigma\vert_{\Lambda})$
   is called a two-sided subshift, written as $\Lambda$ for brevity.

There is a class of subshifts called sofic shifts, 
that are presented by finite labeled graphs.
$\lambda$-graph systems are generalization of finite labeled graphs.
Any subshift is presented by a  
$\lambda$-graph system.
Let ${\frak L} =(V,E,\lambda,\iota)$ be 
a $\lambda$-graph system 
 over $\Sigma$ with vertex set
$
V = \cup_{l \in \Zp} V_{l}
$
and  edge set
$
E = \cup_{l \in \Zp} E_{l,l+1}
$
that is labeled with symbols in $\Sigma$ by a map
$\lambda: E \rightarrow \Sigma$, 
and that is supplied with  surjective maps
$
\iota( = \iota_{l,l+1}):V_{l+1} \rightarrow V_l
$
for
$
l \in \Zp.
$
Here the vertex sets $V_{l},l \in \Zp$
are finite disjoint sets.   
Also  
$E_{l,l+1},l \in \Zp$
are finite disjoint sets.
An edge $e$ in $E_{l,l+1}$ has its source vertex $s(e)$ in $V_{l}$ 
and its terminal  vertex $t(e)$ 
in
$V_{l+1}$
respectively.
Every vertex in $V$ has a successor and  every 
vertex in $V_l$ for $l\in {\Bbb N}$ has a predecessor. 
It is then required that 
for vertices $u \in V_{l-1}$
and $v \in V_{l+1}$,
there exists a bijective correspondence
between the set of edges
$e \in E_{l,l+1}$
such that
$t(e) = v, \iota(s(e)) = u$
and the set of edges 
$f \in E_{l-1,l}$
such that 
$\ s(f) = u, t(f) = \iota(v)$,
preserving thier labels (\cite{Ma3}).
We henceforth assume that $\frak L$ is left-resolving, which means that 
$t(e)\ne t(f)$ whenever $\lambda(e) = \lambda(f)$ for $e,f \in E_{l,l+1}$.
Let us denote by 
$\{v_1^l,\dots, v_{m(l)}^l\}$
the vertex set
$V_l$ at level $l$.
For
$
i=1,2,\dots,m(l),\ j=1,2,\dots,m(l+1), \ \alpha \in \Sigma
$ 
we put
\begin{align*}
A_{l,l+1}(i,\alpha,j)
 & =
\begin{cases}
1 &  
    \text{ if } \ s(e) = {v}_i^l, \lambda(e) = \alpha,
                       t(e) = {v}_j^{l+1} 
    \text{ for some }    e \in E_{l,l+1}, \\
0           & \text{ otherwise,}
\end{cases} \\
I_{l,l+1}(i,j)
 & =
\begin{cases}
1 &  
    \text{ if } \ \iota_{l,l+1}({v}_j^{l+1}) = {v}_i^l, \\
0           & \text{ otherwise.}
\end{cases} 
\end{align*}
The $C^*$-algebra 
${\cal O}_{\frak L}$ 
associated with ${\frak L}$ 
is the  universal    
$C^*$-algebra generated by partial isometries
$S_\alpha, \alpha \in \Sigma$
and projections $E_i^l, i=1,2,\dots,m(l),\ l \in \Zp$  
subject to the following  operator relations called $({\frak L})$:
\begin{align}
& \sum_{\beta \in \Sigma}  S_{\beta}S_{\beta}^*  =  1, \\
  \sum_{i=1}^{m(l)} E_i^l  & =  1, \qquad 
 E_i^l   =  \sum_{j=1}^{m(l+1)}I_{l,l+1}(i,j)E_j^{l+1}, \\
& S_\alpha S_\alpha^* E_i^l =E_i^l 
S_\alpha S_\alpha^*,  \\ 
S_{\alpha}^*E_i^l S_{\alpha} & =  
\sum_{j=1}^{m(l+1)} A_{l,l+1}(i,\alpha,j)E_j^{l+1},
\end{align}
for
$
i=1,2,\dots,m(l),\l\in \Zp, 
 \alpha \in \Sigma.
 $
If $\frak L$ satisfies  $\lambda$-condition $(I)$ 
and is $\lambda$-irreducible,
the $C^*$-algebra ${\cal O}_{\frak L}$
is simple and purely infinite
(\cite{Ma11}, \cite{Ma7}).


Let $\A_{{\frak L},l}$
be the $C^*$-subalgebra of ${\cal O}_{\frak L}$
generated by the projections
$E_i^l, i=1,\dots,m(l)$.
We denote by 
$\A_{\frak L}$ the $C^*$-subalgebra of ${\cal O}_{\frak L}$
generated by the all projections
$E_i^l, i=1,\dots,m(l), l\in \Zp.$
We denote by 
$\iota: \A_{{\frak L},l} \hookrightarrow \A_{{\frak L},l+1}$
the natural inclusion.
Hence the algebra $\A_{\frak L}$
is the inductive limit
$\underset{\iota}{\varinjlim} \A_{{\frak L},l}$
of the inclusions.
For $\alpha \in \Sigma$, put
$$
\rho^{\frak L}_\alpha(X) =  S_\alpha^* X S_\alpha \qquad 
\text{ for }
\quad
X \in \A_{\frak L}.
$$ 
Then
$
\{\rho^{\frak L}_\alpha\}_{\alpha \in \Sigma}
$
yields a family of $*$-endomorphisms of 
$\A_{\frak L}$
such that
$\rho^{\frak L}_\alpha(1) \ne 0,$
$\sum_{\alpha \in \Sigma}\rho^{\frak L}_\alpha(1) \ge 1$
and for any nonzero $x \in {\A}_{\frak L}$, 
$\rho^{\frak L}_\alpha(x) \ne 0$ for some $\alpha \in \Sigma$.

\medskip

The situations above are generalized to
$C^*$-symbolic dynamical systems as follows.

Let ${\cal A}$ be a unital $C^*$-algebra.
In what follows,
an endomorphism of $\cal A$ means 
a $*$-endomorphism of $\cal A$ that does not necessarily preserve the unit
$1_\A$ of 
$\cal A$.
The unit $1_\A$ is denoted by $1$ unless we specify.
For an alphabet $\Sigma$,
a finite family of endomorphisms 
$
\rho_\alpha, \alpha \in \Sigma
$
of $\A$    
is said to be {\it essential}\ if
$\rho_{\alpha}(1) \ne 0$ for all $\alpha \in \Sigma$
 and   
the closed ideal generated by $\rho_\alpha(1), \alpha \in \Sigma$
coincides with $\A$.
It is said to be {\it faithful}\ if for any nonzero $x \in \cal A$ 
there exists a symbol $\alpha\in \Sigma$ such that $\rho_{\alpha}(x) \ne 0$.

\noindent
{\bf Definition (\cite{MaCrelle}).}
A $C^*$-{\it symbolic dynamical system}\ 
is a triplet $({\cal A}, \rho, \Sigma)$ 
consisting of a unital $C^*$-algebra $\A$ 
and 
an essential and faithful finite family 
$\{ \rho_{\alpha} \}_{\alpha \in \Sigma}$  of endomorphisms
of ${\cal A}$.
A $C^*$-symbolic dynamical system $({\cal A}, \rho, \Sigma)$
 yields  a subshift $\Lambda_\rho$
over $\Sigma$ such that a word 
$\alpha_1\cdots\alpha_k$ of $\Sigma$ is admissible for $\Lambda_\rho$
 if and only if 
$(\rho_{\alpha_k}\circ \cdots\circ \rho_{\alpha_1})(1) \ne 0$
(\cite[Proposition 2.1]{MaCrelle}).
Denote by 
$B_k(\Lambda_\rho)$
 the set of admissible words of
$\Lambda_\rho$ respectively with length $k$.
Put
$
B_*(\Lambda_\rho) = \cup_{k=0}^{\infty}B_k(\Lambda_\rho),
$
where $B_0(\Lambda_\rho), B_0(\Lambda_\eta) $ denote the empty word. 
We say that {\it a subshift}\ $\Lambda$ {\it acts on a}\ $C^*$-{\it algebra}\ 
${\cal A}$ if there exists a $C^*$-symbolic dynamical system 
$({\cal A}, \rho, \Sigma)$ 
such that the associated subshift $\Lambda_{\rho}$
is $\Lambda$.
A $C^*$-symbolic dynamical system 
$(\A,\rho,\Sigma)$ is said to be {\it central}
if $\rho_\alpha(Z_\A) \subset Z_\A$ for all 
$\alpha \in \Sigma$.
In this case,
essentiality of the endomorphisms
$\rho_\alpha, \alpha \in \Sigma$
is equivalent to the condition that
$\rho_\alpha(1) \ne 0, \alpha \in \Sigma$
and
the inequality
$\sum_{\alpha \in \Sigma}\rho_\alpha(1) \ge 1$
holds.
All of the examples appeared in the papers
\cite{MaCrelle}, \cite{MaMZ2010}
are central in this sense.
We will henceforth assume that
$C^*$-symbolic dynamical systems are all central. 

As in the above discussion
we have  a 
$C^*$-symbolic dynamical system 
$({\cal A}_{\frak L}, \rho^{\frak L}, \Sigma)$
from a $\lambda$-graph system ${\frak L}$
such that
 the $C^*$-algebra ${\cal A}_{\frak L}$ is $C(\Omega_{\frak L})$ 
 with $\dim \Omega_{\frak L} =0$, 
and the subshift 
$\Lambda_{\rho^{\frak L}}$
coincides with the subshift $\Lambda_{\frak L}$ presented by $\frak L$. 
Conversely, for a 
$C^*$-symbolic dynamical system
 $({\cal A}, \rho, \Sigma)$,
 if the algebra ${\cal A}$ is $C(X)$ with $\dim X =0$,
there exists a $\lambda$-graph system $\frak L$ over $\Sigma$
such that the associated $C^*$-symbolic dynamical system 
$({\cal A}_{\frak L}, \rho^{\frak L}, \Sigma)$
is isomorphic to 
$({\cal A}, \rho, \Sigma)$ (\cite[Theorem 2.4]{MaCrelle}).

The $C^*$-algebra $\ORHO$ 
associated with a $C^*$-symbolic dynamical system
$(\A,\rho,\Sigma)$
has been originally constructed in \cite{MaCrelle} 
as a $C^*$-algebra 
by using the Pimsner's general construction of 
$C^*$-algebras from  Hilbert $C^*$-bimodules
\cite{Pim} 
(cf. \cite{KPW} etc.).
It is called 
{\it the}\ $C^*$-{\it symbolic crossed product of}\ ${\cal A}$ {\it by the subshift}\ 
$\Lambda_\rho$,
and realized as the universal $C^*$-algebra
$C^*(x, S_{\alpha}; x \in \A, \alpha \in \Sigma)$  
generated by 
$x \in \A$ 
and partial isometries $S_{\alpha}, \alpha \in \Sigma$
subject to the following relations called $(\rho)$:
$$
\sum_{\beta \in \Sigma}S_{\beta}S_{\beta}^* =1,\qquad
x S_\alpha S_\alpha^* =  S_\alpha S_\alpha^*  x,\qquad
S_\alpha^* x S_\alpha = \rho_\alpha(x)  
$$
for all $x \in \cal A$ and $\alpha \in \Sigma.$
Furthermore for  $\alpha_1,\dots,\alpha_k \in \Sigma$,  
a word $(\alpha_1,\dots,\alpha_k)$ is admissible for the subshift
$\Lambda_\rho$ 
if and only if
$S_{\alpha_1}\cdots S_{\alpha_k} \ne 0$
(\cite[Proposition 3.1]{MaCrelle}).
The $C^*$-algebra $\ORHO$ is a generalization of 
the $C^*$-algebra ${\cal O}_{\frak L}$
associated with the $\lambda$-graph system $\frak L$.

Let $\alpha$ be an automorphism of a unital $C^*$-algebra $\A$.
Put $\Sigma = \{ \alpha \}$ and $\rho_{\alpha} = \alpha$.
The $C^*$-algebra $\ORHO$
for the $C^*$-symbolic dynamical system $(\A, \rho, \Sigma)$ 
is  the ordinary crossed product $\A\times_{\alpha}{\Bbb Z}$.

\section{$C^*$-textile dynamical systems and their $C^*$-algebras}
Let 
$\CTDS$ be a 
$C^*$-textile dynamical system.
It  consists of two 
$C^*$-symbolic dynamical systems 
$({\cal A}, \rho,  {\Sigma^\rho})$
and 
$({\cal A}, \eta, {\Sigma^\eta})$
with the following  commutation relations
through $\kappa$.
Set 
$$
\Sigma_{\rho \eta} = \{ (\alpha, b ) \in \Sigma^\rho \times \Sigma^\eta \mid 
\eta_b \circ \rho_\alpha  \ne 0 \},\quad
\Sigma_{\eta \rho} = \{ (a, \beta) \in \Sigma^\eta \times \Sigma^\rho \mid 
\rho_\beta \circ \eta_a  \ne 0 \}.
$$
Let 
$
\kappa : \Sigma_{\rho \eta} \longrightarrow \Sigma_{\eta\rho}
$
be a bijection, which is called a specification.
Then the required commutation relations are 
\begin{equation}
\eta_b \circ \rho_\alpha = \rho_\beta \circ \eta_a
\qquad
\text{ if } 
\quad
\kappa(\alpha, b) = (a,\beta).
\end{equation}
$C^*$-textile dynamical systems will yield 
a two-dimensional subshift 
$X_{\rho,\eta}^\kappa$
and 
a $C^*$-algebra $\ORE$.

\medskip

Let $\Sigma$ be a finite set.
 The two-dimensional full shift over $\Sigma$ is defined to be 
\begin{equation*}
\Sigma^{{\Bbb Z}^2}
= \{ (x_{i,j})_{(i,j) \in {\Bbb Z}^2} \mid x_{i,j} \in \Sigma \}.
\end{equation*}
An element $x \in \Sigma^{{\Bbb Z}^2}$
is regarded as a function
$x: {\Bbb Z}^2 \longrightarrow \Sigma
$
which is called a configuration on ${\Bbb Z}^2$.
For $x \in \Sigma^{{\Bbb Z}^2}$
and
$F \subset {\Bbb Z}^2$,
let
$x_F$ denote the restriction of $x$ to $F$.
For a vector $m=(m_1,m_2) \in {\Bbb Z}^2$,
let 
$
\sigma^m :  \Sigma^{{\Bbb Z}^2} \longrightarrow \Sigma^{{\Bbb Z}^2}
$
be
the translation along vector $m$ defined by
\begin{equation*}
\sigma^m ((x_{i,j})_{(i,j) \in {\Bbb Z}^2})
= (x_{i+{m_1}, j+{m_2}})_{(i,j) \in {\Bbb Z}^2}. 
\end{equation*}
A subset 
$X \subset \Sigma^{{\Bbb Z}^2}$ 
is said to be  translation invariant if
$\sigma^m(X) = X$ for all $m \in {\Bbb Z}^2$.
 It is obvious to see that 
a subset 
$X \subset \Sigma^{{\Bbb Z}^2}$ 
is  translation invariant if
 ond only if 
 $X$ is invariant only both horizontaly and vertically,
 that is,
 $\sigma^{(1,0)}(X) = X$ and
$\sigma^{(0,1)}(X) = X$.
For $k \in {\Bbb Z}$, put
 \begin{equation*}
[-k,k]^2 
= \{ (i,j) \in {\Bbb Z}^2 \mid -k\le i,j\le k \}
= [-k,k] \times [-k,k].
 \end{equation*}
A metric $d$ on
$\Sigma^{{\Bbb Z}^2}$ is defined by
for $x,y \in \Sigma^{{\Bbb Z}^2}$
with $x \ne y$ 
 \begin{equation*}
d(x,y)  
= 
\frac{1}{2^k}  \quad \text{ if } \quad x_{(0,0)} = y_{(0,0)},
 \end{equation*}
where
$k = \max \{k \in {\Bbb Z}_+ \mid  
x_{[-k,k]^2} = y_{[-k,k]^2} \}.
$ 
If $x_{(0,0)} \ne y_{(0,0)}$, put $k= -1$ on the above definition.
If $x = y$, we set $d(x,y) =0$.
A two-dimensional subshift $X$ is a closed,
translation invariant subset of $\Sigma^{{\Bbb Z}^2}$
(cf. \cite[p.467]{LM}).
There is an equivalent definition of two dimensional subshift based on
lists of forbidden patterns as follows:
A shape is a finite subset
$F \subset {\Bbb Z}^2$.
A pattern $f$ on a shape $F$ is a function
$f: F \longrightarrow \Sigma$.
For a list ${\frak F}$ of patterns,
put
\begin{equation*}
X_{\frak F} 
= \{ (x_{i,j})_{(i,j)\in {\Bbb Z}^2} \mid 
\sigma^m(x)_F \not\in {\frak F} \text{ for all } m \in {\Bbb Z}^2
\text{ and } F \subset {\Bbb Z}^2 \}. 
\end{equation*}
It is well-known that 
a subset $X \subset \Sigma^{{\Bbb Z}^2}$
is a two-dimensional subshift if and only if 
there exists a list of patterns ${\frak F}$ such that
$X=  X_{\frak F}$.

We will define a certain property of two-dimensional subshift as follows:

\noindent
{\bf Definition.}
A two-dimensional subshift $X$ is said to 
have {\it diagonal property}\
if for 
$ (x_{i,j})_{(i,j)\in {\Bbb Z}^2}, (y_{i,j})_{(i,j)\in {\Bbb Z}^2}\in X$,
the conditions
$
x_{i,j} = y_{i,j},
x_{i+1,j-1} = y_{i+1,j-1}
$
imply 
$
x_{i,j-1} = y_{i,j-1},
x_{i+1,j} = y_{i+1,j}.
$
A two-dimensional subshift having diagonal property is called 
{\it textile dynamical system}.

\begin{lem}
If a two dimensional subshift $X$ has diagonal propety,
then for $x \in X$ and $(i,j) \in {\Bbb Z}^2$,
the configuration $x$ is determined by the diagonal line
$(x_{i+n, j-n})_{n \in {\Bbb Z}}$ 
through   
$(i,j)$.
\end{lem}
\begin{pf}
By the diagonal property,
the sequence
$(x_{i+n, j-n})_{n \in {\Bbb Z}}$ 
determines both the sequences
$(x_{i+1+n, j-n})_{n \in {\Bbb Z}}$ 
and
$(x_{i-1+n, j-n})_{n \in {\Bbb Z}}$.
Repeating this way,
one sees that
the sequence
$(x_{i+n, j-n})_{n \in {\Bbb Z}}$ 
determines 
$x_{n,m}$ for all $(n,m) \in {\Bbb Z}^2$.
\end{pf}

Let $\CTDS$
be a $C^*$-textile dynamical system.
We set
\begin{equation*}
\Sigma_\kappa = \{ \omega= (\alpha, b, a,\beta)
\in \Sigma^\rho\times \Sigma^\eta\times \Sigma^\eta\times \Sigma^\rho
 \mid
\kappa(\alpha, b) = (a, \beta) \}  
\end{equation*}
For 
$\omega= (\alpha, b, a,\beta)$,
since
$
\eta_b \circ \rho_\alpha = \rho_\beta \circ \eta_a 
$
as endomorphism on $\A$,
one may identify 
the quadruplet 
$ (\alpha, b, a,\beta)$
with the endomorphism
$
\eta_b \circ \rho_\alpha(= \rho_\beta \circ \eta_a) 
$
on $\A$
which we will denote by
simply $\omega$.
Define maps
$
t,b :\Sigma_\kappa \longrightarrow \Sigma^\rho
$
and
$
l,r :\Sigma_\kappa \longrightarrow \Sigma^\rho
$
by 
setting 
\begin{equation*}
t(\omega) = \alpha, \quad b(\omega) = \beta, \quad
l(\omega) = a, \quad r(\omega) = b
\end{equation*}
$$
\begin{CD}
\cdot @>\alpha=t(\omega)>> \cdot \\
@V{a=l(\omega)}VV  @VV{b=r(\omega)}V \\
\cdot @>>\beta= b(\omega)> \cdot 
\end{CD}
$$
A configuration
$(\omega_{i,j})_{(i,j)\in{\Bbb Z}^2} \in \Sigma_\kappa^{{\Bbb Z}^2}$
is said to be {\it paived}\
if the following conditions hold
\begin{equation*}
t(\omega_{i,j}) = b(\omega_{i,j+1}), \quad 
r(\omega_{i,j}) = l(\omega_{i+1,j}), \quad
l(\omega_{i,j}) = r(\omega_{i-1,j}), \quad
b(\omega_{i,j}) = t(\omega_{i,j-1})
\end{equation*}
for all $(i,j) \in {\Bbb Z}^2$.

For 
a textile dynamical system 
$\CTDS$,
we set
\begin{align*}
X_{\rho,\eta}^\kappa
=
\{ & (\omega_{i,j})_{(i,j) \in {\Bbb Z}^2} \in \Sigma_\kappa^{{\Bbb Z}^2} 
\mid 
(\omega_{i,j})_{(i,j) \in {\Bbb Z}^2} \text{ is paved and } \\
& \omega_{i+n,j-n}\circ\omega_{i+n-1,j-n+1}\circ\cdots
\circ\omega_{i+1,j-1}\circ \omega_{i,j} \ne 0 
\text{ for all } (i,j) \in {\Bbb Z}^2, n \in {\Bbb N}
\}
\end{align*}
\begin{lem}
Suppose that $(\omega_{i,j})_{(i,j)\in {\Bbb Z}^2}$ 
is paved.
Then 
$(\omega_{i,j})_{(i,j)\in {\Bbb Z}^2} \in X_{\rho,\eta}^\kappa$
if and only if 
\begin{equation*}
\rho_{b(\omega_{i+n,j-m})} \circ \cdots \circ
\rho_{b(\omega_{i+1,j-m})} \circ
\rho_{b(\omega_{i,j-m})} \circ
\eta_{l(\omega_{i,j-m})} \circ \cdots
\eta_{l(\omega_{i,j-1})} \circ
\eta_{l(\omega_{i,j})} \ne 0
\end{equation*}
for all $(i,j) \in {\Bbb Z}^2$, $n,m \in \Zp$.
\end{lem}
$$
\begin{CD}
\cdot                @.     @.   @.  @. \\
@V{l(\omega_{i,j})}VV     @.  @.   @.    \\
\cdot                @.     @.   @.  @. \\
@V{l(\omega_{i,j-1})}VV   @.  @.   @.    \\
\cdot                @.     @.   @.  @. \\
\vdots                   @.  @.   @.    \\
\cdot                @.     @.   @.  @. \\
@V{l(\omega_{i,j-m})}VV   @.  @.   @.    \\
\cdot               @>>{b(\omega_{i,j-m})}> \cdot@>>{b(\omega_{i+1,j-m})}> \cdots
@>>{b(\omega_{i+n,j-m})}> \cdot  
\end{CD}
$$
\begin{pf}
Suppose that
$(\omega_{i,j})_{(i,j)\in {\Bbb Z}^2} \in X_{\rho,\eta}^\kappa$.
For $(i,j) \in {\Bbb Z}^2$, $n,m \in \Zp$,
we may assume that $m \ge n$.
Since
\begin{align*}
0 \ne & \omega_{i+m,j-m}\circ \cdots \circ  
        \omega_{i+n+1,j-m}\circ 
        \omega_{i+n,j-m}\circ \cdots \circ
        \omega_{i,j-m}\circ \cdots \circ
        \omega_{i+1,j-1}\circ 
        \omega_{i,j}\\
   =  & \omega_{i+m,j-m}\circ \cdots \circ  
        \omega_{i+n+1,j-m}\circ
        \rho_{b(\omega_{i+n,j-m})}\circ \cdots \circ  
        \rho_{b(\omega_{i+1,j-m})}\circ  
        \rho_{b(\omega_{i,j-m})} \\
     &  \circ \eta_{l(\omega_{i,j-m})}\cdots \circ
        \eta_{l(\omega_{i,j-m})}\circ \cdots \circ
        \eta_{l(\omega_{i,j-1})}\circ 
        \eta_{l(\omega_{i,j})},
\end{align*}
one has 
\begin{equation*}
\rho_{b(\omega_{i+n,j-m})} \circ \cdots \circ
\rho_{b(\omega_{i+1,j-m})} \circ
\rho_{b(\omega_{i,j-m})} \circ
\eta_{l(\omega_{i,j-m})} \circ \cdots
\eta_{l(\omega_{i,j-1})} \circ
\eta_{l(\omega_{i,j})} \ne 0.
\end{equation*}
Converse implications is clear by the equality:
\begin{align*}
& \omega_{i+n,j-n}\circ \cdots \circ  
        \omega_{i,j-n}\circ \cdots \circ
        \omega_{i,j-1}\circ 
        \omega_{i,j}\\
 & =    \rho_{b(\omega_{i+n,j-n})}\circ \cdots \circ  
        \rho_{b(\omega_{i,j-n})}\circ \eta_{l(\omega_{i,j-n})}\cdots \circ
        \eta_{l(\omega_{i,j-1})}\circ 
        \eta_{l(\omega_{i,j})}.
\end{align*}
\end{pf}
\begin{prop}
For $\CTDS$,
$X_{\rho,\eta}^\kappa$ 
is a two-dimensional subshift having diagonal property,
that is, 
$X_{\rho,\eta}^\kappa$ is a textile dynamical system.
\end{prop}
\begin{pf}
It is 
easy to see that the set
$$
E = \{
(\omega_{i,j})_{(i,j)\in {\Bbb Z}^2} \in {\Sigma}_\kappa^{{\Bbb Z}^2} 
\mid
(\omega_{i,j})_{(i,j)\in {\Bbb Z}^2} \text{ is paved } \}
$$
is closed, because its complement is open. 
The following set 
\begin{align*}
U = \{
 (\omega_{i,j})_{(i,j)\in {\Bbb Z}^2}\in {\Sigma}_\kappa^{{\Bbb Z}^2} \mid
& \omega_{k+n,l-n}\circ \omega_{k+n-1,l-n+1}\circ \cdots 
   \circ \omega_{k+1,l-1}\circ\omega_{k,l} =0
   \\
& \text{ for some } (k,l) \in {\Bbb Z}^2, n \in {\Bbb N} \}
\end{align*} 
is open in ${\Sigma}_\kappa^{{\Bbb Z}^2}$.
As
the equality
$X_{\rho,\eta}^\kappa = E \cap U^{c}$ holds,
the set 
$X_{\rho,\eta}^\kappa$ 
is closed.
It is also obvious that 
$X_{\rho,\eta}^\kappa$ 
is translation invariant so that 
$X_{\rho,\eta}^\kappa$ 
is a textile dynamical system.
It is easy to see that 
$X_{\rho,\eta}^\kappa$ 
has diagonal property.
\end{pf}
We call 
$X_{\rho,\eta}^\kappa$ 
the textile dynamical system associated with
$\CTDS$.

Let us now define a subshift 
$X_{\delta^\kappa}$ 
over $\Sigma_\kappa$,
which consists of diagonal sequences of 
$X_{\rho,\eta}^\kappa$ 
as follows:
\begin{equation*}
X_{\delta^\kappa} = 
\{ (\omega_{n,-n})_{n \in {\Bbb Z}} \in \Sigma_\kappa^{\Bbb Z} 
\mid 
(\omega_{i,j})_{(i,j)\in {\Bbb Z}^2}\in 
X_{\rho,\eta}^\kappa \}.
\end{equation*}
By Lemma 2.1,
an element 
$(\omega_{n,-n})_{n \in {\Bbb Z}} $ 
of 
$X_{\delta^\kappa}$
may 
be extend to
$
(\omega_{i,j})_{(i,j)\in {\Bbb Z}^2}\in 
X_{\rho,\eta}^\kappa
$
in a unique way.
Hence the one-dimensional subshift 
$X_{\delta^\kappa}$
determines the two-dimensional subshift
$X_{\rho,\eta}^\kappa$.
Therefore we have 
\begin{lem}
For $\CTDS$,
the two-dimensional subshift
$X_{\rho,\eta}^\kappa$ 
is not empty if and only if 
the one-dimensional subshift
$X_{\delta^\kappa}$ 
is not empty. 
\end{lem}
For $\CTDS$, 
we will have a $C^*$-symbolic dynamical system
$(\A, \delta^\kappa,\Sigma_\kappa)$
in Section 5.
It presents the subshift
$X_{\delta^\kappa}$.
Since a subshift presented by a  $C^*$-symbolic dynamical system
is always not empty,
one sees that 
$X_{\rho,\eta}^\kappa$ 
is not empty.
\begin{prop}
For $\CTDS$, 
the two-dimensional subshift
$X_{\rho,\eta}^\kappa$ 
is not empty.
\end{prop}


\medskip

The $C^*$-algbebra $\ORE$ 
is defined to be the universal $C^*$-algebra 
$C^*(x, S_\alpha, T_a ; x \in \A, \alpha \in \Sigma^\rho, a \in \Sigma^\eta)$  
generated by 
$x \in \A$ 
and partial isometries $S_{\alpha}, \alpha \in \Sigma^\rho$,
$T_a, a \in \Sigma^\eta$
subject to the following relations called $(\rho,\eta)$:
\begin{align}
\sum_{\beta \in \Sigma^\rho}S_{\beta}S_{\beta}^*  =1,\qquad
x S_\alpha S_\alpha^* & =  S_\alpha S_\alpha^*  x,\qquad
S_\alpha^* x S_\alpha  = \rho_\alpha(x),\\
\sum_{b \in \Sigma^\eta} T_{b} T_{b}^*  =1,\qquad
x T_a T_a^*  & =  T_a T_a^*  x,\qquad
T_a^* x T_a  = \eta_a(x),\\  
S_\alpha T_b  = T_a S_\beta &
\qquad
\text{ if } 
\quad
\kappa(\alpha, b) = (a,\beta)
\end{align}
for all $x \in {\cal A}$ and $\alpha \in \Sigma^\rho, a \in \Sigma^\eta.$
We will study the algebra $\ORE$.
If 
$ \kappa(\alpha, b) = (a,\beta)$,
we write
as
$ 
(\alpha, b) \overset{\kappa}{\cong} (a,\beta).
$
\begin{lem}
For $\alpha \in \Sigma^\rho, a \in \Sigma^\eta$,
one has 
$T_a^* S_\alpha \ne 0$ 
if and only if
there exist 
$ b \in \Sigma^\eta, \beta \in \Sigma^\rho$
such that
$(\alpha, b) \overset{\kappa}{\cong} (a,\beta)$.
 \end{lem}
\begin{pf}
Suppose that 
$T_a^* S_\alpha \ne 0$.
As
$T_a^* S_\alpha = \sum_{b' \in \Sigma^\eta}T_a^* S_\alpha T_{b'}T_{b'}^*$,
there exists
$b' \in \Sigma^\eta$
such that
$T_a^* S_\alpha T_{b'} \ne 0$.
Hence
$\eta_{b'} \circ \rho_\alpha \ne 0$ 
so that
$(\alpha, b') \in \Sigma^{\rho\eta}$.
Then 
one may find 
$(a', \beta') \in \Sigma^\rho$
such that
$(\alpha, b') \overset{\kappa}{\cong} (a',\beta')$
and hence
$S_\alpha T_{b'} = T_{a'}S_{\beta'}$.
Since
$
0 \ne T_a^* S_\alpha T_{b'} = T_a^* T_{a'}S_{\beta'}$,
one sees that 
$a = a'$.
Putting
$b = b', \beta = \beta'$,
we have  
$ \kappa(\alpha, b) = (a,\beta)$.

Suppose next that
$ \kappa(\alpha, b) = (a,\beta)$.
Since
$\eta_b \circ \rho_\alpha = \rho_\beta \circ \eta_a \ne 0$,
one has
$0 \ne S_\alpha T_b = T_a S_\beta$.
It follows that
 $
 S_\beta^* T_a^* S_\alpha T_b = 
(T_a S_\beta)^* T_a S_\beta
$
so that
$T_a^* S_\alpha \ne 0$.
\end{pf}
\begin{lem}
For $\alpha \in \Sigma^\rho, a \in \Sigma^\eta$,
we have
\begin{equation}  
T_a^* S_\alpha 
=
\sum_{
\stackrel{b,\beta}{\kappa(\alpha, b) = (a,\beta)}}
 S_\beta \eta_b(\rho_\alpha(1))T_b^*
\end{equation} 
and hence 
\begin{equation}  
S_\alpha^* T_a 
=
\sum_{
\stackrel{b,\beta}{\kappa(\alpha, b) = (a,\beta)}}
 T_b \rho_\beta(\eta_a(1))S_\beta^*.
\end{equation} 
\end{lem}
\begin{pf}
We may assume that
$T_a^* S_\alpha \ne 0$. 
One has
$T_a^* S_\alpha = \sum_{b' \in \Sigma^\eta}T_a^* S_\alpha T_{b'}T_{b'}^*$.
For $b'\in \Sigma^\eta$
with $(\alpha,b') \in \Sigma^{\rho,\eta}$, 
and for 
$\beta' \in \Sigma^\rho$
such that
$ \kappa(\alpha, b') = (a',\beta')$
for some $a' \in \Sigma^\eta$,
one has 
\begin{equation*}
T_a^* S_\alpha T_{b'}T_{b'}^*
 =  T_a^* T_{a'} S_{\beta'} T_{b'}^*.
\end{equation*}
Hence
$T_a^* S_\alpha T_{b'}T_{b'}^*\ne 0$
implies 
$a= a'$.
Since $T_a^*T_a = \eta_a(1)$ which commutes with $S_{\beta'}S_{\beta'}^*$,
we have
$$
 T_a^* T_a S_{\beta'} T_{b'}^*
 = S_{\beta'}S_{\beta'}^* T_a^*T_a S_{\beta'} T_{b'}^*
 = S_{\beta'} \rho_{\beta'}(\eta_a(1)) T_{b'}^*
 = S_{\beta'} \eta_{b'}(\rho_\alpha(1)) T_{b'}^*.
$$
It follows that
\begin{equation*}
T_a^* S_\alpha 
=\sum_{
\stackrel{b',\beta'}{\kappa(\alpha, b') = (a,\beta')}}
T_a^* T_a S_{\beta'} T_{b'}^*
=\sum_{
\stackrel{b',\beta'}{\kappa(\alpha, b') = (a,\beta')}}
S_{\beta'} \eta_{b'}(\rho_\alpha(1)) T_{b'}^*.
\end{equation*}
\end{pf}
Hence we have
\begin{lem}
For $\alpha \in \Sigma^\rho, a \in \Sigma^\eta$,
we have
\begin{equation}  
T_a T_a^* S_\alpha S_\alpha^* 
=
\sum_{
\stackrel{b}{\kappa(\alpha, b) = (a,\beta)\text{for some } \beta}}
 S_\alpha T_b T_b^* S_\alpha^*.
\end{equation} 
Hence 
$T_a T_a^* $ 
commutes with
$S_\alpha S_\alpha^*$.
\end{lem}
\begin{pf}
By the preceding lemma,
we have
\begin{align*}  
T_a T_a^* S_\alpha S_\alpha^* 
& =
\sum_{
\stackrel{b,\beta}{\kappa(\alpha, b) = (a,\beta)}}
T_a S_{\beta} \eta_{b}(\rho_\alpha(1)) T_{b}^* S_\alpha^* \\
& =
\sum_{
\stackrel{b,\beta}{\kappa(\alpha, b) = (a,\beta)}}
S_\alpha T_b \eta_{b}(\rho_\alpha(1)) T_{b}^* S_\alpha^* \\
& =
\sum_{
\stackrel{b,\beta}{\kappa(\alpha, b) = (a,\beta)}}
S_\alpha  \rho_\alpha(1)) T_b  T_{b}^* S_\alpha^* \\
 & =
\sum_{
\stackrel{b}{\kappa(\alpha, b) = (a,\beta)\text{for some } \beta}}
 S_\alpha T_b T_b^* S_\alpha^*.
\end{align*} 
\end{pf}
More generally, we have
\begin{lem} Suppose that $\A$ is commutative.
For $\alpha \in \Sigma^\rho, a \in \Sigma^\eta$,
and
$x, y \in \A$, 
we know that  
$T_a y T_a^* $ 
commutes with
$S_\alpha x S_\alpha^*$.
\end{lem}
\begin{pf}
It follows that
\begin{align*}  
T_a y T_a^* S_\alpha x S_\alpha^* 
& = T_a y 
\sum_{
\stackrel{b,\beta}{\kappa(\alpha, b) = (a,\beta)}}
S_{\beta} \eta_{b}(\rho_\alpha(1)) T_{b}^* x  S_\alpha^* \\
& = 
\sum_{
\stackrel{b,\beta}{\kappa(\alpha, b) = (a,\beta)}}
T_a S_{\beta}S_{\beta}^* y S_{\beta}\eta_{b}(\rho_\alpha(1)) T_{b}^* xT_{b}T_{b}^* S_\alpha^* \\
& = 
\sum_{
\stackrel{b,\beta}{\kappa(\alpha, b) = (a,\beta)}}
S_\alpha T_b \rho_{\beta}(y)\eta_{b}(\rho_\alpha(1)) \eta_b(x) S_\beta^* T_a^* \\
& = 
\sum_{
\stackrel{b,\beta}{\kappa(\alpha, b) = (a,\beta)}}
S_\alpha T_b \eta_b(x) \eta_{b}(\rho_\alpha(1))\rho_{\beta}(y) S_\beta^* T_a^* \\
& = 
\sum_{
\stackrel{b,\beta}{\kappa(\alpha, b) = (a,\beta)}}
S_\alpha x \rho_\alpha(1) T_b S_\beta^* y T_a^* \\
& = 
\sum_{
\stackrel{b,\beta}{\kappa(\alpha, b) = (a,\beta)}}
S_\alpha x S_\alpha^* S_\alpha T_b S_\beta^* T_a^* T_a y T_a^* \\
& = 
\sum_{
\stackrel{b,\beta}{\kappa(\alpha, b) = (a,\beta)}}
S_\alpha x ( S_\alpha^* S_\alpha T_b T_b^* S_\alpha^* T_a )y T_a^*. 
\end{align*}
Now 
if $(\alpha, b') \not\in \Sigma^{\rho,\eta}$,
then 
$S_\alpha T_{b'} = 0$.
Hence
$$
\sum_{
\stackrel{b,\beta}{\kappa(\alpha, b) = (a,\beta)}}
 S_\alpha^* S_\alpha T_b T_b^* S_\alpha^* T_a 
=
\sum_{b}
 S_\alpha^* S_\alpha T_b T_b^* S_\alpha^* T_a 
= S_\alpha^* T_a. 
$$
Therefore we have
\begin{equation*}  
T_a y T_a^* S_\alpha x S_\alpha^*
=S_\alpha x S_\alpha^* T_a y T_a^*.
\end{equation*}  
\end{pf}
We set 
\begin{align*}
 {\cal D}_{\rho,\eta} 
= &  C^*( S_{\mu} T_\zeta x T_\zeta^* S_{\mu}^{*}:
 \mu \in B_*(\Lambda_\rho), \zeta \in B_*(\Lambda_\eta), x \in \A), \\
 {\cal D}_{j,k} 
= &  C^*( S_{\mu} T_\zeta x T_\zeta^* S_{\mu}^{*}: 
\mu \in B_j(\Lambda_\rho), \zeta \in B_k(\Lambda_\eta), x \in \A) 
          \text{ for } j, k \in \Zp. 
\end{align*}
By the commutation relation (3.5),
one sees that
\begin{equation*}
 {\cal D}_{j,k} 
=   C^*(  T_\xi S_{\nu} x  S_{\nu}^{*} T_\xi^* : 
\nu \in B_j(\Lambda_\rho),  \xi \in B_k(\Lambda_\eta), x \in \A) 
\end{equation*}
The identities 
\begin{align*}
S_{\mu} T_\zeta x T_\zeta^* S_{\mu}^{*}
& = \sum_{a \in \Sigma^\eta}
S_{\mu} T_{\zeta a} \eta_a(x) T_{\zeta a}^* S_{\mu}^{*},\\
 T_\xi S_{\nu} x  S_{\nu}^{*} T_\xi^*
& = \sum_{\alpha \in \Sigma^\rho} 
T_\xi S_{\nu \alpha} \rho_\alpha(x) S_{\nu \alpha}^*T_\xi^*
\end{align*}
for $x \in \A$ and 
$\mu, \nu \in B_j(\Lambda_\rho), \zeta, \xi \in B_k(\Lambda_\eta)$
yield the embeddings
\begin{equation*}
 {\cal D}_{j,k} 
\hookrightarrow  {\cal D}_{j,k+1}, 
\qquad 
 {\cal D}_{j,k} 
\hookrightarrow  {\cal D}_{j+1,k}
\end{equation*}
respectively such that
$
\cup_{j,k\in \Zp}{\cal D}_{j,k}
$ is dense in ${\cal D}_{\rho,\eta}$.
\begin{prop}
If $\A$ is commutative, so is ${\cal D}_{\rho,\eta}$.
\end{prop}
\begin{pf}
The preceding lemma tells us that
${\cal D}_{1,1}$ is commutative.
Suppose that
the algebra
${\cal D}_{j,k}$
is commutative
for a fixed $j,k \in {\Bbb N}$.
We will show that
the both algebras 
${\cal D}_{j+1,k}$
and
${\cal D}_{j,k+1}$
are commutative.
For the algebra
${\cal D}_{j+1,k}$,
it consists of linear span of elements of the form:
$$
S_\alpha x S_\alpha^*\qquad \text{ for } 
x \in {\cal D}_{j,k}, \alpha \in \Sigma^\rho.
$$ 
Let
$ x,y\in {\cal D}_{j,k}, \alpha, \beta \in \Sigma^\rho.
$ 
We will show that
$S_\alpha x S_\alpha^*$ 
commutes with
both 
$
S_\beta y S_\beta^*
$
and
$y$.
If $\alpha = \beta$,
it is easy to see that
$
S_\alpha x S_\alpha^*
$
commutes with
$
S_\alpha y S_\alpha^*,
$
because 
$\rho_\alpha(1) \in \A \subset {\cal D}_{j,k}$.
If 
$\alpha \ne \beta$,
both
$
S_\alpha x S_\alpha^* 
S_\beta y S_\beta^* 
$
and
$
S_\beta y S_\beta^* S_\alpha x S_\alpha^*
$
are zeros.
Since
$
S_\alpha^* y S_\alpha 
\in {\cal D}_{j-1,k} \subset {\cal D}_{j,k},
$
one sees
$
S_\alpha^* y S_\alpha 
$ 
commutes with
$x$.
One also sees that
$S_\alpha S_\alpha^* \in {\cal D}_{j,k}$
commutes with $y$.
It follows that
\begin{equation*}
 S_\alpha x S_\alpha^* y 
=S_\alpha x S_\alpha^* y S_\alpha  S_\alpha^*   
=S_\alpha  S_\alpha^* y S_\alpha x S_\alpha^* 
=y S_\alpha x S_\alpha^*.
\end{equation*}
Hence 
the algebra
${\cal D}_{j+1,k}$
is commutative,
and similarly so is
${\cal D}_{j,k+1}$.
By induction,
one knows that
the algebras 
${\cal D}_{j,k}$
are all commutative for all
$j,k \in {\Bbb N}$.
Since
$\cup_{j,k\in{\Bbb N}}{\cal D}_{j,k}$
is dense in 
${\cal D}_{\rho,\eta}$,
${\cal D}_{\rho,\eta}$
is commuatative.
\end{pf}

\begin{prop}
Let
${\cal O}_{\rho,\eta}^{alg}$
be the dense $*$-subalgebra algebraically generated by elements 
$x \in \A$, $S_\alpha, \alpha \in \Sigma^\rho$
and
$T_a, a \in \Sigma^\eta$.
Then each element of 
${\cal O}_{\rho,\eta}^{alg}$ 
is a finite linear combination of elements of the form:
\begin{equation}
S_\mu T_\zeta x T_\xi^* S_\nu^*
\qquad
\text{ for }
x \in \A, \mu,\nu \in B_*(\Lambda_\rho),\zeta,\xi \in B_*(\Lambda_\eta)
\end{equation}
where
$
S_\mu = S_{\mu_1}\cdots S_{\mu_k},
S_\nu = S_{\nu_1}\cdots S_{\nu_n}
$ 
for 
$\mu = \mu_1\cdots \mu_k, \nu = \nu_1\cdots \nu_n$
and
$
T_\zeta = T_{\zeta_1}\cdots T_{\zeta_l},
T_\xi = T_{\xi_1}\cdots T_{\xi_m}$ 
for 
$\zeta = \zeta_1\cdots \zeta_l,
\xi = \xi_1\cdots \xi_m$.
 \end{prop}
\begin{pf}
For $\alpha, \beta\in \Sigma^\rho$, $a,b \in \Sigma^\eta$ and $x \in \A$,
we have
\begin{align*}
S_\alpha^* S_\beta 
& = 
{\begin{cases}
\rho_\alpha(1) \in \A & \text{ if } \alpha = \beta,\\
0 & \text{otherwise},
\end{cases}} \\
S_\alpha^* T_a 
& = 
\sum_{
\stackrel{b,\beta}{\kappa(\alpha, b) = (a,\beta)}}
 T_b \rho_\beta(\eta_a(1))S_\beta^*, \\
S_\alpha^* x 
& = \rho_\alpha(x) S_\alpha,\qquad \quad 
S_\beta^* T_a^* 
 = T_b^* S_\alpha^*.
\end{align*}
And also 
\begin{align*}
T_a^* T_b 
& = 
{\begin{cases}
\eta_a(1) \in \A & \text{ if } a = b,\\
0 & \text{otherwise},
\end{cases}} \\
 T_a^*S_\alpha 
& = 
\sum_{
\stackrel{b,\beta}{\kappa(\alpha, b) = (a,\beta)}}
S_\beta \eta_b (\rho_\alpha(1))  T_b^*, \\
T_a^* x 
& = \eta_a(x) T_a^*.
\end{align*}
Therefore we conclude that
any element of 
${\cal O}_{\rho,\eta}^{alg}$ 
is a finite linear combination of elements of the form:
$
S_\mu T_\zeta x T_\xi^* S_\nu^*.
$
\end{pf}
Similarly we have
\begin{prop}
Each element of 
${\cal O}_{\rho,\eta}^{alg}$ 
is a finite linear combination of elements of the form:
\begin{equation}
T_\zeta S_\mu  x S_\nu^* T_\xi^* 
\qquad
\text{ for }
x \in \A, \mu,\nu \in B_*(\Lambda_\rho),\zeta,\xi \in B_*(\Lambda_\eta).
\end{equation}
\end{prop}

In the rest of this section,
we will have a $C^*$-symbolic dynamical system
$(\A,\delta^\kappa,\Sigma_\kappa)$
from
$\CTDS$,
which presents the one-dimensional subshift
$X_{\delta^\kappa}$
described in the previous section.
For $\CTDS$, define an endomorphism
$\delta_\omega^\kappa$ on $\A$ for
$\omega \in \Sigma_\kappa$ 
by setting
\begin{equation*}
\delta_\omega^\kappa(x)
=\eta_b(\rho_\alpha(x)) (=\rho_\beta(\eta_a(x))), \qquad x \in \A,
\quad
\omega=(\alpha,b,a,\beta) \in \Sigma_\kappa.
\end{equation*}
\begin{lem}
$(\A,\delta^\kappa,\Sigma_\kappa)$
is a $C^*$-symbolic dynamical system that presents $X_{\delta^\kappa}$.
\end{lem}
\begin{pf}
We will show that 
$\delta^\kappa$ is essential and faithful.
Now both $C^*$-symbolic dynamical systems 
$(\A, \eta,\Sigma^\eta)$ and $(\A, \rho,\Sigma^\eta)$
are essential.
We are further assuming that 
both $C^*$-symbolic dynamical systems 
$(\A, \eta,\Sigma^\eta)$ and $(\A, \rho,\Sigma^\eta)$
are central.
Hence it is clear that
$\delta_\omega^\kappa(Z_\A) \subset Z_\A$.
By the inequalities
\begin{equation*}
\sum_{\omega \in \Sigma_\kappa}
\delta_\omega^\kappa(1) = 
\sum_{b \in \Sigma^\eta} \sum_{\alpha \in \Sigma^\rho}\eta_b(\rho_\alpha(1)) 
\ge
\sum_{b \in \Sigma^\eta} \eta_b(1) \ge 1 
\end{equation*}
$\{ \delta^\kappa\}_{\omega\in \Sigma_\kappa}$
is essential.
For any nonzero $x \in \A$, 
there exists $\alpha \in \Sigma^\rho$
such that $\rho_\alpha(x) \ne 0 $ and  there exists $b \in \Sigma^\eta$
such that $\eta_b(\rho_\alpha(x)) \ne 0$. 
This means that
$\delta_\omega^\kappa(x) \ne 0$ for
$\omega=(\alpha,b,a,\beta) \in \Sigma_\kappa$. 
Hence $\delta^\kappa$ is faithful
so that
$(\A,\delta^\kappa,\Sigma_\kappa)$
is a $C^*$-symbolic dynamical system.
It is obvious  that 
the presented subshift
by
$(\A,\delta^\kappa,\Sigma_\kappa)$
is $X_{\delta^\kappa}$.
\end{pf}
Put
\begin{align*}
\widehat{X}_{\rho,\eta}^\kappa
& =\{ (\omega_{i,-j})_{(i,j)\in {\Bbb N}^2} \in \Sigma_\kappa^{{\Bbb N}^2} \mid
(\omega_{i,j})_{(i,j)\in {\Bbb Z}^2} \in X_{\rho,\eta}^\kappa \} \\
\intertext{and}
\widehat{X}_{\delta^\kappa}
& =\{ (\omega_{n,-n})_{n\in {\Bbb N}} \in \Sigma_\kappa^{\Bbb N} \mid
(\omega_{i,j})_{(i,j)\in {\Bbb N}^2} \in \widehat{X}_{\rho,\eta}^\kappa \}. 
\end{align*}
The latter set
$\widehat{X}_{\delta^\kappa}$
is the right one-sided subshift for
$X_{\delta^\kappa}$.

\begin{lem}
A configuration 
$
(\omega_{i,-j})_{(i,j)\in {\Bbb N}^2} \in
\widehat{X}_{\rho,\eta}^\kappa
$
can extend to a whole configuration
$
(\omega_{i,j})_{(i,j)\in {\Bbb Z}^2} \in X_{\rho,\eta}^\kappa.
$
\end{lem}
\begin{pf}
For
$
(\omega_{i,-j})_{(i,j)\in {\Bbb N}^2} \in
\widehat{X}_{\rho,\eta}^\kappa,
$
put
$x_i = \omega_{i,-i}, i \in {\Bbb N}$
so that
$x = (x_i)_{i \in {\Bbb N}}
\in \widehat{X}_{\delta^\kappa}.
$ 
Since
$\widehat{X}_{\delta^\kappa}$
is a one-sided subshift,
there exists an extension
$\tilde{x} \in X_{\delta^\kappa}$
to two-sided sequence such that
$\tilde{x}_{[1,\infty)} = x$.
By the diagonal property,
$\tilde{x}$ 
determines a whole configuration
$\tilde{\omega}$ 
to
${\Bbb Z}^2$
such that
$\tilde{\omega} \in X_{\delta,\eta}^\kappa$
and
$(\tilde{\omega}_{i,-i})_{i \in {\Bbb N}} = \tilde{x}$.
Hence
$\tilde{\omega}_{i,-j} = \omega_{i,-j}$
for all $i.j \in {\Bbb N}$.
\end{pf}

Let 
${\frak D}_{\rho,\eta}$
be the $C^*$-subalgebra of $\DRE$
defined by
\begin{align*}
{\frak D}_{\rho,\eta} 
& =   C^*( S_{\mu} T_\zeta  T_\zeta^* S_{\mu}^{*}:
 \mu \in B_*(\Lambda_\rho), \zeta \in B_*(\Lambda_\eta)\\
& =   C^*( T_\xi S_{\nu} S_\nu^* T_\xi^*:
 \nu \in B_*(\Lambda_\rho), \xi \in B_*(\Lambda_\eta)
\end{align*}
which is a commutative $C^*$-subalgebra of $\DRE$.
Put
for $\mu=\mu_1\cdots\mu_n \in B_*(\Lambda_\rho)$,
$\zeta=\zeta_1\cdots\zeta_m \in B_*(\Lambda_\eta)$
the cylinder set
$$
U_{\mu,\zeta}
=\{ (\omega_{i,-j})_{(i,j) \in {\Bbb N}^2} \in \widehat{X}_{\rho,\eta}^\kappa \mid
t(\omega_{i,-1}) =\mu_i, i=1,\cdots,n,
r(\omega_{n,-j}) =\zeta_j, j=1,\cdots,m
\}
$$
The following lemma is direct.
\begin{lem}
${\frak D}_{\rho,\eta}$ is isomorphic to
$C(\widehat{X}_{\rho,\eta}^\kappa)$
through the corespondence
such that 
 $S_{\mu} T_\zeta  T_\zeta^* S_{\mu}^{*}$ sends to
$
\chi_{U_{\mu,\zeta}},
$
where
$
\chi_{U_{\mu,\zeta}}
$
is the characteristic function for the cylinder set 
$U_{\mu,\zeta}$
on
$\widehat{X}_{\rho,\eta}^\kappa$.
\end{lem}

\section{Condition (I) for  $C^*$-textile  dynamical systems}

The notion of condition (I) for finite square matrices 
with entries in $\{0,1\}$ has been introduced in \cite{CK}. 
The condition has been generalized by many authors to 
corresponding conditions for generalizations of the Cuntz-Krieger algebras,
 for instance,  infinite directed graphs (\cite{KPRR}), 
infinite matrices with entries in $\{0,1\}$ (\cite{EL}), 
Hilbert $C^*$-bimodules (\cite{KPW},see also \cite{Re}, etc.).
 The condition (I) for 
$C^*$-symbolic dynamical systems
(including $\lambda$-graph systems) has been also defined in 
\cite{MaContem}(cf. \cite{Ma7}, \cite{Ma11}).
All of these conditions 
give rise to the uniqueness of the associated $C^*$-algebras 
subject to some operator relations of the generating elements.

 In this section, 
 we will introduce the notion of condition (I) 
 for $C^*$-textile dynamical systems
 to prove 
 the uniqueness of the $C^*$-algebras
$\ORE$
under the relation $(\rho,\eta;\kappa)$.
In what follows, for a subset $F$ of a $C^*$-algebra $\B$,
we will denote by $C^*(F)$ the $C^*$-subalgebra of $\B$ generated by $F$.

Let $\CTDS$ 
be a $C^*$-symbolic dynamical system over $\Sigma$
and $X_{\rho,\eta}^\kappa$ the associated two-dimensional subshift.
Denote by
$\Lambda_\rho, \Lambda_\eta$
the associated subshifts to
the
 $C^*$-symbolic dynamical systems
$(\A, \rho, \Sigma^\rho), (\A, \eta, \Sigma^\eta) 
$
 respectively.
For 
$\mu =(\mu_1,\dots,\mu_j) \in B_j(\Lambda_\rho),
\zeta =(\zeta_1,\dots,\zeta_k) \in B_k(\Lambda_\eta)$, 
we put
$S_\mu = S_{\mu_1}\cdots S_{\mu_j},
T_\zeta = T_{\zeta_1}\cdots T_{\zeta_k}$
and
$\rho_\mu = \rho_{\mu_j}\circ \cdots \circ \rho_{\mu_1},
\eta_\mu = \eta_{\zeta_k}\circ \cdots \circ \eta_{\zeta_1}$
respectively.
We denote by
$|\mu|, |\zeta|$
the lengths $j,k$ respectively.

In the algebra $\ORE$, 
we set the subalgebras
\begin{align*}
 {\cal F}_{\rho,\eta} 
= &  C^*( S_{\mu} T_\zeta x T_\xi^* S_{\nu}^{*}:
 \mu, \nu \in B_*(\Lambda_\rho), \zeta, \xi \in B_*(\Lambda_\eta), 
|\mu| = |\nu|,|\zeta| = |\xi|, x \in \A) \\
\intertext{and for $ j, k \in \Zp$ } 
{\cal F}_{j,k} 
= &  C^*( S_{\mu} T_\zeta x T_\xi^* S_{\nu}^{*}: 
\mu, \nu \in B_j(\Lambda_\rho), \zeta, \xi \in B_k(\Lambda_\eta), x \in \A). 
\end{align*}
We notice that
\begin{equation*}
 {\cal F}_{j,k} 
=   C^*(  T_\zeta S_{\mu} x  S_{\nu}^{*} T_\xi^* : 
\mu, \nu \in B_j(\Lambda_\rho), \zeta, \xi \in B_k(\Lambda_\eta), x \in \A) 
\end{equation*}

The identities 
\begin{align}
S_{\mu} T_\zeta x T_\xi^* S_{\nu}^{*}
& = \sum_{a \in \Sigma^\eta}
S_{\mu} T_{\zeta a} \eta_a(x) T_{\xi a}^* S_{\nu}^{*},\\
 T_\zeta S_{\mu} x  S_{\nu}^{*} T_\xi^*
& = \sum_{\alpha \in \Sigma^\rho} 
T_\zeta S_{\mu \alpha} \rho_\alpha(x) S_{\nu \alpha}^*T_\xi^*
\end{align}
for $x \in \A$ and 
$\mu, \nu \in B_j(\Lambda_\rho), \zeta, \xi \in B_k(\Lambda_\eta)$
yield the embeddings
\begin{equation*}
 {\cal F}_{j,k} 
\hookrightarrow  {\cal F}_{j,k+1}, 
\qquad 
 {\cal F}_{j,k} 
\hookrightarrow  {\cal F}_{j+1,k}
\end{equation*}
such that
$
\cup_{j,k\in \Zp}{\cal F}_{j,k}
$ is dense in ${\cal F}_{\rho,\eta}$.

By the universality
of $\ORE$,
we may define
an action
$
\widehat{\kappa} : {\Bbb T}^2 \longrightarrow \Aut(\ORE)
$
of the $2$-dimensiona torus group
${\Bbb T}^2 = \{ (z,w) \in {\Bbb C}^2 \mid |z| = |w| =1 \}$
to
$\ORE$
by setting
\begin{equation*}
\widehat{\kappa}_{z,w}(S_\alpha)  = z S_\alpha, \quad
\widehat{\kappa}_{z,w}(T_a)  = w T_a, \quad
\widehat{\kappa}_{z,w}(x) = x
\end{equation*}
for $\alpha \in \Sigma^\rho, a \in \Sigma^\eta$,
$ x\in \A$
and
$z,w \in {\Bbb T}$.
We call the action
$
\widehat{\kappa} : {\Bbb T}^2 \longrightarrow \Aut(\ORE)
$
the gauge action of ${\Bbb T}^2$ on $\ORE$.
The fixed point algebra of $\ORE$ under $\widehat{\kappa}$ is denoted by
$(\ORE)^{\widehat{\kappa}}$.
Let $\E_{\rho,\eta}: \ORE \longrightarrow (\ORE)^{\widehat{\kappa}}$ 
be the conditional expectaton defined by
$$
\E_{\rho,\eta}(X) 
= \int_{(z, w) \in {\Bbb T}^2} \widehat{\kappa}_(z,w)(X)\, dzdw, \qquad X \in \ORE.
$$
The following lemma is routine.
\begin{lem}
$(\ORE)^{\widehat{\kappa}} = {\cal F}_{\rho,\eta}$.
\end{lem}
Put
$ 
\phi_{\rho},  \phi_{\eta} : 
{\cal D}_{\rho,\eta}\longrightarrow {\cal D}_{\rho,\eta}
$
by setting
\begin{equation*}
\phi_\rho(X)  = \sum_{\alpha \in \Sigma^\rho} S_\alpha X S_\alpha^*,\qquad
\phi_\eta(X)  = \sum_{a \in \Sigma^\eta} T_a X T_a^*,\qquad
X \in {\cal D}_{\rho,\eta}.  
\end{equation*}
It is easy to see 
$$
\phi_\rho \circ \phi_\eta = \phi_\eta \circ \phi_\rho
\quad \text{ on } {\cal D}_{\rho,\eta}.
$$

\noindent
{\bf Definition. }
A $C^*$-textile dynamical system $\CTDS$
satisfies {\it condition }(I) if
there exists a unital increasing sequence
$$
\A_0 \subset \A_1 \subset \cdots \subset \A
$$
of $C^*$-subalgebras of $\A$ such that 

(1)
$\rho_\alpha(\A_l) \subset \A_{l+1}$, 
$\eta_a(\A_l) \subset \A_{l+1}$ 
for all $l \in \Zp, \alpha \in \Sigma^\rho, a \in \Sigma^\eta$,
  
(2)
 $\cup_{l\in \Zp} \A_l $ is dense in $\A$,

(3) 
for $\epsilon > 0$, $j,k,l\in {\Bbb N}$ with $j+k \le l$
and
$ X_0 \in {\cal F}_{j,k}^{l}
= C^*( S_{\mu} T_\zeta x T_\xi^* S_{\nu}^{*}: 
\mu, \nu \in B_j(\Lambda_\rho), 
\zeta,\xi \in B_k(\Lambda_\eta),
x \in \A_l)$,
there exists an element 
$g \in {\cal D}_{\rho,\eta} \cap {\A_l}^\prime 
(= \{ y \in {\cal D}_{\rho,\eta} \mid y a = a y \ \text{ for } a \in \A_l \})
$
with
$0 \le g \le 1$
such that
\begin{enumerate}
\renewcommand{\labelenumi}{(\roman{enumi})}
\item
 $\| X_0 \phi_{\rho}^j \circ \phi_{\eta}^k(g) \| \ge \| X_0 \| - \epsilon$,
\item
 $
 \phi_{\rho}^n(g) \phi_{\eta}^m(g) 
=\phi_{\rho}^n(( \phi_{\eta}^m(g)))g 
=\phi_{\rho}^n(g) g 
=\phi_{\eta}^m(g) g = 0
$ 
for all $n=1,2,\dots,j$, $m= 1,2,\dots, k.$
\end{enumerate} 

If in particular,
one may take the above subalgebras
$\A_l \subset \A$, $l=0,1,2, \dots $
to be of finite dimensional,
then 
$\CTDS$ is said to 
satisfy {\it AF-condition }(I).
In this case,
$\A = \overline{\cup_{l=0}^\infty \A_l}$ is an AF-algebra.

As the element $g$ above  belongs to 
the diagonal subalgebra ${\cal D}_{\rho,\eta}$ 
of ${\cal F}_{\rho,\eta}$,
 the condition (I) of $\CTDS$
 is intrinsically determined by itself
 by virtue of Lemma 4.3 below.

 

We will also introduce the following  condition called {\it free}, 
which will be stronger than condition (I) but easier to confirm than condition (I).

\noindent
{\bf Definition. }
A $C^*$-textile dynamical system $\CTDS$
is said to be {\it free } if
there exists a unital increasing sequence
$$
\A_0 \subset \A_1 \subset \cdots \subset \A
$$
of $C^*$-subalgebras of $\A$ such that 

(1)
$\rho_\alpha(\A_l) \subset \A_{l+1}$, 
$\eta_a(\A_l) \subset \A_{l+1}$ 
for all $l \in \Zp, \alpha \in \Sigma^\rho, a \in \Sigma^\eta$,
  
(2)
 $\cup_{l\in \Zp} \A_l $ is dense in $\A$,

(3) 
for  $j,k,l\in {\Bbb N}$ with $j+k \le l$
there exists a projection 
$q \in {\cal D}_{\rho,\eta} \cap {\A_l}^\prime 
$
such that
\begin{enumerate}
\renewcommand{\labelenumi}{(\roman{enumi})}
\item
 $ q a \ne 0$ for $0\ne a \in \A_l$,
\item
 $
 \phi_{\rho}^n(q) \phi_{\eta}^m(q) 
=\phi_{\rho}^n(( \phi_{\eta}^m(q)))q 
=\phi_{\rho}^n(q) q 
=\phi_{\eta}^m(q) q = 0
$ 
for all $n=1,2,\dots,j$, $m= 1,2,\dots, k.$
\end{enumerate} 
If in particular,
one may take the above subalgebras
$\A_l \subset \A$, $l=0,1,2, \dots $
to be of finite dimensional,
then 
$\CTDS$ is said to 
be {\it AF-free }.

\begin{prop}
If a $C^*$-textile dynamical system $\CTDS$
is free (resp. AF-free), 
then it satisfies condition (I) (resp. AF-condition (I)). 
\end{prop}
\begin{pf}
Assume that
$\CTDS$
is free.
Take 
an increasing sequence 
$\A_l, l \in {\Bbb N}$
of $C^*$-subalgebras of $\A$
satisfying the above conditions (1),(2),(3) of freeness.
For  
 $j,k,l\in {\Bbb N}$ with $j+k \le l$
there exists a projection 
$q \in {\cal D}_{\rho,\eta} \cap {\A_l}^\prime 
$
satisfying the above two conditions (i) and (ii) of (3).
Put
$Q_{j,k}^l =\phi_{\rho}^j(\phi_{\eta}^k(q))$.
For
$x \in \A_l, \mu,\nu \in B_j(\Lambda_\rho), \xi,\zeta \in B_k(\Lambda_\eta), $
one has the equality
$$
Q_{j,k}^l S_\mu T_\zeta x T_\xi^* S_\nu^* = S_\mu T_\zeta x T_\xi^* S_\nu^*
$$
so that
$Q_{j,k}^l$
commutes with all of elements of 
${\cal F}_{j,k}^l$.
By using the condition (i) of (3) for $q$
one directly sees that
$
S_\mu T_\zeta x T_\xi^* S_\nu^*  \ne 0
$
if and only if
$
Q_{j,k}^l S_\mu T_\zeta x T_\xi^* S_\nu^*\ne 0.
$
Hence the map 
\begin{equation*}
X \in {\cal F}_{j,k}^l \longrightarrow 
X Q_{j,k}^l \in {\cal F}_{j,k}^l Q_{j,k}^l
\end{equation*}
defines a homomorphism, that is proved to be injective 
by a similar proof to the proof of 
\cite[Proposition 3.7]{MaMZ2010}.
Hence we have
$
\| X Q_{j,k}^l \| =\| X \| 
 \ge \| X \| - \epsilon 
$
for all $X \in {\cal F}_{j,k}^l$. 
\end{pf}

\medskip

Let $\B$ be a unital $C^*$-algebra.
Suppose that there exist an injective $*$-homomorphism
$\pi: \A \longrightarrow \B$ preserving their units
and 
two families $s_\alpha \in \B, \alpha \in \Sigma^\rho$
and
$t_a \in \B, a \in \Sigma^\eta$ 
of partial isometries
satisfying 
\begin{align*}
\sum_{\beta \in \Sigma^\rho} s_{\beta}s_{\beta}^* & =1,\qquad
\pi(x) s_\alpha s_\alpha^* =  s_\alpha s_\alpha^*  \pi(x),\qquad
s_\alpha^* \pi(x) s_\alpha = \pi(\rho_\alpha(x)), \\
\sum_{b \in \Sigma^\eta} t_{b}t_{b}^* & =1,\qquad
\pi(x) t_a t_a^* =  t_a t_a^*  \pi(x),\qquad
t_a^* \pi(x) t_a = \pi(\eta_a(x)), \\  
& s_\alpha t_b =  t_a s_\beta \quad \text{ if } \quad \kappa(\alpha, b) = (a,\beta)
\end{align*}
for all $x \in \A$ and $\alpha \in \Sigma^\rho$, 
$a \in \Sigma^\eta$.
 Put
 $\widetilde{\A} = \pi(\A)$
 and
$
\tilde{\rho}_\alpha(\pi(x)) = \pi(\rho_\alpha(x)), 
\tilde{\eta}_a(\pi(x)) = \pi(\eta_a(x)), x \in \A$.
It is easy to see that 
$(\widetilde{\A},\tilde{\rho}, \tilde{\eta}, \Sigma^\rho, \Sigma^\eta, \kappa)$ 
is a $C^*$-textile dynamical system
such that the presented two-dimensional textile dynamical system 
$X_{\tilde{\rho},\tilde{\eta}}^\kappa$
is the same as the one $X_{\rho,\eta}^\kappa$
presented by
$\CTDS$.
Let
${\cal O}_{\pi,s,t}$
be the $C^*$-subalgebra of $\B$ generated by  
$\pi(x)$ and $s_\alpha$, $t_a$
 for $ x \in \A, \alpha \in \Sigma^\rho, a \in \Sigma^\eta$.
Let
$
{\cal F}_{\pi,s,t} 
$
be the $C^*$-subalgebra of 
${\cal O}_{\pi,s,t}$
generated by
$ s_\mu t_\zeta \pi(x) t_\xi^* s_{\nu}^*$
for
$ 
x \in \A
$
and
$
\mu, \nu \in B_*(\Lambda_\rho), \zeta, \xi \in B_*(\Lambda_\eta)
$
with
$
|\mu| = |\nu|, |\zeta| = |\xi|.
$
By the universality of the algebra $\ORE$,
the correspondence
$$ 
x \in \A \longrightarrow \pi(x) \in \widetilde{A},
\qquad
S_\alpha \longrightarrow s_\alpha, \ \ \alpha \in \Sigma^\rho,
\qquad
T_a \longrightarrow t_a, \ \ a \in \Sigma^\eta
$$
extends to a surjective $*$-homomorphism
$\tilde{\pi}:\ORE \longrightarrow {\cal O}_{\pi,s,t}$. 
\begin{lem}
The restriction of $\tilde{\pi}$ to the subalgebra
${\cal F}_{\rho,\eta}$
is a $*$-isomorphism from
${\cal F}_{\rho,\eta}$ to ${\cal F}_{\pi,s,t}$.
Hence 
if $\CTDS$ satisfies condition (I), 
so does 
$(\widetilde{\A},\tilde{\rho},\tilde{\eta}, \Sigma^\rho,\Sigma^\eta,\kappa)$.
\end{lem}
\begin{pf}
It suffices to show that $\tilde{\pi}$ is injective on
${\cal F}_{j,k}$
for all $j,k \in {\Bbb Z}$.
Suppose  
$$
\sum_{\mu,\nu\in B_j(\Lambda_\rho), \zeta, \xi \in B_k(\Lambda_\eta)}
s_\mu t_\zeta \pi(x_{\mu,\zeta,\xi,\nu}) t_\xi^* s_\nu^*=0
$$
for
$
\sum_{\mu,\nu\in B_j(\Lambda_\rho), \zeta, \xi \in B_k(\Lambda_\eta)}
S_\mu T_\zeta x_{\mu,\zeta,\xi,\nu} T_\xi^* S_\nu^*
\in
{\cal F}_{j,k}
$
with
$x_{\mu,\zeta,\xi,\nu}\in \A$.
For 
$
\mu',\nu'\in B_j(\Lambda_\rho),
\zeta', \xi' \in B_k(\Lambda_\eta),
$
one has 
\begin{align*}
&\pi(\eta_{\zeta'}(\rho_{\mu'}(1)) x_{\mu',\zeta', \xi',\nu'}\eta_{\xi'}(\rho_{\nu'}(1)))\\
= 
&
t_{\zeta'}^* s_{\mu'}^*(
\sum_{\mu,\nu\in B_j(\Lambda_\rho), \zeta, \xi \in B_k(\Lambda_\eta)}
s_\mu t_\zeta \pi(x_{\mu,\zeta,\xi,\nu}) t_\xi^* s_\nu^*)
s_{\nu'} t_{\xi'}
=0.
\end{align*}
As $\pi: \A \longrightarrow \B$ is injective,
one sees
$$
\eta_{\zeta'}(\rho_{\mu'}(1)) x_{\mu',\zeta', \xi',\nu'}\eta_{\xi'}(\rho_{\nu'}(1))=0
$$
so that
$$
 S_{\mu'} T_{\zeta'} 
x_{\mu',\zeta', \xi',\nu'} T_{\xi'}^* S_{\nu'}^*=0.
$$
Hence we have
$$
\sum_{\mu,\nu\in B_j(\Lambda_\rho), \zeta, \xi \in B_k(\Lambda_\eta)}
S_\mu T_\zeta x_{\mu,\zeta,\xi,\nu} T_\xi^* S_\nu^* =0.
$$
Therefore
 $\tilde{\pi}$ is injective on
${\cal F}_{j,k}$.
\end{pf}

We henceforth assume that $\CTDS$ satisfies condition (I) defined above.
Take a unital increasing sequence $\{ \A_l \}_{l \in \Zp}$
of $C^*$-subalgebras of $\A$
as in the definition of condition (I).
Recall that  the algebra 
$
 {\cal F}_{j,k}^l
$
for $j,k \le l$
is defined as 
$$
 {\cal F}_{j,k}^l 
 =   
  C^*( S_{\mu} T_\zeta x T_\xi^* S_{\nu}^{*}: 
\mu, \nu \in B_j(\Lambda_\rho), \zeta, \xi \in B_k(\Lambda_\eta), x \in \A_l). 
$$
There exists an inclusion relation 
${\cal F}_{j,k}^l \subset {\cal F}_{j',k'}^{l'}$ 
for $j\le j', k \le k'$ and $l \le l'$
through the identities (4.1), (4.2). 

Let 
${\cal P}_{\pi,s,t}$ be the $*$-subalgebra of ${\cal O}_{\pi,s,t}$
algebraically generated by 
$\pi(x),s_\alpha, t_a$ for $x \in \A_l, l \in \Zp,\ \alpha\in \Sigma^\rho, a \in \Sigma^\eta$.
\begin{lem}
Any element  $x \in {\cal P}_{\pi,s,t}$ can be expressed in a unique way as  
\begin{align*}
x =&  \sum_{|\nu|,   |\xi| \ge 1} x_{-\xi,-\nu} t_\xi^* s_\nu^* 
  +   \sum_{|\zeta|, |\nu| \ge 1} t_\zeta x_{\zeta,-\nu}  s_\nu^* 
  +   \sum_{|\mu|,   |\xi| \ge 1} s_\mu x_{\mu,-\xi} t_\xi^* 
  +   \sum_{|\mu|,  \zeta| \ge 1} s_\mu t_\zeta x_{\mu,\zeta} \\ 
  + & \sum_{ |\xi| \ge 1} x_{-\xi} t_\xi^*  
  +   \sum_{ |\nu| \ge 1} x_{-\nu} s_\nu^* 
  +   \sum_{|\mu| \ge 1} s_\mu x_{\mu}
  +   \sum_{|\zeta| \ge 1} t_\zeta  x_{\zeta}
  +   x_0 
\end{align*}
where
$
x_{-\xi,-\nu}, x_{\zeta,-\nu}, x_{\mu,-\xi}, x_{\mu,\zeta},
x_{-\xi}, x_{-\nu}, x_{\mu}, x_{\zeta},
x_0  \in {\cal P}_{\pi,s,t} \cap {\cal F}_{\pi,s,t}
$
for
$
 \mu,\nu \in B_*(\Lambda_\rho),  \zeta, \xi \in B_*(\Lambda_\eta),
$
which satisfy
\begin{align*}
& x_{-\xi,-\nu} = x_{-\xi,-\nu} \eta_\xi(\rho_\nu(1)), \quad 
  x_{\zeta,-\nu} = \eta_\zeta(1) x_{\zeta,-\nu}\rho_\nu(1), \\ 
& x_{\mu,-\xi} = \rho_\mu(1) x_{\mu,-\xi} \eta_\xi(1), \quad
  x_{\mu,\zeta} = \eta_\zeta(\rho_\mu(1)) x_{\mu,\zeta},\\ 
& x_{-\xi} = x_{-\xi} \eta_\xi(1),\quad  
  x_{-\nu} = x_{-\nu} \rho_\nu(1), \quad 
  x_{\mu} = \rho_\mu(1) x_{\mu}, \quad 
  x_{\zeta} = \eta_\zeta(1) x_{\zeta}.
\end{align*}
\end{lem}
\begin{pf}
Put
\begin{align*}
& x_{-\xi,-\nu} = {\cal E}_{\rho,\eta}(x s_{\nu}t_\xi), \quad 
  x_{\zeta,-\nu} = {\cal E}_{\rho,\eta}(t_\zeta^* x s_{\nu}), \\ 
& x_{\mu,-\xi} = {\cal E}_{\rho,\eta}(s_\mu^* x t_{\xi}), \quad
  x_{\mu,\zeta} = {\cal E}_{\rho,\eta}(t_\zeta^* s_\mu^* x),\\ 
& x_{-\xi} = {\cal E}_{\rho,\eta}(x t_\xi),\quad  
  x_{-\nu} = {\cal E}_{\rho,\eta}(x s_{\nu}), \quad 
  x_{\mu} = {\cal E}_{\rho,\eta}(s_\mu^* x), \quad 
  x_{\zeta} = {\cal E}_{\rho,\eta}(t_\zeta^* x), \\ 
& x_0 =  {\cal E}_{\rho,\eta}(x).
\end{align*}
Then we have a desired expression of $x$.
\end{pf}
\begin{lem}
For $h \in {\cal D}_{\rho,\eta}\cap \A_l^{\prime}$
and
$j, k \in {\Bbb Z}$ 
with 
$j + k \le l$,
put
$h^{j,k} = \phi_\rho^j\circ \phi_\eta^k(h)$.
Then we have
\begin{enumerate}
\renewcommand{\labelenumi}{(\roman{enumi})}
\item
$h^{j,k} s_\mu = s_\mu h^{j-|\mu|,k}$ for
$\mu \in B_*(\Lambda_\rho)$ with $|\mu|\le j$.
\item
$h^{j,k} t_\zeta = t_\zeta h^{j,k-|\zeta|}$ for
$\zeta \in B_*(\Lambda_\eta)$ with $|\zeta|\le k$.
\item
$h^{j,k}$ commutes with any element of 
${\cal F}_{j,k}^l$.
\end{enumerate} 
\end{lem}
\begin{pf}
(i)
It follows that for
$\mu \in B_*(\Lambda_\rho)$ with $|\mu|\le j$
\begin{equation*}
h^{j,k} s_\mu 
 = \sum_{|\mu'| = |\mu|} s_{\mu'}\phi_\rho^{j - |\mu|}(\phi_\eta^k(h)) s_{\mu'}^* s_\mu 
 = s_{\mu}\phi_\rho^{j - |\mu|}(\phi_\eta^k(h)) s_{\mu}^* s_\mu.
\end{equation*}
Since
$h \in \A_l^\prime$
and
$\A_{j+k}\subset \A_l$,
one has
\begin{align*}
\phi_\rho^{j - |\mu|}(\phi_\eta^k(h)) s_{\mu}^* s_\mu
& = \sum_{\nu \in B_{j-|\mu|}(\Lambda_\rho)}
    \sum_{\xi\in B_k(\Lambda_\eta)} 
    s_{\nu} t_{\xi} h t_{\xi}^* s_{\nu}^*  s_\mu^* s_\mu \\ 
& = \sum_{\nu \in B_{j-|\mu|}(\Lambda_\rho)}
    \sum_{\xi\in B_k(\Lambda_\eta)} 
    s_{\nu} t_{\xi} h t_{\xi}^* s_{\nu}^*  s_\mu^* s_\mu s_{\nu} t_{\xi}t_{\xi}^* s_{\nu}^* \\ 
& = \sum_{\nu \in B_{j-|\mu|}(\Lambda_\rho)}
    \sum_{\xi\in B_k(\Lambda_\eta)} 
    s_{\nu} t_{\xi} h \eta_{\xi}(\rho_{\mu\nu}(1)) t_{\xi}^* s_{\nu}^* \\ 
& = \sum_{\nu \in B_{j-|\mu|}(\Lambda_\rho)}
    \sum_{\xi\in B_k(\Lambda_\eta)} 
    s_{\nu} t_{\xi} \eta_{\xi}(\rho_{\mu\nu}(1))  h t_{\xi}^* s_{\nu}^* \\ 
& = \sum_{\nu \in B_{j-|\mu|}(\Lambda_\rho)}
    \sum_{\xi\in B_k(\Lambda_\eta)} 
    s_{\nu} \rho_{\mu\nu}(1)  t_{\xi} h t_{\xi}^* s_{\nu}^* \\ 
&  = s_{\mu}^*s_{\mu}\phi_\rho^{j - |\mu|}(\phi_\eta^k(h)) 
   = s_{\mu}^*s_\mu h^{j-|\mu|,k}
\end{align*}
so that
$h^{j,k} s_\mu = s_\mu h^{j-|\mu|,k}$.

(ii) Similarly we have 
$h^{j,k} t_\zeta = t_\zeta h^{j,k-|\zeta|}$ for
$\zeta \in B_*(\Lambda_\eta)$ with $|\zeta|\le k$.

(iii)
For 
$x \in \A_l, 
 \mu,\nu \in B_j(\Lambda_\rho),  \zeta, \xi \in B_k(\Lambda_\eta),
$
we have
$$
h^{j,k} s_\mu t_\zeta = s_\mu h^{0,k}t_\zeta = s_\mu t_\zeta h^{0,0}= s_\mu t_\zeta h.
$$
It follows that
 \begin{equation*}
h^{j,k} s_\mu t_\zeta x t_\xi^* s_\nu^* 
= s_\mu t_\zeta h x t_\xi^* s_\nu^* 
 = s_\mu t_\zeta x h t_\xi^* s_\nu^* 
 = s_\mu t_\zeta x t_\xi^* s_\nu^* h^{j,k}.
\end{equation*}
Hence 
$h^{j,k}$ commutes with any element of 
${\cal F}_{j,k}^l$.
\end{pf}

\begin{lem}
Assume that $\CTDS$ satisfies condition (I).
Let  $x \in {\cal P}_{\pi,s,t}$  be expressed as in the preceding lemma. 
Then we have
$$
\| x_0 \| \le \| x \|.
$$
\end{lem}
\begin{pf}
We may assume that 
for $x \in {\cal P}_{\pi,s,t}$,
$$
x_{-\xi,-\nu}, x_{\zeta,-\nu}, x_{\mu,-\xi}, x_{\mu,\zeta},
x_{-\xi}, x_{-\nu}, x_{\mu}, x_{\zeta},
x_0  \in \tilde{\pi}( {\cal F}_{j_1,k_1}^{l_1})
$$
for some $j_1, k_1,l_1$
and
$\mu, \nu \in \cup_{n=0}^{j_0} B_n(\Lambda_\rho)$,
$\zeta, \xi \in \cup_{n=0}^{k_0} B_n(\Lambda_\eta)$
for some $j_0, k_0$.
Take $j,k,l \in \Zp$
such as
\begin{equation*}
j \ge j_0 + j_1, \qquad
k \ge k_0 + k_1, \qquad
l \ge \max\{j+k, l_1\}.
\end{equation*}
By Lemma 4.1, 
$(\widetilde{\A},\tilde{\rho}, \tilde{\eta},\Sigma^\rho,\Sigma^\eta,\kappa)$
satisfies condition (I).
For any $\epsilon > 0$ and the numbers
$j,k,l$, the element 
$
x_0 \in \tilde{\pi}( {\cal F}_{j_1,k_1}^{l_1}),
$
one may find 
$g \in \tilde{\pi}({\cal D}_{\rho,\eta}) \cap \pi(\A_l)^\prime$
with $0 \le g \le 1$
such that
\begin{enumerate}
\renewcommand{\labelenumi}{(\roman{enumi})}
\item
 $\| x_0 \phi_{\rho}^j \circ \phi_{\eta}^k(g) \| \ge \| x_0 \| - \epsilon$,
\item
 $
 \phi_{\rho}^n(g) \phi_{\eta}^m(g) 
=\phi_{\rho}^n(( \phi_{\eta}^m(g)))g 
=\phi_{\rho}^n(g) g 
=\phi_{\eta}^m(g) g = 0
$ 
for all $n=1,2,\dots,j$, $m= 1,2,\dots, k.$
\end{enumerate}
Put
$
h = g^{\frac{1}{2}}
$
and
$
h^{j,k} =  \phi_{\rho}^j\circ \phi_{\eta}^k(h).
$ 
It follows that
\begin{align*}
 \| x \| 
\ge & \|  h^{j,k} x h^{j,k} \| \\
= &  \| \sum_{|\nu|,   |\xi| \ge 1} h^{j,k} x_{-\xi,-\nu} t_\xi^* s_\nu^* h^{j,k} \qquad (1)\\ 
  & +  \sum_{|\zeta|, |\nu| \ge 1} h^{j,k} t_\zeta x_{\zeta,-\nu} s_\nu^*h^{j,k} \qquad (2)\\ 
  & +  \sum_{|\mu|,   |\xi| \ge 1} h^{j,k} s_\mu x_{\mu,-\xi} t_\xi^* h^{j,k} \qquad (3)\\ 
  & +  \sum_{|\mu|,  \zeta| \ge 1} h^{j,k} s_\mu t_\zeta x_{\mu,\zeta} h^{j,k} \qquad (4)\\ 
  & +  \sum_{ |\xi| \ge 1} h^{j,k} x_{-\xi} t_\xi^*h^{j,k}  
  +   \sum_{ |\nu| \ge 1} h^{j,k} x_{-\nu} s_\nu^*h^{j,k} 
  +   \sum_{|\mu| \ge 1} h^{j,k}  s_\mu x_{\mu}h^{j,k}
  +   \sum_{|\zeta| \ge 1} h^{j,k} t_\zeta  x_{\zeta}h^{j,k} \quad (5) \\
  & +   h^{j,k} x_0 h^{j,k} \|
\end{align*}
For (1),
as 
$x_{-\xi,-\nu} \in \tilde{\pi}( {\cal F}_{j_1,k_1}^{l_1}) 
\subset \tilde{\pi}( {\cal F}_{j,k}^{l})$,
one sees that $x_{-\xi,-\nu}$
commutes with
$h^{j,k}$.
Hence we have
$$
h^{j,k} x_{-\xi,-\nu} t_\xi^* s_\nu^* h^{j,k}
= x_{-\xi,-\nu} h^{j,k}  t_\xi^* s_\nu^* h^{j,k}
= x_{-\xi,-\nu} h^{j,k}  h^{j-|\nu|,k-|\xi|} t_\xi^* s_\nu^*
$$
and
\begin{align*}
h^{j,k}  h^{j-|\nu|,k-|\xi|}(h^{j,k}  h^{j-|\nu|,k-|\xi|})^*
=& \phi_{\rho}^j ( \phi_{\eta}^k(g)) \cdot 
   \phi_{\rho}^{j-|\nu|}(\phi_{\eta}^{k-|\xi|}(g))\\
=& \phi_{\rho}^{j-|\nu|}\circ \phi_{\eta}^{k-|\xi|}(
\phi_{\eta}^{|\xi|}(\phi_{\rho}^{|\nu|}(g)g)) =0
\end{align*}
so that
\begin{equation*}
h^{j,k} x_{-\xi,-\nu} t_\xi^* s_\nu^* h^{j,k}
= 0.
\end{equation*}

For (2),
as 
$x_{\xi,-\nu} \in \tilde{\pi}( {\cal F}_{j_1,k_1}^{l_1}) 
\subset \tilde{\pi}( {\cal F}_{j,k-|\xi|}^{l})$,
one sees $x_{\xi,-\nu}$
that
commutes with
$h^{j,k-|\xi|}$.
Hence we have
$$
h^{j,k} t_\xi x_{\xi,-\nu}  s_\nu^* h^{j,k}
= t_\xi  h^{j,k-|\xi|} x_{\xi,-\nu}  h^{j-|\nu|,k}s_\nu^* 
= t_\xi  x_{\xi,-\nu} h^{j,k-|\xi|}  h^{j-|\nu|,k}s_\nu^*
$$
and
\begin{align*}
h^{j,k-|\xi|}   h^{j-|\nu|,k}(h^{j,k-|\xi|}   h^{j-|\nu|,k})^*
=& \phi_{\rho}^j (\phi_{\eta}^{k-|\zeta|}(g))\cdot
\phi_{\rho}^{j-|\nu|}( \phi_{\eta}^{k}(g) ) \\
=& \phi_{\rho}^{j-|\nu|}\circ \phi_{\eta}^{k-|\zeta|}(
\phi_{\rho}^{|\nu|}(g) \phi_{\eta}^{|\zeta|}(g)) =0
\end{align*}
so that
\begin{equation*}
h^{j,k} t_\xi x_{\xi,-\nu}  s_\nu^* h^{j,k}
= 0.
\end{equation*}

For (3),
as 
$x_{\mu,-\xi} \in 
\tilde{\pi}( {\cal F}_{j_1,k_1}^{l_1})
 \subset 
\tilde{\pi}( {\cal F}_{j-|\mu|,k}^{l})$,
one sees $x_{\mu,-\xi}$ that
commutes with
$h^{j-|\mu|,k}$.
Hence we have
$$
h^{j,k} s_\mu  x_{\mu,-\xi} t_\xi^* h^{j,k}
= s_\mu  h^{j-|\mu|,k}  x_{\mu,-\xi} h^{j,k-|\xi|} t_\xi^*
= s_\mu    x_{\mu,-\xi}h^{j-|\mu|,k} h^{j,k-|\xi|} t_\xi^*
$$
and
\begin{align*}
h^{j-|\mu|,k} h^{j,k-|\xi|} (h^{j-|\mu|,k} h^{j,k-|\xi|})^* 
= & \phi_{\rho}^{j-|\mu|} ( \phi_{\eta}^k(g)) \cdot
    \phi_{\rho}^{j} ( \phi_{\eta}^{k-|\xi|}(g)) \\
 =& \phi_{\rho}^{j-|\mu|}\circ \phi_{\eta}^{k-|\xi|}(
    \phi_{\eta}^{|\xi|}(g) \phi_{\rho}^{|\mu|}(g)) =0
\end{align*}
so that
\begin{equation*}
h^{j,k} s_\mu  x_{\mu,-\xi} t_\xi^* h^{j,k}
= 0.
\end{equation*}

For (4),
as 
$x_{\mu, \zeta} \in 
\tilde{\pi}( {\cal F}_{j_1,k_1}^{l_1})
\subset 
\tilde{\pi}( {\cal F}_{j-|\mu|, k-|\zeta|}^{l})$,
one sees $x_{\mu,\zeta}$
that
commutes with
$h^{j-|\mu|,k-|\zeta|}$.
Hence we have
$$
h^{j,k} s_\mu  t_\zeta x_{\mu, \zeta}  h^{j,k}
= s_\mu  t_\zeta h^{j-|\mu|,k-|\zeta|}  x_{\mu, \zeta}  h^{j,k}
= s_\mu  t_\zeta   x_{\mu, \zeta}  h^{j-|\mu|,k-|\zeta|} h^{j,k}
$$
and
\begin{align*}
h^{j-|\mu|,k-|\zeta|}h^{j,k}(h^{j-|\mu|,k-|\zeta|}h^{j,k})^* 
= & \phi_{\rho}^{j-|\mu|} ( \phi_{\eta}^{k-|\zeta|}(g))
    \phi_{\rho}^{j}(\phi_{\eta}^{k}(g)) \\
 =& \phi_{\rho}^{j-|\mu|}\circ \phi_{\eta}^{k-|\xi|}
    (g  \phi_{\rho}^{|\mu|} \circ \phi_{\eta}^{|\xi|}(g)) =0
\end{align*}
so that
\begin{equation*}
h^{j,k} s_\mu  t_\zeta x_{\mu, \zeta}  h^{j,k}
= 0.
\end{equation*}

For (5)
as 
$x_{-\xi}$
commutes with
$h^{j, k}$,
we have
$$
h^{j,k}  x_{-\xi}  t_\xi^*  h^{j,k}
= x_{-\xi} h^{j,k} h^{j,k-|\xi|}t_\xi^*
$$
and
\begin{align*}
h^{j,k} h^{j,k-|\xi|} (h^{j,k}h^{j,k-|\xi|})^* 
= & \phi_{\rho}^{j} ( \phi_{\eta}^{k|}(g))
    \phi_{\rho}^{j}(\phi_{\eta}^{k-|\xi|}(g)) \\
 =& \phi_{\rho}^{j}\circ \phi_{\eta}^{k-|\xi|}
    ( \phi_{\eta}^{|\xi|}(g)) =0
\end{align*}
so that
\begin{equation*}
h^{j,k}  x_{-\xi}  t_\xi^*  h^{j,k}
= 0.
\end{equation*}
We similarly see that
\begin{equation*}
h^{j,k} x_{-\nu} s_\nu^*h^{j,k} 
= h^{j,k}  s_\mu x_{\mu}h^{j,k}
= h^{j,k} t_\zeta  x_{\zeta}h^{j,k}
= 0.
\end{equation*}
Therefore we have
$$
\| x \| 
\ge \| h^{j,k} x_0 h^{j,k} \|
= \|  x_0 (h^{j,k})^2 \|
= \| x_0 \phi_\rho^j \circ \phi_\eta^k(g) \| 
\ge \| x_0 \| - \epsilon.
$$
Hence we get
$ \| x \| \ge \| x_0 \|$.
\end{pf}

By a similar argument ot \cite[2.8 Proposition]{CK}, one sees
\begin{cor}
Assume that
$\CTDS$ satisfies condition (I).
There exists a conditional expectation
$\E_{\pi,s,t}:{\cal O}_{\pi,s,t} \longrightarrow {\cal F}_{\pi,s,t}
$
such that
$\E_{\pi,s,t} \circ \tilde{\pi} =  \tilde{\pi}\circ \E_{\rho,\eta}.
$
\end{cor}

Therefore we have
\begin{prop}
Assume that
$\CTDS$ satisfies condition (I).
The $*$-homomorphism $\tilde{\pi}:\ORE \longrightarrow {\cal O}_{\pi,s,t}$
defined by
$$
\tilde{\pi}(x) = \pi(x),\quad x \in \A, \qquad
\tilde{\pi}(S_\alpha) = s_\alpha, \quad \alpha \in \Sigma^\rho,
\qquad
\tilde{\pi}(T_a) = t_a, \quad a \in \Sigma^\eta
$$
becomes a surjective $*$-isomorphism, and hence
the $C^*$-algebras $\ORE$ and ${\cal O}_{\pi,s}$ 
are canonically $*$-isomorphic through $\tilde{\pi}$.
\end{prop}
\begin{pf}
The map $\tilde{\pi}: {\cal F}_{\rho,\eta} \rightarrow {\cal F}_{\pi,s,t}$
is $*$-isomorphic and satisfies 
$\E_{\pi,s,t} \circ \tilde{\pi} =  \tilde{\pi}\circ \E_{\rho,\eta}.
$
Since 
$\E_\rho:\ORE \longrightarrow {\cal F}_{\rho,\eta}$
is faithful,
a routine argument shows that the $*$-homomorphism 
$\tilde{\pi}:\ORE \longrightarrow {\cal O}_{\pi,s,t}$
is actually a $*$-isomorphism.
\end{pf}
Hence the following uniqueness of the $C^*$-algebra $\ORE$ holds.
\begin{thm}
Assume that
$\CTDS$ satisfies condition (I).
The $C^*$-algebra $\ORE$ is the unique $C^*$-algebra subject to the relation 
$(\rho,\eta;\kappa)$.
This means that if there exist a unital $C^*$-algebra $\B$ and an
injective $*$-homomorphism $\pi: \A \longrightarrow \B$ and two families 
of partial isometries 
$s_{\alpha}, \alpha \in \Sigma^\rho$,
$t_a, a \in \Sigma^\eta$
satisfying the following relations :
\begin{align*}
\sum_{\beta \in \Sigma^\rho}  s_{\beta}s_{\beta}^*   =1,& \qquad
\pi(x) s_\alpha s_\alpha^*  =  s_\alpha s_\alpha^* \pi(x),\qquad
s_\alpha^* \pi(x) s_\alpha  = \pi(\rho_\alpha(x)),\\
\sum_{b \in \Sigma^\eta}  t_{b} t_{b}^* =1,& \qquad
\pi(x) t_a t_a^*  =  t_a t_a^* \pi(x),\qquad
t_a^* \pi(x) t_a  = \pi(\eta_a(x)) \\
s_\alpha t_b  & = t_a s_\beta  \qquad \text{ if } \quad \kappa(\alpha,b) = (a,\beta)  
\end{align*}
for
$(\alpha,b) \in \Sigma^{\rho\eta},
 (a,\beta) \in \Sigma^{\eta\rho}$
and $x \in {\cal A}$, $\alpha \in \Sigma^\rho, a \in \Sigma^\eta$,
then the correspondence 
$$
 x \in \A \longrightarrow \pi(x) \in \B, \quad
S_\alpha \longrightarrow s_\alpha \in \B, \qquad
T_a \longrightarrow t_a \in \B
$$ 
extends to a $*$-isomorphism
$\tilde{\pi}$ from $\ORE$ onto the $C^*$-subalgebra
${\cal O}_{\pi,s,t}$ of $\B$ 
generated by $\pi(x), x \in \A$ and 
$s_\alpha, \alpha \in \Sigma$,
$t_a, a \in \Sigma^\eta$.
\end{thm}

For a $C^*$-textile dynamical system 
$\CTDS$,
let
$\lambda_{\rho,\eta}: \A \rightarrow \A$
be the positive map on $\A$ defined by
\begin{equation*}
\lambda_{\rho, \eta}(x) = \sum_{\alpha \in \Sigma^\rho, a \in \Sigma^\eta}
\eta_a\circ \rho_\alpha(x), \qquad x \in \A.
\end{equation*} 
Then
$\CTDS$
is said to be  {\it irreducible\/} if
there exists no nontrivial ideal of $\A$ invariant under $\lambda_{\rho,\eta}$.
\begin{cor}
If $\CTDS$ 
satisfies condition (I)
and  is irreducible,
the $C^*$-algebra $\ORE$ is simple.
\end{cor}
\begin{pf}
Assume that
there exists a nontrivial ideal
${\cal I}$ of $\ORE$.
Now suppose that  
${\cal I} \cap \A = \{ 0 \}$.
As 
$S_\alpha^* S_\alpha = \rho_\alpha(1),
T_a^* T_a = \eta_a(1) \in \A$
one knows that
$S_\alpha, T_a  \not\in {\cal I}$ 
for all $\alpha \in \Sigma^\rho, a \in \Sigma^\eta$.
By the preceding theorem, 
the quotient map
$q: \ORE \longrightarrow \ORE/{\cal I}$
must be injective so that ${\cal I}$ is trivial.
Hence one sees that
${\cal I} \cap \A \ne \{ 0\}$  
and it is invariant under $\lambda_{\rho,\eta}$.
\end{pf}

\section{Cocrete realization}
In this section we will realize the $C^*$-algebra
$\ORE$ for 
$\CTDS$ in a concrete way.
For
$\gamma_i \in \Sigma^\rho\cup\Sigma^\eta$,
put\begin{equation*}
\xi_{\gamma_i}
=
\begin{cases}
\rho_{\gamma_i} & \text{ if } {\gamma_i} \in \Sigma^\rho,\\
\eta_{\gamma_i} & \text{ if } {\gamma_i} \in \Sigma^\eta.
\end{cases}
\end{equation*}

\noindent

{\bf Definition.}
A finite sequence of labeles
$
(\gamma_1,\gamma_2,\dots,\gamma_k) \in  (\Sigma^\rho\cup\Sigma^\eta)^k
$
is said to be {\it concatenated labeled path}\
if
$\xi_{\gamma_k} \circ \cdots\circ 
\xi_{\gamma_2} \circ \xi_{\gamma_1} (1) \ne 0.
$ 
For
$m,n \in \Zp$,
let
$L_{(n,m)}$
be the set of
concatenated labeled paths
$
(\gamma_1,\gamma_2,\dots,\gamma_{m+n}) 
$
such that 
symbols in $\Sigma^\rho$ appear in $
(\gamma_1,\gamma_2,\dots,\gamma_{m+n}) 
$
$n$-times and
symbols in $\Sigma^\eta$ appear in 
$
(\gamma_1,\gamma_2,\dots,\gamma_{m+n}) 
$
$m$-times.
We define a relation in $L_{(n,m)}$ for
$i=1,2,\dots, n+m-1$.
We write
\begin{equation*}
(\gamma_1,\dots,\gamma_{i-1}, \gamma_i, \gamma_{i+1},\gamma_{i+2},\dots,\gamma_{m+n}) 
\underset{i}{\approx}
(\gamma_1,\dots,\gamma_{i-1}, \gamma'_i, \gamma'_{i+1},\gamma_{i+2},\dots,\gamma_{m+n}) 
\end{equation*}
if one of the following two conditions holds:

(1)
$(\gamma_i, \gamma_{i+1} )\in \Sigma^{\rho\eta}, 
(\gamma'_i, \gamma'_{i+1} )\in \Sigma^{\eta\rho}
$
and
$
\kappa(\gamma_i, \gamma_{i+1} )= 
(\gamma'_i, \gamma'_{i+1} ),
$

(2)
$(\gamma_i, \gamma_{i+1} )\in \Sigma^{\eta\rho}, 
(\gamma'_i, \gamma'_{i+1} )\in \Sigma^{\rho\eta}
$
and
$
\kappa(\gamma'_i, \gamma'_{i+1} )= 
(\gamma_i, \gamma_{i+1} ).
$

Denote by
$\approx$
the equivalence relation in
$L_{(n,m)}$ generated by the relations
$
\underset{i}{\approx}, i=1,2,\dots, n+m-1.
$
Let
${\frak T}_{(n,m)} = L_{(n,m)} / \approx
$
be the set of equivalence classes of 
$
L_{(n,m)}
$ 
under 
$\approx$.
Denote by $[\gamma] \in {\frak T}_{(n,m)} $ 
the equivalence class of $\gamma \in L_{(n,m)}$. 
Put the vectors
$e = (1, 0), f = (0, -1)$ in ${\Bbb R}^2$.
Consider the set of all paths 
consisting of sequences of vectors $e,f$
starting at the  point
$(-n, m) \in {\Bbb R}^2$
for $n,m\in \Zp$
and ending at the origin.
Such a path consists of $n$ $e$-vectors and $m$ $f$-vectors. 
Let
${\frak P}_{(n,m)}$ be the set of all such paths from
$(-n,m)$ to the origin.
We consider the correspondence
\begin{equation*}
\rho_\alpha \longrightarrow e \quad (\alpha \in \Sigma^\rho),
\qquad
\eta_a \longrightarrow f \quad (a \in \Sigma^\eta),
\end{equation*}
denoted by $\pi$.
It extends from
$L_{(n,m)}$ to ${\frak P}_{(n,m)}$
in a natural way.
The following lemma is obvious.
\begin{lem}
For any path $p \in {\frak P}_{(n,m)}$ of vectors,
there uniquely exists
a concatenated labeled path
$\gamma \in L_{(n,m)}$
such that
$\pi(\gamma) = p$.
\end{lem}
For 
a concatenated labeled path
$\gamma=(\gamma_1,\gamma_2,\dots,\gamma_{n+m}) \in L_{(n,m)}$,
put
the projection in $\A$
$$
P_\gamma = \xi_{\gamma_k} \circ \cdots\circ 
\xi_{\gamma_2} \circ \xi_{\gamma_1} (1).
$$ 
We note that 
$P_\gamma \ne 0$
for all $\gamma \in L_{(n,m)}$.
\begin{lem}
For $\gamma, \gamma' \in L_{(n,m)}$,
if 
$\gamma \approx \gamma'$, we have
$P_\gamma = P_{\gamma'}$.
Hence the projection 
$P_{[\gamma]}$ for $[\gamma]\in {\frak T}_{(n,m)}$
is well-defined. 
\end{lem}
\begin{pf}
If $\kappa(\alpha,b) = (a,\beta)$,
one has
$\eta_b \circ \rho_\alpha(1) = \rho_\beta \circ \eta_a(1) \ne 0$.
Hence we know the assertion.
\end{pf}
Denote by $| {\frak T}_{(n,m)} |$
the cardinal number of the finite set 
${\frak T}_{(n,m)}$.
Let $e_t, t \in {\frak T}_{(n,m)}$
be the standard complete orthonomal basis of
${\Bbb C}^{|{\frak T}_{(n,m)}|}$.
Define 
\begin{align*}
H_{(n,m)} & = \sum_{t \in {\frak T}_{(n,m)}}{}^\oplus {\Bbb C}e_t \otimes P_t\A  \\
         (& = \sum_{t \in {\frak T}_{(n,m)}}{}^\oplus 
  \Span\{ c e_t \otimes P_tx \mid c \in {\Bbb C}, x \in \A \}  )
\end{align*}
the direct sum of ${\Bbb C}e_t \otimes P_t\A $ over $t \in {\frak T}_{(n,m)}$.
$H_{(n,m)}$ has a structure of Hilbert $C^*$-bimodule over $\A$ by setting
\begin{align*}
(e_t \otimes P_t x)y & := e_t \otimes P_t xy, \\
\phi(y)(e_t \otimes P_t x) & := e_t \otimes \xi_\gamma(y)x (= e_t \otimes P_t \xi_\gamma(y) x), 
\end{align*}
where $t= [\gamma]$ for $ \gamma = (\gamma_1,\dots,\gamma_{n+m}) $
and
$\xi_\gamma(y) = 
\xi_{\gamma_{n+m}} \circ \cdots\circ 
\xi_{\gamma_2} \circ \xi_{\gamma_1} (y)
$
and
\begin{equation*}
\langle e_t \otimes P_t x \mid e_s \otimes P_s y \rangle :=
\begin{cases}
x^* P_t y & \text{ if } t = s,\\
0              & \text{ otherwise }
\end{cases}
\end{equation*}
for $t, s \in {\frak T}_{(n,m)}$ and $x, y \in \A$.
Put
$H_{(0,0)} = \A$.
Denote by 
$F(\rho,\eta)$
the Hilbert $C^*$-bimodule over $\A$
defined by the direct sum: 
\begin{equation*}
F(\rho,\eta) = \sum_{(n,m) \in {\Bbb Z}^2}{}^\oplus H_{(n,m)}. 
\end{equation*}
For $\alpha \in \Sigma^\rho, a \in \Sigma^\eta$,
the creation operators 
$s_\alpha, t_a $ on $F(\rho,\eta):$
\begin{equation*}
s_\alpha  : H_{(n,m)} \longrightarrow  H_{(n+1,m)},\qquad
t_a       : H_{(n,m)} \longrightarrow  H_{(n,m+1)}
\end{equation*}
are defined by
\begin{align*}
s_\alpha x & = e_\alpha \otimes P_\alpha x, \qquad \text{ for } x \in H_{(0,0)}(=\A), \\
s_\alpha(e_{[\gamma]} \otimes P_{[\gamma]} x)
& =
{\begin{cases}
e_{[\alpha\gamma]} \otimes P_{[\alpha\gamma]} x & 
\text{ if } \alpha \gamma \in L_{(n+1,m)},\\
0  & \text{ otherwise},
\end{cases}} \\
t_a x & = e_a \otimes P_a x, \qquad \text{ for } x \in H_{(0,0)}(=\A), \\
t_a(e_{[\gamma]} \otimes P_{[\gamma]} x)
& =
{\begin{cases}
e_{[a\gamma]} \otimes P_{[a\gamma]} x & \text{ if } a \gamma \in L_{(n,m+1
)},\\
0  & \text{ otherwise}.
\end{cases}}
\end{align*}
For $y  \in \A$
an operator 
$i_{F(\rho,\eta)}(y)$ 
on $F(\rho,\eta)$:
\begin{equation*}
i_{F(\rho,\eta)}(y): H_{(n,m)}\longrightarrow H_{(n,m)}
\end{equation*}
is defined by
\begin{align*}
i_{F(\rho,\eta)}(y)x 
&=yx \qquad \text{ for } x \in H_{(0,0)}(=\A), \\
i_{F(\rho,\eta)}(y)( e_{[\gamma]} \otimes P_{[\gamma]} x)
 &= \phi(y) ( e_{[\gamma]} \otimes P_{[\gamma]} x)
( =e_{[\gamma]} \otimes \xi_\gamma(y) x).
\end{align*}
Define the Cuntz-Toeplitz $C^*$-algebra for $\CTDS$
$$
{\cal T}_{(\rho,\eta)}^\kappa
= C^*( s_\alpha, t_a, i_{F(\rho,\eta)}(y) \mid 
\alpha \in \Sigma^\rho, a \in \Sigma^\eta, y \in \A)
$$
as the $C^*$-algebra on $F_{\rho,\eta}$
generated by  
$
s_\alpha, t_a, i_{F(\rho,\eta)}(y)
$
for
$ 
\alpha \in \Sigma^\rho, a \in \Sigma^\eta, y \in \A.
$
\begin{lem}
For $\alpha \in \Sigma^\rho, a \in \Sigma^\eta$,
we have
\begin{enumerate}
\renewcommand{\labelenumi}{(\roman{enumi})}
\item
$
s_\alpha^*(e_{[\gamma]} \otimes P_{[\gamma]} x)
=
{\begin{cases}
\phi(\rho_\alpha(1))(e_{[\gamma']} \otimes P_{[\gamma']} x) 
& \text{ if } \gamma\approx \alpha \gamma',\\
0 & \text{otherwise.}
\end{cases}}
$
\item
$
t_a^*(e_{[\gamma]} \otimes P_{[\gamma]} x)
=
{\begin{cases}
\phi(\eta_a(1))(e_{[\gamma']} \otimes P_{[\gamma']} x) & \text{ if } \gamma\approx a \gamma',\\
0 & \text{otherwise.}
\end{cases}}
$
\end{enumerate}
\end{lem}
\begin{pf}
(i)
Suppose that
$\gamma\approx \alpha \gamma'$.
\begin{align*}
\langle s_\alpha^*(e_{[\gamma]} \otimes P_{[\gamma]} x) \mid 
e_{[\gamma']} \otimes P_{[\gamma']} x' \rangle 
& = \langle e_{[\gamma]} \otimes P_{[\gamma]} x \mid 
e_{[\alpha\gamma']} \otimes P_{[\alpha\gamma']} x' \rangle \\ 
& =
{\begin{cases}
x^* P_{[\alpha \gamma']} x & \text{ if } \gamma\approx \alpha \gamma',\\
0 & \text{otherwise.} 
\end{cases}}
\end{align*}
On the other hand,
$$
\phi(\rho_\alpha(1))(e_{[\gamma']} \otimes P_{[\gamma']} x)
= e_{[\gamma']} \otimes P_{[\alpha\gamma']} P_{\gamma'} x 
= e_{[\gamma']} \otimes P_{[\alpha\gamma']} x 
$$
so that
\begin{equation*}
\langle \phi(\rho_\alpha(1))(e_{[\gamma']} \otimes P_{[\gamma']} x) 
\mid 
e_{[\gamma']} \otimes P_{[\gamma']} x' \rangle 
=  x^* P_{[\alpha \gamma']} x'.
\end{equation*}
Hence we obtain the desired equality.
Similarly we see (ii).
\end{pf}

\begin{lem}
For $\alpha \in \Sigma^\rho, a \in \Sigma^\eta$,
we have
\begin{enumerate}
\renewcommand{\labelenumi}{(\roman{enumi})}
\item
$ s_\alpha^* s_\alpha = \phi(\rho_\alpha(1))$ 
and
$$
s_\alpha s_\alpha^*(e_{[\gamma]} \otimes P_{[\gamma]} x)
=
{\begin{cases}
e_{[\gamma]} \otimes P_{[\gamma]} x) 
  & \text{ if } \gamma\approx \alpha \gamma' \text{ for some } \gamma',\\
0 & \text{otherwise.}
\end{cases}}
$$
\item
$ t_a^* t_a = \phi(\eta_a(1))$
and
$$
t_a t_a^*(e_{[\gamma]} \otimes P_{[\gamma]} x)
=
{\begin{cases}
e_{[\gamma]} \otimes P_{[\gamma]} x) 
  & \text{ if } \gamma\approx a \gamma' \text{ for some } \gamma',\\
0 & \text{otherwise.}
\end{cases}}
$$
\end{enumerate}
\end{lem}
\begin{pf}
(i) 
It follows that for $\gamma \in L(n,m)$
\begin{equation*}
s_\alpha^* s_\alpha (e_{[\gamma]} \otimes P_{[\gamma]} x) 
 = {\begin{cases}
\phi(\rho_\alpha(1))(e_{[\gamma]} \otimes P_{[\gamma]} x) 
  & \text{ if }  \alpha \gamma \in L_{(n+1,m)},\\
0 & \text{otherwise.}
\end{cases}} 
\end{equation*}
On the other hand,
$$
\phi(\rho_\alpha(1))(e_{[\gamma]} \otimes P_{[\gamma]} x)
= 
e_{[\gamma]} \otimes P_{[\alpha \gamma]}P_{[\gamma]} x 
= 
e_{[\gamma]} \otimes P_{[\alpha \gamma]} x.
$$
Hence
$
\phi(\rho_\alpha(1))(e_{[\gamma]} \otimes P_{[\gamma]} x) =0
$
if
$\alpha \gamma \not\in L_{(n+1,m)}$.
Therefore we have
$$
s_\alpha^* s_\alpha (e_{[\gamma]} \otimes P_{[\gamma]} x) =
\phi(\rho_\alpha(1))(e_{[\gamma]} \otimes P_{[\gamma]} x)
$$
and hence
$
s_\alpha^* s_\alpha =\phi(\rho_\alpha(1)).
$

The equality
\begin{equation*}
s_\alpha s_\alpha^* (e_{[\gamma]} \otimes P_{[\gamma]} x)
 =
{\begin{cases}
e_{[\gamma]} \otimes P_{[\gamma]} x  
  & \text{ if }  \gamma \approx \alpha \gamma' \text{ for some } \gamma',\\
0 & \text{otherwise}
\end{cases}} 
\end{equation*}
is direct.

(ii) The assertion is similar to (i).
\end{pf}
\begin{lem} \hspace{5cm}
\begin{enumerate}
\renewcommand{\labelenumi}{(\roman{enumi})}
\item
$1 - \sum_{\alpha \in \Sigma^\rho}s_\alpha s_\alpha^* =$
the projection onto the subspace spanned by the vectors
$ e_{[\gamma]} \otimes P_{[\gamma]} x$
for $\ \gamma \in \cup_{m=0}^\infty L_{(0,m)}, x \in \A. $
\item
$1 - \sum_{a \in \Sigma^\eta}t_a t_a^* =$
the projection onto the subspace spanned by the vectors
$e_{[\gamma]} \otimes P_{[\gamma]} x $
for $\gamma \in \cup_{n=0}^\infty L_{(n,0)}, x \in \A.$
\end{enumerate}
\end{lem}
\begin{lem}
For $\alpha \in \Sigma^\rho, a \in \Sigma^\eta$ 
and 
$x \in \A$,
we have
\begin{enumerate}
\renewcommand{\labelenumi}{(\roman{enumi})}
\item
$ s_\alpha^* x s_\alpha = \phi(\rho_\alpha(x))$.
\item
$ t_a^* x  t_a = \phi(\eta_a(x))$.
\end{enumerate}
\end{lem}
\begin{pf}
For $y \in \A$, we have

(i)
\begin{align*}
s_\alpha^* x s_\alpha ( e_{[\gamma]}\otimes P_\gamma y) 
& 
= s_\alpha^*( e_{[\alpha\gamma]}\otimes P_{\alpha\gamma} y \xi_{\alpha \gamma}(x)) \\
& 
= e_{[\gamma]}\otimes P_{\gamma} y \xi_{\gamma}(\rho_\alpha(x))) \\
& 
= 
\phi(\rho_\alpha(x))(e_{[\gamma]}\otimes P_{\gamma} y ). 
\end{align*}

(ii)

\begin{align*}
t_a^* x t_a ( e_{[\gamma]}\otimes P_\gamma y) 
& 
= t_a^*( e_{[a\gamma]}\otimes P_{a\gamma} y \xi_{\alpha \gamma}(x)) \\
& 
= e_{[\gamma]}\otimes P_{\gamma} y \xi_{\gamma}(\eta_a(x))) \\
& 
= 
\phi(\eta_a(x))(e_{[\gamma]}\otimes P_{\gamma} y ). 
\end{align*}
\end{pf}

\begin{lem}
For $\alpha, \beta \in \Sigma^\rho$, $a,b \in \Sigma^\eta$
we have
\begin{equation}
s_\alpha t_b = t_a s_\beta
\qquad \text{if }
\kappa(\alpha,b) = (a,\beta).
\end{equation} 
\end{lem}
\begin{pf}
For $\gamma \in L_{(n,m)}$, 
suppose that
$\alpha b \gamma, a \beta \gamma \in L_{(n+1,m+1)}$.
It follows that
\begin{align*}
s_\alpha t_b (e_{[\gamma]}\otimes P_{\gamma} x )
& = e_{[\alpha b \gamma]}\otimes P_{\alpha b \gamma} y ),\\
t_a s_\beta (e_{[\gamma]}\otimes P_{\gamma} x) 
& =(e_{[a \beta \gamma]}\otimes P_{a \beta \gamma} x ).
\end{align*}
Since
$\kappa(\alpha,b) = (a,\beta)$,
the condition
$\alpha b \gamma \in L_{(n+1,m+1)}$
is equivalent to the condition
$a \beta \gamma \in L_{(n+1,m+1)}$.
We then have
$[\alpha b \gamma] = [a \beta \gamma]$
and
$P_{\alpha b \gamma} = P_{a \beta \gamma}.
$
\end{pf}

Let
${\cal I}_{(\rho,\eta)}$ be the ideal of
${\cal T}_{(\rho,\eta)}^\kappa$ generated by
the projections:
$
1 - \sum_{\alpha \in \Sigma^\rho}s_\alpha s_\alpha^*,
$
and
$1 - \sum_{a \in \Sigma^\eta}t_a t_a^*.
$
Let
$\OHRE$
be
the quatient $C^*$-algebra
\begin{equation*}
\OHRE = {\cal T}_{(\rho,\eta)}^\kappa /{\cal I}_{(\rho,\eta)}.
\end{equation*}
Let
$\pi_{(\rho,\eta)}: {\cal T}_{(\rho,\eta)}^\kappa \longrightarrow \OHRE$
be the quatient map.
Put
\begin{equation*}
\widehat{S}_\alpha = \pi_{(\rho,\eta)}(s_\alpha),\quad
\widehat{T}_a = \pi_{(\rho,\eta)}(t_a),\quad
\hat{i}(x) = \pi_{(\rho,\eta)}(i_{(F_{(\rho,\eta )}})(x)
\end{equation*}
for $\alpha \in \Sigma^\rho, a \in \Sigma^\eta$ and $x \in \A$. 
By the above discussions, the following relations hold:
\begin{align*}
\sum_{\beta \in \Sigma^\rho} \widehat{S}_{\beta} \widehat{S}_{\beta}^* & =1,
\qquad
\hat{i}(x)\widehat{S}_\alpha \widehat{S}_\alpha^* 
 =  \widehat{S}_\alpha \widehat{S}_\alpha^* \hat{i}(x),
\qquad
\widehat{S}_\alpha^* \hat{i}(x) \widehat{S}_\alpha 
 = \hat{i}(\rho_\alpha(x)),\\
\sum_{b \in \Sigma^\eta} \widehat{T}_{b} \widehat{T}_{b}^* & =1,
\qquad
\hat{i}(x) \widehat{T}_a \widehat{T}_a^*  
=  \widehat{T}_a \widehat{T}_a^* \hat{i}(x),
\qquad
\widehat{T}_a^* \hat{i}(x) \widehat{T}_a  = \hat{i}(\eta_a(x)),\\  
\widehat{S}_\alpha \widehat{T}_b&  = \widehat{T}_a \widehat{S}_\beta 
 \quad \text{ if }
\quad
\kappa(\alpha,b) = (a,\beta)
 \end{align*}
for all $x \in {\cal A}$ and $\alpha \in \Sigma^\rho, a \in \Sigma^\eta.$
Therefore we have 
\begin{prop}
Suppose that $\CTDS$ satisfies condition (I).
Then the algebra
$\OHRE$ is canonically isomorphic to 
the $C^*$-algebra $\ORE$
through the correspondences:
\begin{equation*}
S_{\alpha} \longrightarrow 
\widehat{S}_\alpha,
\qquad
T_a \longrightarrow  
\widehat{T}_a,
\qquad
x \longrightarrow \hat{i}(x)
\end{equation*}
for   $\alpha \in \Sigma^\rho, a \in \Sigma^\eta$ and 
$x \in {\cal A}$.
\end{prop}

\section{K-Theory Machinery}

In this section, we will study K-theory groups 
$K_*(\ORE)$ for the $C^*$-algebra $\ORE$.
We fix a $C^*$-textile dynamical system
$\CTDS$. 
We define  two actions
$$
\hat{\rho} :
{\Bbb T} 
\longrightarrow
\Aut(\ORE), \quad
 \hat{\eta} :
{\Bbb T} 
\longrightarrow
\Aut(\ORE)
$$
of the circle group ${\Bbb T} = \{ z \in {\Bbb C} \mid |z| = 1 \}$
to  $\ORE$ 
by setting
\begin{equation*}
\hat{\rho}_z = \widehat{\kappa}_{(z,1)},\qquad
\hat{\eta}_w = \widehat{\kappa}_{(1,w)}, \qquad z, w \in {\Bbb T}.  
\end{equation*}
They satisfy
\begin{equation*}
\hat{\rho}_z \circ \hat{\eta}_w =
\hat{\eta}_w \circ \hat{\rho}_z = \widehat{\kappa}_{(z,w)}, \qquad z, w \in {\Bbb T}.  
\end{equation*}

Set the fixed point algebras
\begin{align*}
(\ORE)^{\hat{\rho}}
& = \{ x \in \ORE \mid \hat{\rho}_z(x) = x \text{ for all } z \in {\Bbb T} \}, \\
(\ORE)^{\hat{\eta}}
& = \{ x \in \ORE \mid \hat{\eta}_z(x) = x \text{ for all } z \in {\Bbb T} \}.
\end{align*}
For 
$x \in (\ORE)^{\hat{\rho}}$,
define the constant function
$
\widehat{x} \in L^1({\Bbb T},\ORE)\subset\ORE\times_{\hat{\rho}}{\Bbb T}    
$
by setting
$
\widehat{x}(z) = x, z \in {\Bbb T}.
$
Put
$p_0 = \widehat{1}$.
By \cite{Rosen},
the algebra
$(\ORE)^{\hat{\rho}}$
is canonically isomorphic to
$p_0(\ORE\times_{\hat{\rho}}{\Bbb T})p_0$
through the map
\begin{equation*}
j_{\rho}: x \in (\ORE)^{\hat{\rho}}
\longrightarrow
\widehat{x} \in p_0(\ORE\times_{\hat{\rho}}{\Bbb T} )p_0   
\end{equation*}
which induces an isomorphism 
\begin{equation}
j_{\rho_*}: K_i((\ORE)^{\hat{\rho}}) 
\longrightarrow
K_i( p_0(\ORE\times_{\hat{\rho}}{\Bbb T} )p_0),
\qquad
i=0,1   
\end{equation}
on their K-groups.

\begin{lem}\hspace{5cm}
\begin{enumerate}
\renewcommand{\labelenumi}{(\roman{enumi})}
\item
There exists an isometry  
$v \in M((\ORE\times_{\hat{\rho}}{\Bbb T})\otimes {\cal K})$
such that
$v v^* = p_0 \otimes 1, v^* v =1$.
\item
$\ORE\times_{\hat{\rho}}{\Bbb T}$ is stably isomorphic to 
$(\ORE)^{\hat{\rho}}$,
and similarly 
$\ORE\times_{\hat{\eta}}{\Bbb T}$ is stably isomorphic to 
$(\ORE)^{\hat{\eta}}$.
\item
The inclusion
$\iota_{\hat{\rho}}:
p_0(\ORE\times_{\hat{\rho}}{\Bbb T})p_0
\hookrightarrow
\ORE\times_{\hat{\rho}}{\Bbb T}
$
induces an isomorphism 
\begin{equation}
\iota_{\hat{\rho}*}:
K_0(p_0(\ORE\times_{\hat{\rho}}{\Bbb T})p_0)
\cong
K_0(\ORE\times_{\hat{\rho}}{\Bbb T})
\end{equation}
on their K-groups.
\end{enumerate}
\end{lem}
\begin{pf}
(i)
We will prove that 
$p_0$ is a full projection in
$\ORE\times_{\hat{\rho}}{\Bbb T}$.
Suppose that there exists 
an irreducible representation
$\pi$ of $\ORE\times_{\hat{\rho}}{\Bbb T}$
such that 
$\pi(p_0) =0$.
Denote by $*$ the $\hat{\rho}$-twisted convolution product
in $ L^1({\Bbb T},\ORE)$
(= the product in the algebra $\ORE\times_{\hat{\rho}}{\Bbb T}$).
For $Y \in \ORE$,
put
$\widehat{Y}(z) = Y$ for $z \in {\Bbb T}$.
The equality
$\widehat{Y} * p_0 = \widehat{Y}$
implies
$\widehat{Y} \in \ker(\pi)$.
For $Y, Z \in \ORE$
by using the equality
$\widehat{Z}^*(z) =\hat{\rho}_z(Z^*)$,
we have
$(\widehat{Y} * \widehat{Z}^*)(z) = Y \hat{\rho}_z(Z^*)$.
For $\mu \in B_k(\Lambda_\rho)$,
we have
\begin{equation*}
(\widehat{Y S_\mu}*\widehat{S}_\mu^*)(z) = z^{-k}Y S_\mu S_\mu^*
\end{equation*}
and hence  
\begin{equation*}
(\sum_{\mu \in B_k(\Lambda_\rho)}\widehat{Y S_\mu}*\widehat{S}_\mu^*)(z) 
= z^{-k}Y. 
\end{equation*}
As
$\widehat{Y S_\mu}, \widehat{S}_\mu^* \in \ker(\pi)$,
the function
$z \in {\Bbb T} \longrightarrow z^{-k}Y \in \ORE$
belongs to $\ker(\pi)$ for $k=0,1,2,\dots$.
Let
$E_i^l, i = 1,2,\dots,m(k)$
be the minimal projections
in the commutative $C^*$-algebra
$C^*(\rho_\mu(1)\mid \mu \in B_k(\Lambda_\rho))$
generated by
the projections
$\rho_\mu(1), \mu \in B_k(\Lambda_\rho)$.
Hence
$\sum_{i=1}^{m(k)}E_i^k =1$
and for $i=1,\dots,m(k)$,
there exists $\mu(i) \in B_k(\Lambda_\rho)$
such that
$E_i^k \le S_{\mu(i)}^*S_{\mu(i)}$.
Since
 for $ Y \in \ORE$,
\begin{equation*}
(\widehat{Y E_i^k S_\mu^*}*\widehat{S^*_\mu}^*)(z)
  = z^{k}Y E_i^k S_\mu^* S_\mu= z^{k}Y E_i^k, 
\end{equation*}
we have
\begin{equation*}
(\sum_{i=1}^{m(k)}\widehat{Y E_i^k S_\mu^*}*\widehat{S^*_\mu}^*)(z)
  = z^{k}Y. 
\end{equation*}
As
$\widehat{Y E_i^k S_\mu^*}, \widehat{S^*}_\mu^* \in \ker(\pi)$,
the function
$z \in {\Bbb T} \longrightarrow z^k Y \in \ORE$
belongs to $\ker(\pi)$ for $k=0,1,2,\dots$.
Therefore we know that
the functions 
$z \in {\Bbb T} \longrightarrow z^k Y \in \ORE$
belongs to $\ker(\pi)$ for all $k\in {\Bbb Z}$.
In particular, for $Y=1$
the functions
$z \in {\Bbb T} \longrightarrow z^k  \in \ORE$
belongs to $\ker(\pi)$ for all $k\in {\Bbb Z}$
so that
$C({\Bbb T})$ is contained in $\ker(\pi)$ .
Take an approximate identity
$\varphi_n \in C({\Bbb T}), n \in {\Bbb N}$
for the usual convolution product in $L^1({\Bbb T})$.
Then for $X \in L^1({\Bbb T},\ORE)$,
one has
$\| X * \varphi_n - X \|_1 \rightarrow 0$ 
as 
$n\rightarrow \infty$.
Since $ X * \varphi_n \in \ker(\pi)$,
one has 
$X \in \ker(\pi)$.
Hence we have
$ L^1({\Bbb T},\ORE) \subset \ker(\pi)$
so that 
$\ker(\pi) = \ORE\times_{\hat{\rho}}{\Bbb T}$.
Therefore 
$p_0 \in \ORE\times_{\hat{\rho}}{\Bbb T}$
is a full projection of $\ORE\times_{\hat{\rho}}{\Bbb T}$,
By \cite[Corollary 2.6]{Bro},
there exists 
$v \in M((\ORE\times_{\hat{\rho}}{\Bbb T})\otimes {\cal K})$
such that
$v v^* = p_0 \otimes 1, v^* v =1$.

(ii)
As 
$
Ad(v^*):  
x \in p_0(\ORE\times_{\hat{\rho}}{\Bbb T})p_0\otimes {\cal K}
\rightarrow
v^* x v \in \ORE\times_{\hat{\rho}}{\Bbb T}\otimes {\cal K}
$
is an isomorphism
and
$(\ORE)^{\hat{\rho}}$
is isomorphic 
to
$
p_0(\ORE\times_{\hat{\rho}}{\Bbb T})p_0
$
through 
$j_\rho$,
we have
$(\ORE)^{\hat{\rho}}$
is stably isomorphic 
to
$
\ORE\times_{\hat{\rho}}{\Bbb T}.
$

(iii)
Let
$v \in M((\ORE\times_{\hat{\rho}}{\Bbb T})\otimes {\cal K})$
be the isometry as above such that
$v v^* = p_0 \otimes 1, v^* v =1$.
For a projection 
$q \in  
p_0(\ORE\times_{\hat{\rho}}{\Bbb T})p_0\otimes {\cal K}$,
we have 
$ [v^*qv] = [q] = \iota_{\hat{\rho}*}([q])$
and hence
$
\iota_{\hat{\rho}*}= Ad(v^*)_*:
K_0(p_0(\ORE\times_{\hat{\rho}}{\Bbb T})p_0)
\cong
K_0(\ORE\times_{\hat{\rho}}{\Bbb T})
$
is an isomorphism.
\end{pf}
Thanks to the lemma above,
$
Ad(v^*):  
x \in p_0(\ORE\times_{\hat{\rho}}{\Bbb T})p_0\otimes {\cal K}
\rightarrow
v^* x v \in \ORE\times_{\hat{\rho}}{\Bbb T}\otimes {\cal K}
$
induces isomorphisms
\begin{equation}
 Ad(v^*)_*:
K_i(p_0(\ORE\times_{\hat{\rho}}{\Bbb T})p_0)
\longrightarrow
K_i(\ORE\times_{\hat{\rho}}{\Bbb T}),
\qquad i=0,1.
\end{equation}

Let $\hat{\hat{\rho}}$ 
be the automorphism on
$\ORE\times_{\hat{\rho}}{\Bbb T} $
for the positive generator of
${\Bbb Z}$ for
the dual action of 
$\ORE\times_{\hat{\rho}}{\Bbb T} $.
By (6.1) and (6.3), we may define an isomorphism 
\begin{equation}
\beta_{\rho,i} = j_{\rho*}^{-1}\circ Ad(v^*)_*^{-1}\circ \hat{\hat{\rho}}_* 
\circ Ad(v^*)_* \circ j_{\rho*}: 
K_i((\ORE)^{\hat{\rho}}) \longrightarrow 
K_i((\ORE)^{\hat{\rho}}), \quad i=0,1
\end{equation}
so that the diagram is commutative:
\begin{equation*}
\begin{CD}
K_i(\ORE\times_{\hat{\rho}}{\Bbb T}) @>{\hat{\hat{\rho}}_*}>>
K_i(\ORE\times_{\hat{\rho}}{\Bbb T})\\
@AA{Ad(v^*)_*}A 
@AA{Ad(v^*)_*}A \\
K_i(p_0(\ORE\times_{\hat{\rho}}{\Bbb T})p_0) @.
K_i(p_0(\ORE\times_{\hat{\rho}}{\Bbb T})p_0)\\
@AA{j_{\rho*}}A 
@AA{j_{\rho*}}A \\
K_i((\ORE)^{\hat{\rho}})
@>{\beta_{\rho,i}}>>
K_i((\ORE)^{\hat{\rho}})
\end{CD}
\end{equation*}

By \cite{Pim} (cf. \cite{KPW}), one has the six term
exact sequence of K-theory:
\begin{equation*}
\begin{CD}
 K_0(\ORE\times_{\hat{\rho}}{\Bbb T}) @>\id - \hat{\hat{\rho}}_{*}>> 
 K_0(\ORE\times_{\hat{\rho}}{\Bbb T}) @>\iota_{*}>> 
 K_0((\ORE\times_{\hat{\rho}}{\Bbb T})\times_{\hat{\hat{\rho}}}{\Bbb Z}) \\
@A{\delta}AA @. @V{\exp}VV \\
 K_1((\ORE\times_{\hat{\rho}}{\Bbb T})\times_{\hat{\hat{\rho}}}{\Bbb Z})
@<<\iota_{*}<
 K_1(\ORE\times_{\hat{\rho}}{\Bbb T}) @<<\id - \hat{\hat{\rho}}_{*}< 
 K_1(\ORE\times_{\hat{\rho}}{\Bbb T})  
\end{CD}
\end{equation*}
Since
$
(\ORE\times_{\hat{\rho}}{\Bbb T})\times_{\hat{\hat{\rho}}}{\Bbb Z} 
\cong \ORE\otimes{\cal K}
$
and
$
 K_*(\ORE\times_{\hat{\rho}}{\Bbb T})\cong K_*((\ORE)^{\hat{\rho}}),
$
one has 
\begin{lem}
The following six term
exact sequence of K-theory holds:
\begin{equation*}
\begin{CD}
 K_0((\ORE)^{\hat{\rho}}) @>\id - \beta_{\rho,0}>> 
 K_0((\ORE)^{\hat{\rho}}) @>\iota_{*}>> 
 K_0(\ORE) \\
@A{\delta}AA @. @V{\exp}VV \\
 K_1(\ORE)
@<<\iota_{*}<
 K_1((\ORE)^{\hat{\rho}}) @<<\id - \beta_{\rho,1}< 
 K_1((\ORE)^{\hat{\rho}}) . 
\end{CD}
\end{equation*}
Hence there exist short exact sequences for $i=0,1$:
\begin{align*}
0 &
\longrightarrow 
\Coker(\id  - \beta_{\rho,i}) \text{ in } K_i((\ORE)^{\hat{\rho}})\\
& \longrightarrow 
K_i(\ORE)\\
& \longrightarrow
\Ker(\id - \beta_{\rho,i+1}) \text{ in } K_{i+1}((\ORE)^{\hat{\rho}}) \\
& \longrightarrow 
0.
\end{align*}
\end{lem}
We will then study the following groups that appear in the above sequences
$$  
\Coker(\id  - \beta_{\rho,i}) \text{ in } K_i((\ORE)^{\hat{\rho}}),
\qquad
\Ker(\id - \beta_{\rho,i+1}) \text{ in } K_{i+1}((\ORE)^{\hat{\rho}})
$$
for
$i=0,1$.
The action $\hat{\eta}$ acts on the subalgebra 
$(\ORE)^{\hat{\rho}}$,
which we still denote by $\hat{\eta}$.
Then 
the fixed point algebra 
$((\ORE)^{\hat{\rho}})^{\hat{\eta}}$
of
$(\ORE)^{\hat{\rho}}$
under the action
$\hat{\eta}$
coincides with
$\FRE$.
The above discussions 
for the action $\hat{\rho}: {\Bbb T}\longrightarrow \ORE$
works for the action
$\hat{\eta}: {\Bbb T}\longrightarrow (\ORE)^{\hat{\rho}}$
as in the following way.
For 
$y \in ((\ORE)^{\hat{\rho}})^{\hat{\eta}}$,
define the constant function
$
\widehat{y} \in L^1({\Bbb T},(\ORE)^{\hat{\rho}})
\subset(\ORE)^{\hat{\rho}}\times_{\hat{\eta}}{\Bbb T}    
$
by setting
$
\widehat{y}(z) = y, z \in {\Bbb T}.
$
Putting
$q_0 = \widehat{1}$,
the algebra
$((\ORE)^{\hat{\rho}})^{\hat{\eta}}$
is canonically isomorphic to
$q_0((\ORE)^{\hat{\rho}}\times_{\hat{\eta}}{\Bbb T})q_0$
through the map
\begin{equation*}
j^{\rho}_{\eta}: y \in ((\ORE)^{\hat{\rho}})^{\hat{\eta}}
\longrightarrow
\hat{y} \in q_0((\ORE)^{\hat{\rho}}\times_{\hat{\eta}}{\Bbb T} )q_0   
\end{equation*}
which induces an isomorphism 
\begin{equation}
j^{\rho}_{\eta*}: K_i(((\ORE)^{\hat{\rho}})^{\hat{\eta}})
\longrightarrow
K_i(q_0((\ORE)^{\hat{\rho}}\times_{\hat{\eta}}{\Bbb T} )q_0)   
\end{equation}
on their K-groups.
Similarly to Lemma 6.1, we have
\begin{lem}\hspace{5cm}
\begin{enumerate}
\renewcommand{\labelenumi}{(\roman{enumi})}
\item
There exists an isometry  
$u \in M(((\ORE)^{\hat{\rho}}\times_{\hat{\eta}}{\Bbb T})\otimes {\cal K})$
such that
$u u^* = q_0 \otimes 1, u^* u =1$.
\item
$(\ORE)^{\hat{\rho}}\times_{\hat{\eta}}{\Bbb T}$ 
is stably isomorphic to $((\ORE)^{\hat{\rho}})^{\hat{\eta}}$.
\item
The inclusion
$\iota^{\hat{\rho}}_{\hat{\eta}}:
q_0((\ORE)^{\hat{\rho}}\times_{\hat{\eta}}{\Bbb T})q_0  
(= ((\ORE)^{\hat{\rho}})^{\hat{\eta}} = \FRE)
\hookrightarrow
(\ORE)^{\hat{\rho}}\times_{\hat{\eta}}{\Bbb T}
$
induces an isomorphism 
\begin{equation}
\iota^{\hat{\rho}}_{\hat{\eta}*}:
K_0(q_0((\ORE)^{\hat{\rho}}\times_{\hat{\eta}}{\Bbb T})q_0)
\cong
K_0((\ORE)^{\hat{\rho}}\times_{\hat{\eta}}{\Bbb T})
\end{equation}
on their K-groups.
\end{enumerate}
\end{lem}

The isomorphism
$$
Ad(u^*):
y \in q_0((\ORE)^{\hat{\rho}}\times_{\hat{\eta}}{\Bbb T})q_0
\longrightarrow
u^* y u \in 
(\ORE)^{\hat{\rho}}\times_{\hat{\eta}}{\Bbb T}
$$
induces isomorphismss
\begin{equation}
Ad(u^*)_*:
K_i(q_0((\ORE)^{\hat{\rho}}\times_{\hat{\eta}}{\Bbb T})q_0)
\cong
K_i((\ORE)^{\hat{\rho}}\times_{\hat{\eta}}{\Bbb T}),
\qquad i=0,1.
\end{equation}
Let $\hat{\hat{\eta}}_\rho$ 
be the automorphism of the positive generator of ${\Bbb Z}$ 
for the dual action of 
$(\ORE)^{\hat{\rho}}\times_{\hat{\eta}}{\Bbb T} $.
Define an isomorphism
\begin{equation}
\gamma_{\eta,i}  
= j_{\eta*}^{\rho -1}\circ Ad(u^*)_*^{-1} \circ \hat{\hat{\eta}}_{\rho *} 
\circ Ad(u^*)_* \circ j^\rho_{\eta*}: 
K_i(\FRE ) \longrightarrow 
K_i(\FRE )
\end{equation}
such that the diagram is commutative for $i=0,1$:
\begin{equation*}
\begin{CD}
K_i((\ORE)^{\hat{\rho}}\times_{\hat{\eta}}{\Bbb T}) @>{\hat{\hat{\eta}}_{\rho *}}>>
K_i((\ORE)^{\hat{\rho}}\times_{\hat{\eta}}{\Bbb T})\\
@AA{Ad(u^*)_*}A 
@AA{Ad(u^*)_*}A \\
K_i(q_0((\ORE)^{\hat{\rho}}\times_{\hat{\eta}}{\Bbb T})q_0) @.
K_i(q_0((\ORE)^{\hat{\rho}}\times_{\hat{\eta}}{\Bbb T})q_0)\\
@AA{j^\rho_{\eta*}}A 
@AA{j^\rho_{\eta*}}A \\
K_i(((\ORE)^{\hat{\rho}})^{\hat{\eta}})
@.
K_i(((\ORE)^{\hat{\rho}})^{\hat{\eta}})\\
@|  @| \\
K_i(\FRE)
@>{\gamma_{\eta,i}}>>
K_i(\FRE)
\end{CD}
\end{equation*}
We similarly define
an endomorphism
$\gamma_{\rho,i}: K_i(\FRE) \longrightarrow K_i(\FRE)$.


Under the equality  $((\ORE)^{\hat{\rho}})^{\hat{\eta}} = \FRE$,
we have the following lemma which is similar to Lemma 6.2
\begin{lem}
The following six term
exact sequence of K-theory holds:
\begin{equation*}
\begin{CD}
 K_0(\FRE) @>\id - \gamma_{\eta,0} >> 
 K_0(\FRE) @>\iota_{*}>> 
 K_0((\ORE)^{\hat{\rho}}) \\
@A{\delta}AA @. @V{\exp}VV \\
 K_1((\ORE)^{\hat{\rho}})
@<<\iota_{*}<
 K_1(\FRE) @<<\id - \gamma_{\eta,1} < 
 K_1(\FRE)  
\end{CD}
\end{equation*}
In particular,
if $K_1(\FRE) =0$, we have
\begin{align}
K_0((\ORE)^{\hat{\rho}})
& = \Coker(\id - \gamma_{\eta,0}) \quad \text{ in } K_0(\FRE),\\
K_1((\ORE)^{\hat{\rho}})
& = \Ker(\id - \gamma_{\eta,0}) \quad \text{ in } K_0(\FRE).
\end{align}
\end{lem}

The following lemmas hold.
\begin{lem}
For a projection $ q \in M_n((\ORE)^{\rho})$
and a partial isometry $S \in \ORE$ 
such that
\begin{equation*}
{\hat{\rho}}_z(S) = z S \quad \text{for } z \in {\Bbb T},
\qquad
q (S S^*\otimes 1_n) = (S S^*\otimes 1_n) q,
\end{equation*}
we have
\begin{equation*}
\beta_{\rho,0}^{-1} ([ (S S^* \otimes 1_n) q])  
= [ (S^*\otimes 1_n) q (S \otimes 1_n)] 
\quad \text{ in }
K_0((\ORE)^{\hat{\rho}}).
\end{equation*}
\end{lem}
\begin{pf}
As
$q$ commutes with
$S S^*\otimes 1_n$,
$p= (S^*\otimes 1_n) q (S\otimes 1_n)$ is a projection in 
$(\ORE)^{\hat{\rho}}$.
Since $ p \le S^* S\otimes 1_n$,
By a similar argument to the proof of \cite[Lemma 4.5]{MathScand98},
one sees that
$\beta_{\rho,0}([p]) = [(S\otimes 1_n) p (S^*\otimes 1_n)]$ 
in
$K_0((\ORE)^{\hat{\rho}})$.
\end{pf}

\begin{lem}\hspace{5cm}
\begin{enumerate} 
\renewcommand{\labelenumi}{(\roman{enumi})}
\item
For a projection $ q \in M_n(\FRE)$
and a partial isometry 
$T \in (\ORE)^{\hat{\rho}}$ 
such that
\begin{equation*}
\hat{\eta}_z(T) = z T \quad \text{for } z \in {\Bbb T},
\qquad
q (T T^*\otimes 1_n) = (T T^*\otimes 1_n) q,
\end{equation*}
we have
\begin{equation*}
\gamma_{\eta,0}^{-1} ([ (T T^*\otimes 1_n) q])  = [( T^*\otimes 1_n) q (T\otimes 1_n)] 
\quad \text{ in }
K_0(\FRE).
\end{equation*}
\item
For a projection $ q \in M_n(\FRE) $
and a partial isometry $S \in (\ORE)^{\hat{\eta}}$ 
such that
\begin{equation*}
\hat{\rho}_z(S) = z S \quad \text{for } z \in {\Bbb T},
\qquad
q (S S^*\otimes 1_n) = (S S^*\otimes 1_n) q,
\end{equation*}
we have
\begin{equation*}
\gamma_{\rho,0}^{-1} ([ (S S^*\otimes 1_n) q])  = [ (S^*\otimes 1_n) q (S\otimes 1_n)] 
\quad \text{ in }
K_0(\FRE).
\end{equation*}
\end{enumerate}
\end{lem}
Hence we have 
\begin{lem}
The diagram 
\begin{equation}
\begin{CD}
K_0(\FRE) @>{\id-\gamma_{\rho,0}}>>
K_0(\FRE) \\
@VV{\iota_*}V 
@VV{\iota_*}V \\
K_0((\ORE)^{\hat{\rho}}) @>{\id-\beta_{\rho,0}}>>
K_0((\ORE)^{\hat{\rho}})
\end{CD}
\end{equation}
is commutative.
\end{lem}
\begin{pf}
By \cite[Proposition 3.3]{PasJOT},
the map
$\iota_*:K_0(\FRE) \longrightarrow K_0((\ORE)^{\hat{\rho}})$
is induced by the natural inclusion
$\FRE(= ((\ORE)^{\hat{\rho}})^{\eta}) \hookrightarrow (\ORE)^{\hat{\rho}}$.
For an element
 $[q] \in K_0(\FRE)$
 one may assume that
 $q \in M_n(\FRE)$ for some $n \in {\Bbb N}$
 so that  one has
\begin{align*}
 \gamma_{\rho,0}^{-1}([q]) 
=& \sum_{\alpha \in \Sigma^\rho}[(S_\alpha S_\alpha^* \otimes 1_n)q ] \\
=&\sum_{\alpha \in \Sigma^\rho}
       [(S_\alpha^*\otimes 1_n) q (S_\alpha\otimes 1_n)] \\
=&\sum_{\alpha \in \Sigma^\rho}
       \beta_{\rho,0}^{-1}([ q (S_\alpha S_\alpha^*\otimes 1_n)])
 = \beta_{\rho,0}^{-1}([ q ])
\end{align*}
so that 
$\beta_{\rho,0} |_{K_0(\FRE)} = \gamma_{\rho,0}$.
\end{pf}

In the rest of this section,
we assume that $K_1(\FRE) =0$.
 The following lemma is crucial in our further discussions.
\begin{lem}
In the six term exact sequence in Lemma 6.4 with $K_1(\FRE) =0$, 
we have
the following commutative diagrams:
\begin{equation} \label{eq:key}
\begin{CD}
0 @. 0 \\
@VVV @VVV \\
K_1((\ORE)^{\hat{\rho}}) @>\id - \beta_{\rho,1}>> K_1((\ORE)^{\hat{\rho}}) \\
@V{\delta}VV   @V{\delta}VV \\
K_0(\FRE) @>\id - \gamma_{\rho,0}>> K_0(\FRE)\\
@V{\id - \gamma_{\eta,0}}VV  @V{\id - \gamma_{\eta,0}}VV \\
K_0(\FRE) @>\id - \gamma_{\rho,0}>> K_0(\FRE)\\  
@V{\iota_*}VV  @V{\iota_*}VV \\
K_0((\ORE)^{\hat{\rho}}) @>\id - \beta_{\rho,0} 
>> K_0((\ORE)^{\hat{\rho}}) \\
@VVV @VVV \\
0 @. 0   
\end{CD}
\end{equation}
\end{lem}
\begin{pf}
It is well-known that 
$\delta$-map is functorial (see \cite[Theorem 7.2.5]{WeOl}, \cite[p.266 (LX)]{Bla}).
Hence the diagram of the upper square
\begin{equation*}
\begin{CD}
K_1((\ORE)^{\hat{\rho}}) @>\id - \beta_{\rho,1}>> K_1((\ORE)^{\hat{\rho}}) \\
@V{\delta}VV   @V{\delta}VV \\
K_0(\FRE) @>\id - \gamma_{\rho,0}>> K_0(\FRE)
\end{CD}
\end{equation*}
is commutative.

Since 
$
\gamma_{\rho,0} \circ \gamma_{\eta,0} 
= \gamma_{\eta,0} \circ \gamma_{\rho,0}
$
the diagram of the middle square 
\begin{equation}
\begin{CD}
K_0(\FRE) @>{\id-\gamma_{\rho,0}}>>
K_0(\FRE) \\
@VV{\id-\gamma_{\eta,0}}V 
@VV{\id-\gamma_{\eta,0}}V \\
K_0(\FRE) @>{\id-\gamma_{\rho,0}}>>
K_0(\FRE)
\end{CD}
\end{equation}
is commutative.

The commutativity of the lower square
comes from the preceding lemma. 
\end{pf} 
\begin{lem}
Suppose that $K_1(\FRE) =0$.
The six term exact sequence in Lemma 6.2 with Lemma 6.8
goes to the following  commutative diagrams:
{\small
\begin{equation*} 
\begin{CD}
  @. 0 @. 0  @. @.\\
@. @VVV @VVV @. @.\\
 @. K_1((\ORE)^{\hat{\rho}})@>\id - \beta_{\rho,1}>> 
                           K_1((\ORE)^{\hat{\rho}}) @. @.\\
@. @V{\delta}VV   @V{\delta}VV @. @.\\
 @. K_0(\FRE) @>\id - \gamma_{\rho,0}>> K_0(\FRE)@. @.\\
@. @V{\id - \gamma_{\eta,0}}VV  @V{\id - \gamma_{\eta,0}}VV @. @.\\
 @. K_0(\FRE) @>\id - \gamma_{\rho,0}>> K_0(\FRE)@. @.\\  
@. @V{\iota_*}VV  @V{\iota_*}VV @. @.\\
@>\delta >> K_0((\ORE)^{\hat{\rho}}) @>\id - \beta_{\rho,0}>>
    K_0((\ORE)^{\hat{\rho}}) @>\iota_*>> K_0(\ORE) @>\exp>>  \\
 @. @VVV @VVV @. @.\\
  @. 0 @. 0   @. @.
\end{CD}
\end{equation*}}
{\small
\begin{equation*} 
\begin{CD}
  @. 0 @. 0  @. @.\\
@. @VVV @VVV @. @.\\
@>\exp>> K_1((\ORE)^{\hat{\rho}})@>\id - \beta_{\rho,1}>> 
         K_1((\ORE)^{\hat{\rho}}) @>\iota_*>> K_1(\ORE) @>\delta>>\\
@. @V{\delta}VV   @V{\delta}VV @. @.\\
 @. K_0(\FRE) @>\id - \gamma_{\rho,0}>> K_0(\FRE)@. @.\\
@. @V{\id - \gamma_{\eta,0}}VV  @V{\id - \gamma_{\eta,0}}VV @. @.\\
 @. K_0(\FRE) @>\id - \gamma_{\rho,0}>> K_0(\FRE)@. @.\\  
@. @V{\iota_*}VV  @V{\iota_*}VV @. @.\\
 @. K_0((\ORE)^{\hat{\rho}}) @>\id - \beta_{\rho,0}>>
    K_0((\ORE)^{\hat{\rho}}) @. @.  \\
 @. @VVV @VVV @. @.\\
  @. 0 @. 0   @. @.
\end{CD}
\end{equation*}}
\end{lem}
We will describe the K-theory groups
$K_*(\ORE)$ in terms of the kernels and cokernels of the homomorphisms
$\id - \gamma_{\rho,i}$ and $\id - \gamma_{\eta,i}$
on $K_0(\FRE)$.
Recall that there exist short exact sequences by Lemma 6.2:
\begin{enumerate}
\renewcommand{\labelenumi}{(\roman{enumi})}
\item
\begin{align*}
0 &
\longrightarrow 
\Coker(\id  - \beta_{\rho,0}) \text{ in } K_0((\ORE)^{\hat{\rho}})\\
& \longrightarrow 
K_0(\ORE)\\
& \longrightarrow
\Ker(\id - \beta_{\rho,1}) \text{ in } K_{1}((\ORE)^{\hat{\rho}}) \\
& \longrightarrow 
0.
\end{align*}
\item
\begin{align*}
0 &
\longrightarrow 
\Coker(\id  - \beta_{\rho,1}) \text{ in } K_1((\ORE)^{\hat{\rho}})\\
& \longrightarrow 
K_1(\ORE)\\
& \longrightarrow
\Ker(\id - \beta_{\rho,0}) \text{ in } K_{0}((\ORE)^{\hat{\rho}}) \\
& \longrightarrow 
0.
\end{align*}
\end{enumerate}

As
$
\gamma_{\eta,0} \circ \gamma_{\rho,0} = 
\gamma_{\rho,0} \circ \gamma_{\eta,0}
$
on 
$K_0(\FRE)$,
$\gamma_{\rho,0} $
and
$
\gamma_{\eta,0}
$
naturally act on 
$
\Coker(\id - \gamma_{\eta,0})=
K_0(\FRE)/(\id - \gamma_{\eta,0})K_0(\FRE)
$
and
$
\Coker(\id - \gamma_{\rho,0})=
K_0(\FRE)/(\id - \gamma_{\rho,0})K_0(\FRE)
$
as endomorphisms
respectively,
which we denote by 
$\bar{\gamma}_{\rho,0}$
and
$\bar{\gamma}_{\eta,0}$
respectively.

\begin{lem} \hspace{5cm}
\begin{enumerate}
\renewcommand{\labelenumi}{(\roman{enumi})}
\item For $K_0(\ORE)$, we have 
\begin{align*}
      &\Coker(\id - \beta_{\rho,0})\text{ in } K_0( (\ORE)^{\hat{\rho}}) \\
\cong &
 \Coker(\id - \bar{\gamma}_{\rho,0})
\text{ in }
K_0(\FRE)/(\id - \gamma_{\eta,0})K_0(\FRE) \\
\cong & 
K_0(\FRE)
/ 
((\id - \gamma_{\rho,0})K_0(\FRE) +
(\id - \gamma_{\eta,0})K_0(\FRE) )\\
\intertext{ and}
& \Ker(\id - \beta_{\rho,0})\text{ in } K_1( (\ORE)^{\hat{\rho}}) \\
\cong 
& \Ker(\id - \gamma_{\rho,0}) 
\text{ in } (\Ker(\id - \gamma_{\eta,0}) \text{ in } K_0(\FRE)) \\
\cong 
& \Ker(\id - \gamma_{\rho,0}) \cap \Ker(\id - \gamma_{\eta,0}) 
\text{ in } K_0(\FRE).
\end{align*}
\item For $K_1(\ORE)$, we have 
\begin{align*}
& \Coker(\id - \beta_{\rho,1})\text{ in } K_1( (\ORE)^{\hat{\rho}}) \\
\cong &
(\Ker(\id - \gamma_{\eta,0}) \quad \text{ in } K_0(\FRE))
/(\id - \gamma_{\rho,0})(\Ker(\id - \gamma_{\eta,0})\text{ in } K_0(\FRE))\\
\intertext{ and }
& \Ker(\id - \beta_{\rho,0})\text{ in } K_0( (\ORE)^{\hat{\rho}}) \\
\cong & 
\Ker(\id - \bar{\gamma}_{\rho,0}) 
\text{ in } (K_0(\FRE)/(\id - \gamma_{\eta,0})K_0(\FRE) ).
\end{align*}
\end{enumerate}
\end{lem}
\begin{pf}
(i) 
We will first prove the assertions for the group
$\Coker(\id - \beta_{\rho,0}) \text{ in } K_0((\ORE)^{\hat{\rho}})$.

In the diagram (\ref{eq:key}),
the exactness of the vertical arrows at $K_0(\FRE)$,
one sees that 
$\delta$ is injective and
$ \Im(\delta) = \Ker(\id - \gamma_\eta)$
so that we have
\begin{equation}
K_1((\ORE)^{\hat{\rho}})
\cong
\delta(K_1((\ORE)^{\hat{\rho}}))
\cong
\Ker(\id - \gamma_{\eta,0}) \text{ in } K_0(\FRE).
\end{equation}
By the commutativity in the upper square in the diagram (\ref{eq:key}),
one has
\begin{equation*}
\Ker(\id - \beta_{\rho,0}) \text{ in } K_1((\ORE)^{\hat{\rho}})
\cong 
\Ker(\id - \gamma_{\rho,0}) \text{ in }
(
\Ker(\id - \gamma_{\eta,0}) \text{ in } K_0(\FRE)).
\end{equation*}
Since
$\gamma_{\eta,0}$ 
commutes with 
$\gamma_{\rho,0}$ in $K_0(\FRE)$,
we have
\begin{align*}
& \Ker(\id - \gamma_{\rho,0}) \text{ in }
(
\Ker(\id - \gamma_{\eta,0}) \text{ in } K_0(\FRE)
) \\
\cong
&
\Ker(\id - \gamma_{\rho,0}) \cap
\Ker(\id - \gamma_{\eta,0}) \text{ in } K_0(\FRE).
\end{align*}

We will second prove the assertions for the group
$\Ker(\id - \beta_{\rho,1}) \text{ in } K_1((\ORE)^{\hat{\rho}})$.

In the diagram (\ref{eq:key}),
the exactness of the vertical arrows at $K_0(\FRE)$,
one sees that 
$\iota_*$ is surjective so that
\begin{align*}
K_0((\ORE)^{\hat{\rho}})
& \cong \iota_* (K_0(\FRE)) \\
& \cong
K_0(\FRE) / 
\Ker(\id - \gamma_{\eta,0}) \text{ in } K_0(\FRE).
\end{align*}
By the commutativity in the lower square in the diagram (\ref{eq:key}),
one has
\begin{align*}
 &\Coker(\id - \beta_{\rho,0}) \text{ in } K_1((\ORE)^{\hat{\rho}}) \\
\cong & 
\Coker(\id - \bar{\gamma}_{\rho,0}) \text{ in } 
(
\Coker(\id - \gamma_{\eta,0}) \text{ in } K_0(\FRE)
) 
\end{align*}
We will show that
\begin{align*}
& \Coker(\id - \bar{\gamma}_{\rho,0}) \text{ in }
 K_0(\FRE) / (\id - \gamma_{\eta,0}) K_0(\FRE)  \\
\cong &
K_0(\FRE) / 
(  (\id - \gamma_{\eta,0}) K_0(\FRE) )
 + (\id- \gamma_{\rho,0}) K_0(\FRE) ).
\end{align*}
Put
$
H_{\rho,\eta}
 =
(\id - \gamma_{\eta,0}) K_0(\FRE) )
 + (\id- \gamma_{\rho,0}) K_0(\FRE) )
$
the subgroup 
of 
$K_0(\FRE)$
generated by
$(\id - \gamma_{\rho,0})K_0(\FRE)$
and
$(\id - \gamma_{\eta,0})K_0(\FRE).$
Set the quotient maps
\begin{align*}
 K_0(\FRE) 
\overset{q_\eta}{\longrightarrow} 
& K_0(\FRE) / (\id - \gamma_{\eta,0}) K_0(\FRE) \\
 \overset{q_{(\id-\gamma_{\rho,0})}}{\longrightarrow}
& \Coker(\id - \bar{\gamma}_{\rho,0})  
\text{ in }
 K_0(\FRE) / (\id - \gamma_{\eta,0}) K_0(\FRE)  
\end{align*}
and
$\Phi = q_{(\id-\gamma_{\rho,0})} \circ q_\eta: 
K_0(\FRE) 
\longrightarrow
\Coker(\id- \bar{\gamma}_{\rho,0})  
\text{ in }
 K_0(\FRE) / (\id - \gamma_{\eta,0}) K_0(\FRE).
$
As
$(\id - \gamma_{\rho,0})$
commutes with
$(\id -\gamma_{\eta,0})$,
one has
$$
(\id  - \gamma_{\eta,0}) K_0(\FRE) \subset \Ker(\Phi),
\qquad
(\id - \gamma_{\rho,0}) K_0(\FRE) \subset \Ker(\Phi).
$$
Hence we have
$H_{\rho,\eta} \subset \Ker(\Phi).$

On the other hand,
for $g \in \Ker(\Phi)$,
we have
$g \in 
(\id- \bar{\gamma}_{\rho,0})
( K_0(\FRE) / (\id - \gamma_{\eta,0}) K_0(\FRE) )
$ 
so that
$ g= (\id- \gamma_{\rho,0})[h]$ 
for some
$[h] \in K_0(\FRE) / (\id - \gamma_{\eta,0}) K_0(\FRE)$.
Hence
$g = (\id - \gamma_{\rho,0})h + (\id- \gamma_{\rho,0}) (\id - \gamma_{\eta,0}) K_0(\FRE)$
so that 
$ g \in H_{\rho,\eta}$.
Hence we have
$\Ker(\Phi) \subset H_{\rho,\eta}$
and 
$\Ker(\Phi) = H_{\rho,\eta}$. 
As
\begin{align*}
& 
( K_0(\FRE) / (\id - \gamma_{\eta,0}) K_0(\FRE) )
/ (\id- \bar{\gamma}_{\rho,0})
(( K_0(\FRE) / (\id - \gamma_{\eta,0}) K_0(\FRE) ) \\
\cong &
K_0(\FRE) / 
(\id - \gamma_{\eta,0}) K_0(\FRE) )
 + (\id - \gamma_{\rho,0}) K_0(\FRE) ),
\end{align*}
we have
\begin{align*}
& \Coker(\id - \beta_{\rho,1}) \text{ in } K_1((\ORE)^{\hat{\rho}}) \\
\cong
& K_0(\FRE) / 
(\id - \gamma_{\eta,0}) K_0(\FRE) )
 + (\id - \gamma_{\rho,0}) K_0(\FRE) ).
\end{align*}

(ii) 
The assertions are similarly shown to (i).
\end{pf}
Therefore we have

\begin{thm}
Assume that $K_1(\FRE) = 0$.
There exist short exact sequences:
\begin{enumerate}
\renewcommand{\labelenumi}{(\roman{enumi})}
\item
\begin{align*}
0 
&
\longrightarrow 
K_0(\FRE)
/ 
(\id - \gamma_{\rho,0})K_0(\FRE) +
(\id - \gamma_{\eta,0})K_0(\FRE) \\
& \longrightarrow 
K_0(\ORE)\\
& \longrightarrow
\Ker(\id - \gamma_{\rho,0}) \cap 
\Ker(\id - \gamma_{\eta,0})  \text{ in } K_0(\FRE)) \\
& \longrightarrow 
0.
\end{align*}
\item
\begin{align*}
0 
&
\longrightarrow 
( \Ker(\id - \gamma_{\eta,0}) 
\text{ in }
K_0(\FRE) )
/ 
(\id - \gamma_{\rho,0})
( \Ker(\id - \gamma_{\eta,0}) 
\text{ in }
K_0(\FRE) ) \\
& \longrightarrow 
K_1(\ORE)\\
& \longrightarrow
\Ker(id - \gamma_{\rho,0}) 
\text{ in }
K_0(\FRE)
/ 
(\id - \bar{\gamma}_{\eta,0})K_0(\FRE) \\
& \longrightarrow 
0.
\end{align*}
\end{enumerate}
\end{thm}
As a corollary we have
\begin{cor}
Suppose $K_1(\FRE) =0$. 
We then have
\begin{enumerate}
\renewcommand{\labelenumi}{(\roman{enumi})}
\item
\begin{align*}
0 &
\longrightarrow 
\Coker(\id - \bar{\gamma}_{\rho,0})
\text{ in }
K_0(\FRE)/(\id - \gamma_{\eta,0})K_0(\FRE) \\
& \longrightarrow 
K_0(\ORE)\\
& \longrightarrow
\Ker(\id - \gamma_{\rho,0}) 
\text{ in } (\Ker(\id - \gamma_{\eta,0})  \text{ in } K_0(\FRE)) \\
& \longrightarrow 
0.
\end{align*}
\item
\begin{align*}
0 &
\longrightarrow 
(\Ker(\id - \gamma_{\eta,0})  \text{ in } K_0(\FRE))
/(\id - \gamma_{\rho,0})(\Ker(\id - \gamma_{\eta,0}) \text{ in } K_0(\FRE)) \\
& \longrightarrow 
K_1(\ORE)\\
& \longrightarrow
\Ker(\id - \bar{\gamma}_{\rho,0}) 
\text{ in } (K_0(\FRE)/(\id - \gamma_{\eta,0})K_0(\FRE) ) \\
& \longrightarrow 
0.
\end{align*}
\end{enumerate}
\end{cor}

\section{K-Theory formulae}
We henceforth denote the endomorphisms $\gamma_{\rho,0}, \gamma_{\eta,0}$ 
on $K_0(\FRE)$ by 
$\gamma_\rho, \gamma_\eta$
respectively. 

 In this section, 
 we will prove 
more useful formulae for the K-groups
$K_i(\ORE)$
under certain additional assumption on $\CTDS$.
The assumed condition on $\CTDS$ is the following:

\noindent
{\bf Definition.}
A $C^*$-textile dynamical system
$\CTDS$ is said to {\it form square}
if the 
$C^*$-subalgebra
$C^*(\rho_\alpha(1) : \alpha \in \Sigma^\rho)$
of $\A$
generated by
the projections
$\rho_\alpha(1), \alpha \in \Sigma^\rho$ 
coincides with 
the 
$C^*$-subalgebra
$C^*(\eta_a(1) : a \in \Sigma^\eta)$
of $\A$ generated by
the projections
$\eta_a(1), a \in \Sigma^\eta$.

\begin{lem}
Assume that $\CTDS$ forms square.
Put
for $l \in \Zp$
\begin{equation*}
\A_l^\rho 
 = C^*(\rho_\mu(1) : \mu \in B_l(\Lambda_\rho)),\qquad
\A_l^\eta 
 = C^*(\eta_\xi(1) : \xi \in B_l(\Lambda_\eta)).
\end{equation*}
Then $\A_l^\rho = \A_l^\eta$.
\end{lem}
\begin{pf}
By assumption, 
we have
$\A_1^\rho = \A_1^\eta$.
Hence the desired equality for $l=1$ holds.
Suppose that the equalities hold  for all $l \le k$
for some $k \in {\Bbb N}$. 
For 
$\mu = \mu_1\mu_2 \cdots \mu_k \mu_{k+1} \in B_{k+1}(\Lambda_\rho)$
we have
$\rho_\mu(1) 
= \rho_{\mu_{k+1}}(\rho_{\mu_1\mu_2 \cdots \mu_k}(1))$
so that 
$\rho_\mu(1) \in \rho_{\mu_{k+1}}(\A_k^\rho)$.
By the $\kappa$-commutation relation,
one sees that 
$$
\rho_{\mu_{k+1}}(\A_k^\rho)
\subset C^*(\eta_\xi(\rho_\alpha(1)) : 
\xi \in B_k(\Lambda_\eta),\alpha \in \Sigma^\rho).
$$
Since
$
C^*(\rho_\alpha(1) : \alpha \in \Sigma^\rho)
=
 C^*(\eta_a(1) : a \in \Sigma^\eta),
$
one knows that
the algebra 
$
 C^*(\eta_\xi(\rho_\alpha(1)) : 
\xi \in B_k(\Lambda_\eta),\alpha \in \Sigma^\rho)
$
is contained in 
$\A_{k+1}^\eta$
so that
$\rho_{\mu_{k+1}}(\A_k^\eta) \subset \A_{k+1}^\eta$.
Therefore 
we have
$\rho_\mu(1) \in \A_{k+1}^\eta$ 
so that
$\A_{k+1}^\rho \subset \A_{k+1}^\eta$
and hence
$\A_{k+1}^\rho = \A_{k+1}^\eta$.
\end{pf}
Therefore we have

\begin{lem}
Assume that $\CTDS$ forms square.
Put
for $j,k \in \Zp$
\begin{align*}
\A_{j,k} 
& = C^*(\rho_\mu(\eta_\zeta(1)) : 
\mu \in B_j(\Lambda_\rho), \zeta \in B_k(\Lambda_\eta))\\
&( = C^*(\eta_\xi(\rho_\nu(1)) : 
\xi \in B_k(\Lambda_\eta), \nu \in B_j(\Lambda_\rho))).
\end{align*}
Then $\A_{j,k}$ is commutative and of finite dimensional such that 
\begin{equation*}
\A_{j,k} = \A_{j+k}^\rho ( = \A_{j+k}^\eta). 
\end{equation*}
Hence
$\A_{j,k} = \A_{j',k'}$ if $j+k = j' +k'$.
\end{lem}
\begin{pf}
Since
$\eta_\zeta(1) \in Z_\A$ and
$\rho_\mu(Z_\A) \subset Z_\A$,
the algebra 
$\A_{j,k}$ belongs to the center $Z_\A$ of $\A$.
By the preceding lemma,
we have
\begin{equation*}
\A_{j,k} 
 = C^*(\rho_\mu(\rho_\nu(1)) : 
\mu \in B_j(\Lambda_\rho), \nu \in B_k(\Lambda_\rho))
= \A_{j+k}^\rho.
\end{equation*}
\end{pf}
For $j,k\in \Zp$,
put $l =j+k$.
We denote by 
$\A_l$ 
the commutative finite dimensional algebra $\A_{j,k}$.
Put 
$m(l) = \dim \A_l$.
Take the finite sequence of minimal projections
$E_i^l, i=1,2,\dots,m(l)$ in $\A_l$ 
such that
$\sum_{i=1}^{m(l)} E_i^l = 1$.
 Hence we have
 $\A_l = \sum_{i=1}^{m(l)} {\Bbb C}E_i^l$.
Since 
$\rho_\alpha(\A_l) \subset \A_{l+1}$,
there exists
$A_{l,l+1}^\rho(i,\alpha,n)$, 
which takes $0$ or $1$,
such that
\begin{equation*}
\rho_\alpha(E_i^l) = 
\sum_{n=1}^{m(l+1)}
A_{l,l+1}^\rho(i,\alpha,n) E_n^{l+1},
\qquad
\alpha \in \Sigma^\rho, \, i=1,\dots,m(l).
\end{equation*} 
Similarly,
there exists
$A_{l,l+1}^\eta(i,a,n)$, 
which takes $0$ or $1$,
such that
\begin{equation*}
\eta_a(E_i^l) =
\sum_{n=1}^{m(l+1)}
A_{l,l+1}^\eta(i,a,n) E_n^{l+1},
\qquad
a \in \Sigma^\eta, \, i=1,\dots,m(l).
\end{equation*}

Let 
$
N_{j,k}(i)
$
be the cardinal number of the set
$$
\{ (\mu,\zeta) \in B_j(\Lambda_\rho)\times B_k(\Lambda_\eta) \mid
 \rho_\mu(\eta_\zeta(1)) \ge E_i^l \}. 
$$
Set for $i=1,\dots,m(l)$
\begin{align*}
{\cal F}_{j,k}(i)
& = C^*( S_\mu T_\zeta E_i^l x E_i^l  T_\xi^* S_\nu^* \mid 
         \mu,\nu \in B_j(\Lambda_\rho),\zeta,\xi \in B_k(\Lambda_\eta), 
         x \in \A) \\
& = C^*( T_\zeta S_\mu E_i^l x E_i^l S_\nu^* T_\xi^*  \mid 
         \mu,\nu \in B_j(\Lambda_\rho),\zeta,\xi \in B_k(\Lambda_\eta), 
         x \in \A). 
\end{align*}
Since $E_i^l$ is a central projection in $\A$,
we have
\begin{lem}
\begin{enumerate}
\renewcommand{\labelenumi}{(\roman{enumi})}
\item
 ${\cal F}_{j,k}(i)$ is isomorphic to the matrix algebra 
 $M_{N_{j,k}(i)}(E_i^l  \A E_i^l ) 
 (=M_{N_{j,k}(i)}({\Bbb C})\otimes E_i^l \A E_i^l)$
over $E_i^l  \A E_i^l$.
\item
$
{\cal F}_{j,k} = {\cal F}_{j,k}(1)\oplus \cdots \oplus {\cal F}_{j,k}(m(l)).
$
\end{enumerate} 
\end{lem}
\begin{pf}
(i)
For 
$(\mu,\zeta) \in B_j(\Lambda_\rho)\times B_k(\Lambda_\eta)$
with
$S_\mu T_\zeta E_i^l \ne 0$,
one has
$ \eta_\zeta(\rho_\mu(1)) E_i^l \ne 0$
so that
$ \eta_\zeta(\rho_\mu(1))
 \ge E_i^l.
$
Hence
$
(S_\mu T_\zeta E_i^l)^* S_\mu T_\zeta E_i^l 
=  E_i$.
One sees that 
the set
$$
\{ S_\mu T_\zeta E_i^l 
\mid (\mu,\zeta) \in B_j(\Lambda_\rho)\times B_k(\Lambda_\eta);
S_\mu T_\zeta E_i^l \ne 0
\}
$$
consist of isometries which give rise to  matrix units of 
${\cal F}_{j,k}(i)$
such that
${\cal F}_{j,k}(i)$
is isomorphic to
$M_{N_{j,k}(i)}(E_i^l  \A E_i^l )$.

(ii) Since
$\A = E_1^l \A E_1^l \oplus \cdots \oplus E_{m(l)}^l \A E_{m(l)}^l$
the assertion is easy.
\end{pf} 
Define
$\lambda_{\rho*}, \lambda_{\eta*}: K_0(\A) \longrightarrow K_0(\A)$
by setting
\begin{equation*}
\lambda_{\rho*}([p]) 
= \sum_{\alpha \in \Sigma^\rho} [\rho_\alpha\otimes 1_n (p)],
\qquad
\lambda_{\eta*}([p]) 
= \sum_{a \in \Sigma^\eta} [\eta_a\otimes 1_n (p)]
\end{equation*}
for a projection
$p \in M_n(\A)$
for some
$ n \in {\Bbb N}.$ 
\begin{lem}
Assume that
$\CTDS$ forms square.
There exists an isomorphism
$$
\Phi_{j,k}: K_0({\cal F}_{j,k}) \longrightarrow K_0(\A)
$$
such that the  following diagrams are commutative:
\begin{enumerate}
\renewcommand{\labelenumi}{(\roman{enumi})}
\item
\begin{equation*}
\begin{CD}
K_0({\cal F}_{j,k}) @>\iota_{+1,*}>> K_0({\cal F}_{j+1,k})   \\
@V{\Phi_{j,k}}VV   @V{\Phi_{j+1,k}}VV \\
K_0(\A) @>\lambda_{\rho*}>> K_0(\A) 
\end{CD}
\end{equation*}
\item
\begin{equation*}
\begin{CD}
K_0({\cal F}_{j,k}) @>\iota_{*,+1}>> K_0({\cal F}_{j,k+1})   \\
@V{\Phi_{j,k}}VV   @V{\Phi_{j,k+1}}VV \\
K_0(\A) @>\lambda_{\eta*}>> K_0(\A) 
\end{CD}
\end{equation*}
\end{enumerate}
\end{lem}
\begin{pf}
Put
for $i=1,2,\cdots m(l)$
$$
P_i = \sum_{\mu \in B_j(\Lambda_\rho), \zeta \in B_k(\Lambda_\eta)}
S_\mu T_\zeta E_i^l T_\zeta^* S_\mu^*
$$
Then $P_i$ is a projection which belongs to the center of ${\cal F}_{j,k}$
such that
$\sum_{i=1}^{m(l)} P_i = 1$.
For 
$X \in {\cal F}_{j,k}$,
one has 
$P_i X P_i \in {\cal F}_{j,k}(i)$
such that
$$
X = \sum_{i=1}^{m(l)} P_i X P_i 
\in 
\bigoplus_{i=1}^{m(l)}  {\cal F}_{j,k}(i).  
$$
Define an isomorphism
$$
\varphi_{j,k} : 
X \in {\cal F}_{j,k} 
\longrightarrow 
\sum_{i=1}^{m(l)} P_i X P_i
\in
\bigoplus_{i=1}^{m(l)} {\cal F}_{j,k}(i)
$$
which induces an isomorphism on their K-groups
$$
\varphi_{j,k *} :
K_0({\cal F}_{j,k}) 
\longrightarrow 
\bigoplus_{i=1}^{m(l)} K_0({\cal F}_{j,k}(i)).
$$
Take and fix 
$\nu(i), \mu(i) \in B_j(\Lambda_\rho)$
and
$\zeta(i), \xi(i) \in B_k(\Lambda_\eta)$
such that
$$
T_{\xi(i)} S_{\nu(i)} = S_{\mu(i)} T_{\zeta(i)}
\quad
\text{ and }
\quad
 T_{\xi(i)} S_{\nu(i)} E_i^l \ne 0.
$$
Hence
$S_{\nu(i)}^* T_{\xi(i)}^* T_{\xi(i)} S_{\nu(i)} \ge E_i^l$.
Since ${\cal F}_{j,k}(i)$
is isomorphic to 
$ M_{N_{j,k(i)}}({\Bbb C}) \otimes E_i^l \A E_i^l$,
the embedding
\begin{equation*}
\iota_{j,k}(i): x \in E_i^l \A E_i^l 
\longrightarrow  
T_{\xi(i)} S_{\nu(i)} x S_{\nu(i)}^* T_{\xi(i)}^*
\in
{\cal F}_{j,k}(i)
\end{equation*}
induces an isomorphism on their K-groups
\begin{equation*}
\iota_{j,k}(i)_*: K_0( E_i^l \A E_i^l) 
\longrightarrow  
K_0({\cal F}_{j,k}(i)).
\end{equation*}
Put
\begin{equation*}
\psi_{j,k} = \oplus_{i=1}^{m(l)}\iota_{j,k}(i): 
\oplus_{i=1}^{m(l)}E_i^l \A E_i^l \longrightarrow 
\oplus_{i=1}^{m(l)} {\cal F}_{j,k}(i)
\end{equation*}
and hence
\begin{equation*}
\psi_{j,k*} = \bigoplus_{i=1}^{m(l)} \iota_{i*}: 
\bigoplus_{i=1}^{m(l)} K_0( E_i^l \A E_i^l) 
\longrightarrow  
\bigoplus_{i=1}^{m(l)} K_0({\cal F}_{j,k}(i)).
\end{equation*}
Hence we have isomorphisms
\begin{equation*}
K_0({\cal F}_{j,k})
\overset{\varphi_{j,k*}}{\longrightarrow}
\bigoplus_{i=1}^{m(l)}
K_0({\cal F}_{j,k}(i))
\overset{{\psi_{j,k*}}^{-1}}{\longrightarrow}
\bigoplus_{i=1}^{m(l)} K_0( E_i^l \A E_i^l). 
\end{equation*}
Since
$K_0(\A) = \bigoplus_{i=1}^{m(l)} K_0( E_i^l \A E_i^l)$,
we have an isomorphism
\begin{equation*}
\Phi_{j,k} ={\psi_{j,k*}}^{-1} \circ \varphi_{j,k*} : 
K_0({\cal F}_{j,k}) \longrightarrow 
K_0(\A).
\end{equation*}

(i)
It suffices to show the following diagram
\begin{equation*}
\begin{CD}
K_0({\cal F}_{j,k}) @>\iota_{+1,*}>> K_0({\cal F}_{j+1,k})   \\
@V{\varphi_{j,k*}}VV   @V{\varphi_{j+1,k*}}VV \\
\bigoplus_{i=1}^{m(l)} K_0({\cal F}_{j,k}(i)) @. 
\bigoplus_{i=1}^{m(l)} K_0({\cal F}_{j+1,k}(i)) \\
@A{\psi_{j,k*}}AA   @A{\psi_{j+1,k*}}AA \\
K_0(\A) @>\lambda_{\rho*}>> K_0(\A) 
\end{CD}
\end{equation*}
 is commutative.
For
$a  = \sum_{i=1}^{m(l)} E_i^l a E_i^l \in \A$,
we have
$$
\psi_{j,k}(a)
=\sum_{i=1}^{m(l)} 
T_{\xi(i)} S_{\nu(i)} E_i^l a E_i^l S_{\nu(i)}^*T_{\xi(i)}^*
= \sum_{i=1}^{m(l)} 
S_{\mu(i)} T_{\zeta(i)}E_i^l a E_i^l T_{\zeta(i)}^* S_{\mu(i)}^*. 
$$
Since
$
P_i T_{\xi(i)} S_{\nu(i)} E_i^l a E_i^l S_{\nu(i)}^*T_{\xi(i)}^* P_i 
= T_{\xi(i)} S_{\nu(i)} E_i^l a E_i^l S_{\nu(i)}^*T_{\xi(i)}^*,
$
we have 
$$
\varphi_{j,k}^{-1}\circ \psi_{j,k} (a)
= \sum_{i=1}^{m(l)} 
T_{\xi(i)} S_{\nu(i)} E_i^l a E_i^l S_{\nu(i)}^*T_{\xi(i)}^*
$$
so that
\begin{equation*} 
\iota_{+1,*}\circ \varphi_{j,k}^{-1}\circ \psi_{j,k} (a)
 = \sum_{\alpha \in \Sigma} \sum_{i=1}^{m(l)} 
T_{\xi(i)} S_{\nu(i)\alpha} \rho_\alpha(E_i^l a E_i^l) 
S_{\nu(i)\alpha}^*T_{\xi(i)}^*.
\end{equation*}
Since 
\begin{equation*}
S_{\nu(i) \alpha} \rho_\alpha(E_i^l a E_i^l) S_{\nu(i) \alpha}^*
= 
\sum_{n=1}^{m(l+1)} A_{l,l+1}^\rho(i,\alpha,n)
S_{\nu(i) \alpha} E_n^{l+1} \rho_\alpha(a ) E_n^{l+1} S_{\nu(i) \alpha}^*
\end{equation*}
and
$
A_{l,l+1}^\rho(i,\alpha,n)S_{\nu(i) \alpha} E_n^{l+1}
=S_{\nu(i) \alpha} E_n^{l+1},
$
we have
\begin{equation*}
\sum_{\alpha \in \Sigma^\rho}
S_{\nu(i) \alpha} \rho_\alpha(E_i^l a E_i^l) S_{\nu(i) \alpha}^*
 = 
\sum_{n=1}^{m(l+1)} 
\sum_{\alpha \in \Sigma^\rho} 
S_{\nu(i) \alpha} E_n^{l+1} \rho_\alpha(a ) E_n^{l+1} S_{\nu(i) \alpha}^* 
\end{equation*}
so that
\begin{equation*}
\iota_{+1,*}\circ \varphi_{j,k}^{-1}\circ \psi_{j,k} (a)
 = \sum_{\alpha \in \Sigma} \sum_{i=1}^{m(l)} \sum_{n=1}^{m(l+1)} 
    T_{\xi(i)} S_{\nu(i) \alpha} E_n^{l+1} \rho_\alpha(a ) E_n^{l+1} 
    S_{\nu(i) \alpha}^*T_{\xi(i)}^*.
\end{equation*}
On the other hand,
\begin{align*}
\psi_{j,k} (\lambda_\rho(a) )
& =\psi_{j,k}(\sum_{\alpha \in \Sigma^\rho} \rho_\alpha(a)) \\
& =\psi_{j,k}(\sum_{n=1}^{m(l+1)} 
   \sum_{\alpha \in \Sigma^\rho} E_n^{l+1}\rho_\alpha(a))E_n^{l+1} \\
& = \sum_{\alpha \in \Sigma} \sum_{i=1}^{m(l)} \sum_{n=1}^{m(l+1)} 
    T_{\xi(i)} S_{\nu(i) \alpha} E_n^{l+1} \rho_\alpha(a ) E_n^{l+1} 
    S_{\nu(i) \alpha}^*T_{\xi(i)}^*.
\end{align*}
Therefore we have 
\begin{equation*}
\iota_*\circ \varphi_{j,k}^{-1}\circ \psi_{j,k} (a)
 =
\psi_{j,k} (\lambda_\rho(a) ).
\end{equation*}
(ii) is symmetric to (i).
\end{pf}


Define the abelian groups of inductive limits:
\begin{equation*}
G_\rho = \lim \{\lambda_\rho: K_0(\A) \longrightarrow K_0(\A) \},\qquad
G_\eta = \lim \{\lambda_\eta: K_0(\A) \longrightarrow K_0(\A) \}.
\end{equation*}
Put for $j,k \in \Zp$ the subalgebras of $\FRE$
\begin{align*}
{\cal F}_{\rho,k}
& = C^*( T_\zeta S_\mu  x S_\nu^* T_\xi^* \mid 
         \mu,\nu \in B_*(\Lambda_\rho),|\mu| = |\nu|, 
         \zeta,\xi \in B_k(\Lambda_\eta), 
         x \in \A) \\
& = C^*( T_\zeta y T_\xi^* \mid 
         \zeta,\xi \in B_k(\Lambda_\eta), 
         y \in {\cal F}_\rho ) \\
\intertext{ and }
{\cal F}_{j,\eta}
& = C^*( S_\mu T_\zeta x T_\xi^* S_\nu^*  \mid 
         \mu,\nu \in B_j(\Lambda_\rho),
         \zeta,\xi \in B_*(\Lambda_\eta),|\zeta| = |\xi|, 
         x \in \A) \\
& = C^*( S_\mu y S_\nu^*  \mid 
         \mu,\nu \in B_j(\Lambda_\rho),
         y \in {\cal F}_\eta ). 
\end{align*}

\begin{lem}
For $j,k \in \Zp$,
there exist isomorphisms
\begin{equation*}
\Phi_{\rho,k}: K_0({\cal F}_{\rho,k}) \longrightarrow G_\rho, \qquad
\Phi_{j,\eta}: K_0({\cal F}_{j,\eta}) \longrightarrow G_\eta
\end{equation*} 
such that the following diagrams are commutative:
\begin{enumerate}
\renewcommand{\labelenumi}{(\roman{enumi})}
\item
\begin{equation*}
\begin{CD}
K_0({\cal F}_{j,k}) @>\iota_{+1,*}>> K_0({\cal F}_{j+1,k}) 
@>\iota_{+1,*}>> \cdots @>\iota_{+1,*}>>
K_0({\cal F}_{\rho,k}) \\
@V{\Phi_{j,k}}VV   @V{\Phi_{j+1,k}}VV  @. @V{\Phi_{\rho,k}}VV   \\
K_0(\A) @>\lambda_{\rho*}>> K_0(\A) @>\lambda_{\rho*}>> \cdots @>\lambda_{\rho*}>>
G_\rho  
\end{CD}
\end{equation*}
\item
\begin{equation*}
\begin{CD}
K_0({\cal F}_{j,k}) @>\iota_{*,+1}>> K_0({\cal F}_{j,k+1}) 
@>\iota_{*,+1}>> \cdots @>\iota_{*,+1}>>
K_0({\cal F}_{j,\eta}) \\
@V{\Phi_{j,k}}VV   @V{\Phi_{j,k+1}}VV  @. @V{\Phi_{j,\eta}}VV   \\
K_0(\A) @>\lambda_{\eta*}>> K_0(\A) @>\lambda_{\eta*}>> \cdots @>\lambda_{\eta*}>>
G_\eta  
\end{CD}
\end{equation*}
\end{enumerate}
\end{lem} 
%
\begin{lem}
If 
$\xi = \xi_1 \cdots \xi_k \in B_k(\Lambda_\eta), \nu = \nu_1\cdots \nu_j \in B_j(\Lambda_\rho)$
and
$i=1,\dots,m(l)$
satisfy the condition 
$\rho_\nu(\eta_\xi(1)) \ge E_i^l$
where $ l= j+k$,
then 
$T_{\xi_1}^* T_\xi S_\nu E_i^l =    T_{\bar{\xi}} S_\nu E_i^l$
where
$\bar{\xi} = \xi_2 \cdots \xi_k$. 
\end{lem}
\begin{pf}
Since
$T_{\xi_1}^* T_\xi 
=T_{\xi_1}^*T_{\xi_1} T_{\bar{\xi}}T_{\bar{\xi}}^*T_{\bar{\xi}}
=T_{\bar{\xi}}T_{\bar{\xi}}^* T_{\xi_1}^*T_{\xi_1} T_{\bar{\xi}}
=T_{\bar{\xi}}T_{\xi}^* T_{\xi},
$
we have
\begin{equation*}
T_{\xi_1}^* T_\xi S_\nu E_i^l
 = T_{\bar{\xi}} S_\nu S_\nu^* T_{\xi}^* T_{\xi} S_\nu E_i^l 
 = T_{\bar{\xi}} S_\nu \rho_\nu(\eta_{\xi}(1))  E_i^l 
 = T_{\bar{\xi}} S_\nu E_i^l.
\end{equation*}
\end{pf}

\begin{lem}
For $k,j$ we have 
\begin{enumerate}
\renewcommand{\labelenumi}{(\roman{enumi})}
\item The restriction of $\gamma_\eta^{-1}$
to $K_0({\cal F}_{j,k})$ 
makes the following diagram commutative:
\begin{equation*}
\begin{CD}
K_0({\cal F}_{j,k}) @>\gamma_\eta^{-1}>> K_0({\cal F}_{j,k-1}) 
@>\iota_{*,+1}>> K_0({\cal F}_{j,k}) \\
@V{\Phi_{j,k}}VV   @.  @V{\Phi_{j,k}}VV \\
K_0(\A)   @. \overset{\lambda_{\eta*}}{\longrightarrow} @. K_0(\A).
\end{CD}
\end{equation*}
\item The restriction of $\gamma_\rho^{-1}$
to $K_0({\cal F}_{j,k})$ 
makes the following diagram commutative:
\begin{equation*}
\begin{CD}
K_0({\cal F}_{j,k}) @>\gamma_\rho^{-1}>> K_0({\cal F}_{j-1,k}) 
@>\iota_{+1,*}>> K_0({\cal F}_{j,k}) \\
@V{\Phi_{j,k}}VV   @.  @V{\Phi_{j,k}}VV \\
K_0(\A)  @. \overset{\lambda_{\rho*}}{\longrightarrow} @. K_0(\A).
\end{CD}
\end{equation*}
\end{enumerate}
\end{lem}
\begin{pf}
(i)
Put 
$l = j+k$.
Take a projection 
$p \in M_n(\A)$ 
for some $n \in {\Bbb N}$.
Since
$
\A\otimes M_n({\Bbb C}) = \sum_{i=1}^{m(l)} (E_i^l \otimes 1) (\A \otimes M_n) (E_i^l \otimes 1),
$
by putting   
$
p_i^l 
= (E_i^l \otimes 1) p (E_i^l \otimes 1)\in 
(E_i^l \otimes 1) (\A \otimes M_n) (E_i^l \otimes 1) = M_n(E_i^l \A E_i^l),
$
we have
$p = \sum_{i=1}^{m(l)} p_i^l.$
Take 
$ \xi(i)=\xi_1(i)\cdots \xi_k(i) \in B_k(\Lambda_\eta), 
  \nu(i)= \nu_1(i)\cdots \nu_j(i) \in B_j(\Lambda_\rho)$, 
such that
$\rho_{\nu(i)}(\eta_{\xi(i)}(1)) \ge E_i^l$
and
put
$\bar{\xi}(i) = \xi_2(i)\cdots \xi_k(i)$
so that
$\xi(i)= \xi_1(i)\bar{\xi}(i).$
Since
$$
\psi_{j,k*}([p]) 
= \sum_{i=1}^{m(l)}\oplus 
[ (T_{\xi(i)} S_{\nu(i)} \otimes 1_n) 
p_i^l (S_{\nu(i)}^*T_{\xi(i)}^*  \otimes 1_n) ]
\in
\bigoplus_{i=1}^{m(l)} K_0({\cal F}_{j,k}(i)).
$$
As
$$
(T_{\xi(i)} S_{\nu(i)} \otimes 1_n) p_i^l (S_{\nu(i)}^*T_{\xi(i)}^*  
\otimes 1_n)
\le
T_{\xi_1(i)} T_{\xi_1(i)}^* \otimes 1_n,
$$
by the preceding lemma we have
$$
T_{\xi_1(i)}^* T_{\xi(i)} S_{\nu(i)} E_i^l 
=T_{\bar{\xi}(i)} S_{\nu(i)} E_i^l
$$
so that 
$$
\gamma_\eta^{-1}
([ (T_{\xi(i)} S_{\nu(i)} \otimes 1_n) 
p_i^l (S_{\nu(i)}^*T_{\xi(i)}^*  \otimes 1_n) ]
=
[(T_{\bar{\xi}(i)} S_{\nu(i)}\otimes 1_n) 
p_i^l (S_{\nu(i)}^*T_{\bar{\xi}(i)}^*\otimes 1_n)].
$$
Hence
$K_0({\cal F}_{j,k})$ 
goes to
$K_0({\cal F}_{j,k-1})$
by the homomorphism
$
\gamma_\eta^{-1}.
$
Take
$\mu(i) \in B_j(\Lambda_\rho),
\bar{\zeta}(i) \in B_{k-1}(\Lambda_\eta)
$
such that
$T_{\bar{\xi}(i)} S_{\nu(i)} = S_{\mu(i)} T_{\bar{\zeta}(i)}$
for $i=1,\dots,m(l)$.
The element
\begin{align*}
& \sum_{i=1}^{m(l)}
[(T_{\bar{\xi}(i)} S_{\nu(i)}\otimes 1_n) 
p_i^l (S_{\nu(i)}^*T_{\bar{\xi}(i)}^*\otimes 1_n)] \\
=&
\sum_{i=1}^{m(l)}
[(S_{\mu(i)} T_{\bar{\zeta}(i)}\otimes 1_n) 
p_i^l (T_{\bar{\zeta}(i)}^*S_{\mu(i)}^* \otimes 1_n)]
\in K_0({\cal F}_{j,k-1})
\end{align*}
goes to
$$
\sum_{i=1}^{m(l)} \sum_{a \in \Sigma^{\eta}}
[(S_{\mu(i)} T_{\bar{\zeta}(i) a}\otimes 1_n) (T_a^* \otimes 1_n) 
p_i^l (T_a \otimes 1_n) (T_{\bar{\zeta}(i) a}^*S_{\mu(i)}^* \otimes 1_n)]
\in K_0({\cal F}_{j,k})
$$
by $\iota_{*,+1}$.
The element is expressed as
\begin{equation}
\sum_{h=1}^{m(l)} \oplus 
\sum_{i=1}^{m(l)} \sum_{a \in \Sigma^{\eta}}
[(S_{\mu(i)} T_{\bar{\zeta}(i) a}\otimes 1_n) E_h^l (T_a^* \otimes 1_n) 
p_i^l 
(T_a \otimes 1_n) E_h^l (T_{\bar{\zeta}(i) a}^*S_{\mu(i)}^* \otimes 1_n)]
\end{equation}
in  
$\bigoplus_{h=1}^{m(l)}  K_0({\cal F}_{j,k}(h))$.

On the other hand,
$$
\lambda_{\eta*}([p]) = \sum_{a \in \Sigma^\eta}[(T_a^* \otimes 1_n) p (T_a \otimes 1_n)]  
\in K_0(\A).
$$
The element
$$
\sum_{a \in \Sigma^\eta}[(T_a^* \otimes 1_n) p (T_a \otimes 1_n)]  
=
\sum_{h=1}^{m(l)} \oplus
\sum_{a \in \Sigma^\eta}[E_h^l (T_a^* \otimes 1_n) p (T_a \otimes 1_n)E_h^l]  
\in
\bigoplus_{h=1}^{m(l)} K_0(E_h^l \A E_h^l )
$$
is expressed as
\begin{align*}
& \sum_{h=1}^{m(l)} \oplus
\sum_{a \in \Sigma^\eta}
[(T_{\xi(h)} S_{\nu(h)}E_h^l\otimes 1_n) (T_a^* \otimes 1_n) p (T_a \otimes 1_n)
(E_h^l S_{\nu(h)}^* T_{\xi(h)}^* \otimes 1_n)]  \\
= 
& \sum_{h=1}^{m(l)} \oplus
\sum_{a \in \Sigma^\eta} \sum_{i=1}^{m(l)}
[(T_{\xi(h)} S_{\nu(h)}E_h^l\otimes 1_n) (T_a^* \otimes 1_n) p_i^l (T_a \otimes 1_n)
(E_h^l S_{\nu(h)}^* T_{\xi(h)}^* \otimes 1_n )] 
\end{align*}
in $\bigoplus_{h=1}^{m(l)} K_0({\cal F}_{j,k}(h))$.
Take
$\mu'(h) \in B_j(\Lambda_\rho),
\zeta'(h) \in B_k(\Lambda_\eta)$
such that
$T_{\xi(h)} S_{\nu(h)} =  S_{\mu'(h) } T_{\zeta'(h)}$
so that the above element is
\begin{equation}
 \sum_{h=1}^{m(l)} \oplus
\sum_{i=1}^{m(l)}
\sum_{a \in \Sigma^\eta} 
[(S_{\mu'(h)} T_{\zeta'(h)}E_h^l\otimes 1_n) (T_a^* \otimes 1_n) p_i^l 
(T_a \otimes 1_n)
(E_h^l T_{\zeta'(h)}^* S_{\nu'(h)}^* \otimes 1_n )] 
\end{equation}
in $\bigoplus_{h=1}^{m(l)} K_0({\cal F}_{j,k}(h))$.
 Since for
 $h, i=1,\dots,m(l), a \in \Sigma^\eta$,
 the classes of the K-groups coincide such as 
\begin{align*}
& [(S_{\mu(i)} T_{\bar{\zeta}(i) a}\otimes 1_n) E_h^l (T_a^* \otimes 1_n) 
p_i^l 
(T_a \otimes 1_n) E_h^l (T_{\bar{\zeta}(i) a}^*S_{\mu(i)}^* \otimes 1_n)] \\
=
& [(S_{\mu'(h)} T_{\zeta'(h)}E_h^l\otimes 1_n) (T_a^* \otimes 1_n) p_i^l 
(T_a \otimes 1_n)
(E_h^l T_{\zeta'(h)}^* S_{\nu'(h)}^* \otimes 1_n )] 
\in K_0({\cal F}_{j,k}(h))
\end{align*}
we have 
the element of 
(7.1) is equal to the element of (7.2)
in $K_0({\cal F}_{j,k})$.
Therefore (i) holds.

(ii) is similar to (i).
\end{pf}
The following lemma is direct.
\begin{lem}
For $k,j$ the following diagrams are commutative.
\begin{enumerate}
\renewcommand{\labelenumi}{(\roman{enumi})}
\item
\begin{equation*}
\begin{CD}
K_0({\cal F}_{j,k}) @>\gamma_\eta^{-1}>> K_0({\cal F}_{j,k-1}) \\
@V{\iota_{+1,*}}VV     @V{\iota_{+1,*}}VV \\
K_0({\cal F}_{j+1,k}) @>\gamma_\eta^{-1}>> K_0({\cal F}_{j+1,k-1}). \\
\end{CD}
\end{equation*}
Hence 
$\gamma_\eta^{-1}$ yields  a homomorphism 
from
$
K_0({\cal F}_{\rho, k})
 = \lim_j \{ \iota_{+1,*}: K_0({\cal F}_{j,k})
\longrightarrow 
K_0({\cal F}_{j+1, k})\}
$
to
$
K_0({\cal F}_{\rho, k-1})
 = \lim_j \{ \iota_{+1,*}: K_0({\cal F}_{j,k-1})
\longrightarrow 
K_0({\cal F}_{j+1, k-1})\}.
$
\item
\begin{equation*}
\begin{CD}
K_0({\cal F}_{j,k}) @>\gamma_\rho^{-1}>> K_0({\cal F}_{j-1,k}) \\
@V{\iota_{*,+1}}VV     @V{\iota_{*,+1}}VV \\
K_0({\cal F}_{j,k+1}) @>\gamma_\rho^{-1}>> K_0({\cal F}_{j-1,k+1}) \\
\end{CD}
\end{equation*}
Hence 
$\gamma_\rho^{-1}$ yields  a homomorphism 
from
$
K_0({\cal F}_{j,\eta})
 = \lim_k \{ \iota_{*,+1}: K_0({\cal F}_{j,k})
\longrightarrow 
K_0({\cal F}_{j,k+1})\}
$
to
$
K_0({\cal F}_{j-1,\eta })
 = \lim_k \{ \iota_{*,+1}: K_0({\cal F}_{j-1,k})
\longrightarrow 
K_0({\cal F}_{j-1, k+1})\}.
$
\end{enumerate}
\end{lem}
\begin{lem}
For $k,j$ the following diagrams are commutative.
\begin{enumerate}
\renewcommand{\labelenumi}{(\roman{enumi})}
\item
\begin{equation*}
\begin{CD}
K_0({\cal F}_{\rho,k}) @>\gamma_\eta^{-1}>> K_0({\cal F}_{\rho,k-1}) \\
@V{\iota_{*,+1}}VV     @V{\iota_{*,+1}}VV \\
K_0({\cal F}_{\rho,{k+1}}) @>\gamma_\eta^{-1}>> K_0({\cal F}_{\rho,k}) \\
\end{CD}
\end{equation*}
\item
\begin{equation*}
\begin{CD}
K_0({\cal F}_{j,\eta}) @>\gamma_\rho^{-1}>> K_0({\cal F}_{j-1,\eta}) \\
@V{\iota_{+1,*}}VV     @V{\iota_{+1,*}}VV \\
K_0({\cal F}_{j+1,\eta}) @>\gamma_\rho^{-1}>> K_0({\cal F}_{j,\eta}) \\
\end{CD}
\end{equation*}
\end{enumerate}
\end{lem}
\begin{pf}
(i)
As in the proof of Lemma 7.8,
one may take an element of
$K_0({\cal F}_{\rho,k})$
as in the following form:
\begin{equation*}
 \sum_{i=1}^{m(l)}\oplus 
[ (T_{\xi(i)} S_{\nu(i)} \otimes 1_n) 
p_i^l (S_{\nu(i)}^*T_{\xi(i)}^*  \otimes 1_n) ]
\in
\bigoplus_{i=1}^{m(l)} K_0({\cal F}_{j,k}(i))
\end{equation*}
for some projection
$p \in M_n(\A)$
and
$j,l$ with $l= j+k$,
where
$
p_i^l 
= (E_i^l \otimes 1) p (E_i^l \otimes 1)\in 
(E_i^l \otimes 1) (\A \otimes M_n) (E_i^l \otimes 1) = M_n(E_i^l \A E_i^l).
$
Let
$\xi(i) = \xi_1(i) \bar{\xi}(i)$
with 
$\xi_1(i) \in \Sigma^\eta$, 
$\bar{\xi}(i) \in B_{k-1}(\Lambda_\eta)$.
One may assume that
$
T_{\xi(i)} S_{\nu(i)} \ne 0
$
so that
$
T_{\bar{\xi}(i)} S_{\nu(i)}
= S_{\nu(i)'}T_{\bar{\xi}(i)'}
$
for some
$
\nu(i)'
 \in B_j(\Lambda_\rho), \bar{\xi}(i)' \in B_{k-1}(\Lambda_\eta).
$
As in the proof of Lemma 7.8,
one has
\begin{align*}
& \gamma_\eta^{-1}
([ (T_{\xi(i)} S_{\nu(i)} \otimes 1_n) 
p_i^l (S_{\nu(i)}^*T_{\xi(i)}^*  \otimes 1_n) ] \\
= &
[(T_{\bar{\xi}(i)} S_{\nu(i)}\otimes 1_n) 
p_i^l (S_{\nu(i)}^*T_{\bar{\xi}(i)}^*\otimes 1_n)]\\
= &
[(S_{\nu(i)'}T_{\bar{\xi}(i)'} \otimes 1_n) 
p_i^l (S_{\nu(i)'}^*T_{\bar{\xi}(i)'}^*\otimes 1_n)]
\end{align*}
Hence we have
\begin{align*}
 & \iota_{*,+1} \circ \gamma_\eta^{-1}
([ (T_{\xi(i)} S_{\nu(i)} \otimes 1_n) 
p_i^l (S_{\nu(i)}^*T_{\xi(i)}^*  \otimes 1_n) ] \\
= & \iota_{*,+1}
([S_{\nu(i)'}T_{\bar{\xi}(i)'} \otimes 1_n) 
p_i^l (T_{\bar{\xi}(i)'}^* S_{\nu(i)'}^*\otimes 1_n]) \\
= &
\sum_{b \in \Sigma^\eta}
[(S_{\nu(i)'}T_{\bar{\xi}(i)'b} \otimes 1_n) 
(T_b^* \otimes 1_n) p_i^l (T_b \otimes 1_n)(T_{\bar{\xi}(i)'b}^*S_{\nu(i)'}^* \otimes 1_n)]
\end{align*}

On the other hand,
 we have
$
T_{\xi(i)} S_{\nu(i)} 
= T_{\xi(i)_1}T_{\bar{\xi}(i)} S_{\nu(i)}
= T_{\xi(i)_1}S_{\nu(i)'}T_{\bar{\xi}(i)'}
$
so that
\begin{align*}
 & \iota_{*,+1} 
([ (T_{\xi(i)} S_{\nu(i)} \otimes 1_n) 
p_i^l (S_{\nu(i)}^*T_{\xi(i)}^*  \otimes 1_n) ] \\
 =& \sum_{b \in \Sigma^\eta}
[(T_{\xi(i)_1} S_{\nu(i)'}T_{\bar{\xi}(i)'b} \otimes 1_n) 
(T_b^* \otimes 1_n) p_i^l 
(T_b \otimes 1_n)(T_{\bar{\xi}(i)'b}^*S_{\nu(i)'}^* T_{\xi(i)_1}^*\otimes 1_n)]
\end{align*}
and
hence
\begin{align*}
 & \gamma_\eta^{-1} \circ
 \iota_{*,+1} 
([ (T_{\xi(i)} S_{\nu(i)} \otimes 1_n) 
p_i^l (S_{\nu(i)}^*T_{\xi(i)}^*  \otimes 1_n) ] \\
= & \sum_{b \in \Sigma^\eta}
 \gamma_\eta^{-1}
([(T_{\xi(i)_1} S_{\nu(i)'}T_{\bar{\xi}(i)'b} \otimes 1_n) 
(T_b^* \otimes 1_n) p_i^l 
(T_b \otimes 1_n)(T_{\bar{\xi}(i)'b}^*S_{\nu(i)'}^* T_{\xi(i)_1}^*\otimes 1_n)]) \\
= & \sum_{b \in \Sigma^\eta}
[( S_{\nu(i)'}T_{\bar{\xi}(i)'b} \otimes 1_n) 
(T_b^* \otimes 1_n) p_i^l 
(T_b \otimes 1_n)(T_{\bar{\xi}(i)'b}^*S_{\nu(i)'}^* \otimes 1_n)].
\end{align*}

(ii)  The assertion is completely symmetric to the above proof.
\end{pf}

\begin{lem}
For $k,j$ the following diagrams are commutative.
\begin{enumerate}
\renewcommand{\labelenumi}{(\roman{enumi})}
\item
\begin{equation*}
\begin{CD}
K_0({\cal F}_{\rho,k}) @>\gamma_\eta^{-1}>> K_0({\cal F}_{\rho,k-1}) 
@>\iota_{*,+1}>> K_0({\cal F}_{\rho,k}) \\
@V{\Phi_{\rho,k}}VV   @.  @V{\Phi_{\rho,k}}VV \\
G_\rho   @. \overset{\lambda_{\eta*}}{\longrightarrow} @. G_\rho.
\end{CD}
\end{equation*}
\item
\begin{equation*}
\begin{CD}
K_0({\cal F}_{j,\eta}) @>\gamma_\rho^{-1}>> K_0({\cal F}_{j-1,\eta}) 
@>\iota_{+1,*}>> K_0({\cal F}_{j,\eta}) \\
@V{\Phi_{j,\eta}}VV   @.  @V{\Phi_{j,\eta}}VV \\
G_\eta   @. \overset{\lambda_{\rho*}}{\longrightarrow} @. G_\eta.
\end{CD}
\end{equation*}
\end{enumerate}
\end{lem}
\begin{pf}
(i) 
As in the proof of Lemma 7.8 and Lemma 7.10
one may take an element of
$K_0({\cal F}_{\rho,k})$
as in the following form:
\begin{equation*}
 \sum_{i=1}^{m(l)}\oplus 
[ (T_{\xi(i)} S_{\nu(i)} \otimes 1_n) 
p_i^l (S_{\nu(i)}^*T_{\xi(i)}^*  \otimes 1_n) ]
\in
\bigoplus_{i=1}^{m(l)} K_0({\cal F}_{j,k}(i))
\end{equation*}
for some projection
$p \in M_n(\A)$
and
$j,l$ with $l= j+k$,
where
$
p_i^l 
= (E_i^l \otimes 1) p (E_i^l \otimes 1)\in 
(E_i^l \otimes 1) (\A \otimes M_n) (E_i^l \otimes 1) = M_n(E_i^l \A E_i^l).
$
Keep the notations as in the proof of Lemma 8.10,
we have
\begin{align*}
 & \iota_{*,+1} \circ \gamma_\eta^{-1}
([ (T_{\xi(i)} S_{\nu(i)} \otimes 1_n) 
p_i^l (S_{\nu(i)}^*T_{\xi(i)}^*  \otimes 1_n) ]) \\
= &
\sum_{b \in \Sigma^\eta}
[ (S_{\nu(i)'}T_{\bar{\xi}(i)'b} \otimes 1_n) 
(T_b^* \otimes 1_n) p_i^l (T_b \otimes 1_n)(T_{\bar{\xi}(i)'b}^*S_{\nu(i)'}^* \otimes 1_n)]
\end{align*}
so that
\begin{align*}
 & \Phi_{\rho,k}\circ
\iota_{*,+1} \circ \gamma_\eta^{-1}
([ (T_{\xi(i)} S_{\nu(i)} \otimes 1_n) 
p_i^l (S_{\nu(i)}^*T_{\xi(i)}^*  \otimes 1_n) ] \\
= & 
\sum_{b \in \Sigma^\eta}
\Phi_{\rho,k}\circ(
[S_{\nu(i)'}T_{\bar{\xi}(i)'b} \otimes 1_n) 
(T_b^* \otimes 1_n) p_i^l (T_b \otimes 1_n)(T_{\bar{\xi}(i)'b}^*S_{\nu(i)'}^* \otimes 1_n])) \\
= & 
\sum_{b \in \Sigma^\eta}
[ (T_b^* \otimes 1_n) p_i^l (T_b \otimes 1_n)) \\
= & 
\sum_{b \in \Sigma^\eta}
[ (\eta_b \otimes 1_n) (p_i^l)] \\
= & \lambda_{\eta *}([p_i^l]) \\
= & \lambda_{\eta *}\circ \Phi_{\rho,k} ([ (T_{\xi(i)} S_{\nu(i)} \otimes 1_n) 
p_i^l (S_{\nu(i)}^*T_{\xi(i)}^*  \otimes 1_n) ]). 
\end{align*}
Therefore we have
$\Phi_{\rho,k}\circ \iota_* \circ
\gamma_\eta^{-1} = \lambda_{\eta*} \circ \Phi_{\rho,k}$.

 (ii) The assertion is completely symmetric to the above proof.
\end{pf}


Put
\begin{align*}
G_{\rho,k} & = K_0({\cal F}_{\rho,k}) 
(\cong G_\rho = \lim \{\lambda_{\rho,*}:K_0(\A) \longrightarrow K_0(\A) \}),\\
G_{j,\eta} & = K_0({\cal F}_{j,\eta}) 
(\cong G_\eta = \lim \{\lambda_{\eta,*}:K_0(\A) \longrightarrow K_0(\A) \} ).
\end{align*}

\begin{lem}
\begin{enumerate}
\renewcommand{\labelenumi}{(\roman{enumi})}
The following diagrams are commutative:
\item
\begin{equation*}
\begin{CD}
K_0({\cal F}_{\rho, k}) @>\iota_{*,+1}>> K_0({\cal F}_{\rho,k+1}) \\
@|    @| \\
G_{\rho,k} @>\lambda_{\eta*}>>  G_{\rho,k+1} 
\end{CD}
\end{equation*}
\item
\begin{equation*}
\begin{CD}
K_0({\cal F}_{j,\eta}) @>\iota_{+1,*}>> K_0({\cal F}_{j+1,\eta}) \\
@|    @| \\
G_{j,\eta} @>\lambda_{\rho*}>>  G_{j+1,\eta} 
\end{CD}
\end{equation*}
\end{enumerate}
\end{lem}
Since
\begin{align*}
K_0(\FRE) 
&= \lim_k\{ \iota_{*,+1}: K_0({\cal F}_{\rho, k}) \longrightarrow K_0({\cal F}_{\rho,k+1}) \}\\
&= \lim_j\{ \iota_{+1,*}: K_0({\cal F}_{j,\eta}) \longrightarrow K_0({\cal F}_{j,\eta}) \},
\end{align*}
by putting
$G_{\rho,\eta} = K_0(\FRE) $,
one has
\begin{align*}
G_{\rho,\eta} 
&= \lim_k\{ \lambda_{\eta *}: G_{\rho, k} \longrightarrow G_{\rho,k+1} \}\\
&= \lim_j\{ \lambda_{\rho *}: G_{j,\eta} \longrightarrow G_{j,\eta} \}.
\end{align*}
Define endomorphisms 
\begin{align*}
\sigma_\eta & \text{ on }
G_{\rho,\eta} = \lim_{k}\{ \lambda_{\eta*}:G_{\rho,k}\longrightarrow G_{\rho,k+1} \},\\
\sigma_\rho & \text{ on }
G_{\rho,\eta} = \lim_{j}\{ \lambda_{\eta*}:G_{j,\eta}\longrightarrow G_{j+1,\eta} \}
\end{align*}
by setting
\begin{align*}
\sigma_\rho : & [g,k]\in G_{\rho,k}\longrightarrow [g,k-1]\in G_{\rho,k-1},\\
\sigma_\eta : & [g,j]\in G_{j,\eta}\longrightarrow [g,j-1]\in G_{j-1,\eta}.
\end{align*}
\begin{lem} \hspace{5cm}
\begin{enumerate}
\renewcommand{\labelenumi}{(\roman{enumi})}
\item
There exists an isomorphism
$\Phi_{\rho,\infty}: K_0({\cal F}_{\rho,\eta}) \longrightarrow
G_{\rho,\eta}
$
such that the following diagram is commutative:
\begin{equation*}
\begin{CD}
K_0({\cal F}_{\rho,\eta}) @>\gamma_\eta^{-1}>> K_0({\cal F}_{\rho,\eta}) \\
@V{\Phi_{\rho,\infty}}VV     @V{\Phi_{\rho,\infty}}VV \\
G_{\rho,\eta}
 @>\sigma_\eta>> 
G_{\rho,\eta} 
\end{CD}
\end{equation*}
and hence
\begin{equation*}
\begin{CD}
K_0({\cal F}_{\rho,\eta}) @>id-\gamma_\eta^{-1}>> K_0({\cal F}_{\rho,\eta}) \\
@V{\Phi_{\rho,\infty}}VV     @V{\Phi_{\rho,\infty}}VV \\
G_{\rho,\eta}
 @>id - \sigma_\eta>> 
G_{\rho,\eta}. 
\end{CD}
\end{equation*}
\item
There exists an isomorphism
$\Phi_{\infty,\eta}: K_0({\cal F}_{\rho,\eta})
 \longrightarrow
G_{\rho,\eta} 
$
such that the following diagram is commutative:
\begin{equation*}
\begin{CD}
K_0({\cal F}_{\rho,\eta}) @>\gamma_\rho^{-1}>> K_0({\cal F}_{\rho,\eta}) \\
@V{\Phi_{\infty,\eta}}VV     @V{\Phi_{\infty,\eta}}VV \\
G_{\rho,\eta}
 @>\sigma_\rho>> 
G_{\rho,\eta} 
\end{CD}
\end{equation*}
and hence
\begin{equation*}
\begin{CD}
K_0({\cal F}_{\rho,\eta}) @>id-\gamma_\rho^{-1}>> K_0({\cal F}_{\rho,\eta}) \\
@V{\Phi_{\infty,\eta}}VV     @V{\Phi_{\infty,\eta}}VV \\
G_{\rho,\eta}
 @>id - \sigma_\rho>> 
G_{\rho,\eta}. 
\end{CD}
\end{equation*}
\end{enumerate}
\end{lem}
As
$J_\A: \A = {\cal F}_{0,0}\subset {\cal F}_{\rho,\eta}$
is a subalgebra,
there exists a homomorphism
\begin{equation*}
J_{\A*}: K_0(\A) \longrightarrow K_0({\cal F}_{\rho,\eta}).
\end{equation*}
\begin{lem}
The homomorphism
$
J_{\A*}: K_0(\A) \longrightarrow K_0({\cal F}_{\rho,\eta})
$
is injective such that
\begin{equation*}
J_{A*} \circ \lambda_{\rho *} = \gamma_\rho^{-1} \circ J_{A*}
\quad 
\text{ and }
\quad
J_{A*} \circ \lambda_{\eta *} = \gamma_\eta^{-1} \circ J_{A*}.
\end{equation*}
\end{lem}
\begin{pf}
We will first show that
the endomorphisms
$\lambda_{\rho*},\lambda_{\eta*}$
on $K_0(\A)$ are both injective.
Put
a projection $Q_\alpha = S_\alpha S_\alpha^*$ 
and a subalgebra $\A_\alpha = \rho_\alpha(\A)$ of $\A$ for 
$\alpha \in \Sigma^\rho$.
Then
the endomorphism
$\rho_\alpha$ on $\A$ 
can extend to an isomorphism
from $\A Q_\alpha$ onto $\A_\alpha$
by setting
$\rho_\alpha (x) = S_\alpha^* x S_\alpha, x \in \A Q_\alpha$
whose inverse is 
$\phi_\alpha : \A_\alpha\longrightarrow \A Q_\alpha$
defined by
$\phi_\alpha(y) = S_\alpha y S_\alpha^*, y \in \A_\alpha$.
Hence the induced homomorphism
$\rho_{\alpha*}: K_0(\A Q_\alpha) \longrightarrow K_0(\A_\alpha)$
is an isomorphism.
Since
$\A = \bigoplus_{\alpha \in \Sigma^\rho}Q_\alpha \A$,
the homomorphism
$$
\sum_{\alpha \in \Sigma^\rho}
\phi_{\alpha*}\circ \rho_{\alpha*} : 
K_0(\A) \longrightarrow \bigoplus_{\alpha \in \Sigma^\rho} K_0(Q_\alpha \A)
$$
is an isomorphism,
one may identify
$K_0(\A) =\bigoplus_{\alpha \in \Sigma^\rho} K_0(Q_\alpha \A).
$
Let $g \in K_0(\A)$ satisfy
$\lambda_{\rho*}(g) =0$.
Put
$g_\alpha = \phi_{\alpha*}\circ \rho_{\alpha*}(g) \in  K_0(Q_\alpha \A)$
for
$
\alpha \in \Sigma^\rho 
$
so that
$ g = 
\sum_{\alpha \in \Sigma^\rho} g_\alpha.
$
As 
$
\rho_{\beta*} \circ \phi_{\alpha*} = 0
$
for
$\beta \ne \alpha$,
one sees 
$
\rho_{\beta*}(g_\alpha) = 0
$
for
$\beta \ne \alpha$.
Hence we have
$$
0=
\lambda_{\rho*}(g) 
= 
\sum_{\beta \in \Sigma^\rho}\sum_{\alpha \in \Sigma^\rho} \rho_{\beta*}(g_\alpha)
= 
\sum_{\alpha \in \Sigma^\rho} \rho_{\alpha*}(g_\alpha)
\in \bigoplus_{\alpha \in \Sigma^\rho} K_0(\A_\alpha).
$$
It follows that
$\rho_{\alpha*}(g_\alpha) =0$
in $K_0(\A_\alpha)$.
Since
$\rho_{\alpha*}:K_0(Q_\alpha \A)\longrightarrow K_0(\A_\alpha)$
is isomorophic,
one sees that
$g_\alpha =0$
in $K_0(\A Q_\alpha )$
for all $\alpha \in \Sigma^\rho$.
This implies that
$g = 
\sum_{\alpha \in \Sigma^\rho} g_\alpha =0 $
in
$K_0(\A)$.
Therefore we know that
the endomorphism
$\lambda_{\rho*}$
on $K_0(\A)$ is injective.
and similarly 
so is 
$\lambda_{\eta*}$.

By the previous lemma,
there exists an isomorphism
 $
\Phi_{j,k} : K_0({\cal F}_{j,k}) \longrightarrow K_0(\A)
$
such that the following diagram 
\begin{equation*}
\begin{CD}
K_0({\cal F}_{j,k}) @>\iota_{+1,*}>> K_0({\cal F}_{j+1,k})   \\
@V{\Phi_{j,k}}VV   @V{\Phi_{j+1,k}}VV \\
K_0(\A) @>\lambda_{\rho*}>> K_0(\A) 
\end{CD}
\end{equation*}
is commutative
so that
the embedding
$
\iota_{+1,*} : K_0({\cal F}_{j,k}) \longrightarrow K_0({\cal F}_{j+1,k})
$
is injective,
and similarly
$
\iota_{*,+1} : K_0({\cal F}_{j,k}) \longrightarrow K_0({\cal F}_{j,k+1})
$
is injective.
Hence for $n,m \in {\Bbb N}$,
the homomorphism 
$$
\iota_{n,m}: K_0(\A)=K_0({\cal F}_{0,0}) \longrightarrow K_0({\cal F}_{n,m})
$$
defined by the compositions of 
$ 
\iota_{+1,*}
$
and
$ 
\iota_{*,+1}
$
is injective.
By \cite[Theorem 6.3.2 (iii)]{RLL},
one knows 
$\Ker(J_{\A*}) = \cup_{n,m\in {\Bbb N}} \Ker(\iota_{n,m})$,
so that
$\Ker(J_{\A*}) =0$.
\end{pf}
We henceforth identify the group
$K_0(\A)$ with its imaga 
$J_{A*}(K_0(\A))$ in $K_0(\FRE)$.
As in the above proof,
not only $K_0(\A)(=K_0({\cal F}_{0,0}))$
but also the groups
$K_0({\cal F}_{j,k})$
for $j,k$ are identified with subgroups 
of $K_0(\FRE)$ 
via injective homomorphisms
from 
$K_0({\cal F}_{j,k})$
to
$K_0(\FRE)$
induced by 
the embedding
of ${\cal F}_{j,k}$ 
into
$\FRE$.

We note that 
\begin{align*}
(\id -\gamma_\eta)(K_0({\cal F}_{\rho,\eta})) 
& =
(\id -\gamma_\eta^{-1})(K_0({\cal F}_{\rho,\eta})), \\
(\id -\gamma_\rho)(K_0({\cal F}_{\rho,\eta})) 
& =
(\id -\gamma_\rho^{-1})(K_0({\cal F}_{\rho,\eta}))
\end{align*}
and
\begin{equation*}
(\id -\gamma_\rho) \cap(\id -\gamma_\eta) \text{ in } K_0({\cal F}_{\rho,\eta}) 
=
(\id -\gamma_\rho^{-1}) \cap(\id -\gamma_\eta^{-1}) \text{ in } K_0({\cal F}_{\rho,\eta}). 
\end{equation*}
Denote by
$
(\id -\gamma_\rho)(K_0(\FRE)) + (\id -\gamma_\eta)(K_0(\FRE)) 
$
the subgroup of 
$K_0(\FRE)$
generated by 
$
(\id -\gamma_\rho)(K_0(\FRE))
$
and
$ (\id -\gamma_\eta)(K_0(\FRE)). 
$
\begin{lem}
An element in $K_0(\FRE)$ is equivalent 
to an element of 
$K_0(\A)$
modulo 
the subgroup
$
(\id -\gamma_\rho)(K_0(\FRE)) + (\id -\gamma_\eta)(K_0(\FRE)). 
$
\end{lem}
\begin{pf}
For
$ g \in K_0(\FRE)$,
we may assume that
$g \in K_0({\cal F}_{j,k})$ 
for some $j,k \in \Zp$.
As  
$\gamma_\rho^{-1}$ commutes with
$\gamma_\eta^{-1}$,
one sees that
$ (\gamma_\rho^{-1})^j \circ (\gamma_\eta^{-1})^k(g) \in K_0(\A)$.
Put
$g_1 = \gamma_\rho^{-1}(g)$
so that 
\begin{equation*}
g - (\gamma_\rho^{-1})^j \circ (\gamma_\eta^{-1})^k(g)
= g - \gamma_\rho^{-1}(g)
+ g_1 - (\gamma_\rho^{-1})^{j-1} \circ (\gamma_\eta^{-1})^k(g_1).
\end{equation*}
We inductively see that 
$
g - (\gamma_\rho^{-1})^j \circ (\gamma_\eta^{-1})^k(g)
$
belongs to the subgroup 
$
 (\id -\gamma_\rho)(K_0(\FRE))
+
  (\id -\gamma_\eta)(K_0(\FRE)). 
$
Hence 
$g$ 
is equivalent 
to 
$
(\gamma_\rho^{-1})^j \circ (\gamma_\eta^{-1})^k(g)
$ 
modulo 
$
(\id -\gamma_\rho)(K_0(\FRE)) + (\id -\gamma_\eta)(K_0(\FRE)).
$
\end{pf}
Denote by
$
(\id -\lambda_{\rho*})(K_0(\A)) + (\id -\lambda_{\eta*})(K_0(\A)) 
$
the subgroup of 
$ K_0(\A) $
generated by 
$
(\id -\lambda_{\rho*})(K_0(\A))
$
and
$(\id -\lambda_{\eta*})(K_0(\A)).$ 
\begin{lem}
For
$ g \in K_0(\A)$,
the condition
$
g \in
(\id -\gamma_\rho^{-1})(K_0(\FRE))+ (\id -\gamma_\eta^{-1})(K_0(\FRE)) 
$
implies
$
g \in
(\id -\lambda_{\rho*})(K_0(\A)) +(\id -\lambda_{\eta*})(K_0(\A)). 
$ 
\end{lem}
\begin{pf}
By the assumption that
$
g \in  
(\id -\gamma_\rho^{-1})(K_0(\FRE)) + (\id -\gamma_\eta^{-1})(K_0(\FRE)) 
$
there exist
$
g_1  \in
(\id -\gamma_\rho^{-1})(K_0(\FRE))
$
and
$
g_2 \in (\id -\gamma_\eta^{-1})(K_0(\FRE))
$
such that
$
g = g_1 + g_2,
$
where
$g_1 = (\id -\gamma_\rho^{-1})(h_1)$
and
$g_2 = (\id -\gamma_\eta^{-1})(h_2)$
for some
$h_1, h_2 \in K_0(\FRE)$.
We may assume that
$h_1, h_2 \in K_0({\cal F}_{j,k})$
for large enough $j,k \in \Zp$.
Put
$
e_i =(\gamma_\rho^{-1})^j \circ (\gamma_\eta^{-1})^k(h_i )
$ 
which belongs to $K_0({\cal F}_{0,0})(=K_0(\A))$
for
$i=0,1$.
It follows that
$$
\lambda_\rho^j \circ \lambda_\eta^k(g) 
= 
(\id -\lambda_{\eta*})(e_1) + (\id -\lambda_{\rho*})(e_2).
$$
Now 
$
g \in K_0(\A)$ and 
$
\lambda_\rho^j \circ \lambda_\eta^k(g) \in 
(\id -\lambda_{\eta*})(K_0(\A))
+ (\id -\lambda_{\rho*})(K_0(\A)) \subset K_0(\A).
$
As in the proof of the preceding lemma,
by putting
$
g^{(n)} = \lambda_{\rho*}^n(g),
g^{(n,m)} =\lambda_{\eta*}^m(g^{(n)}) \in K_0(\A)
$
we have
\begin{align*}
 & g - \lambda_{\rho*}^j\circ \lambda_{\eta*}^k(g) \\
=& g - \lambda_{\rho*}(g) 
+ g^{(1)} - \lambda_{\rho*}(g^{(1)}) 
+ g^{(2)} - \lambda_{\rho*}(g^{(2)}) + \cdots
+ g^{(j-1)} - \lambda_{\rho*}(g^{(j-1)}) \\
 & + g^{(j)} -  \lambda_{\eta*}^k(g^{(j)}) \\
=& g - \lambda_{\rho*}(g) 
+ g^{(1)} - \lambda_{\rho*}(g^{(1)}) 
+ g^{(2)} - \lambda_{\rho*}(g^{(2)}) + \cdots
+ g^{(j-1)} - \lambda_{\rho*}(g^{(j-1)})\\
& + g^{(j)} - \lambda_{\eta*}(g^{(j)}) 
+ g^{(j,1)} - \lambda_{\eta*}(g^{(j,1)}) 
+ g^{(j,2)} - \lambda_{\eta*}(g^{(j,2)}) + \cdots \\
& + g^{(j,k-1)} - \lambda_{\eta*}(g^{(j,k-1)})\\
=& (\id - \lambda_{\rho*})(g + g^{(1)} + \cdots + g^{(j-1)})
 + (\id - \lambda_{\eta*})(g^{(j)} + g^{(j,1)} + \cdots + g^{(j,k-1)}) 
\end{align*}
Since
$ 
\lambda_{\rho*}^j\circ \lambda_{\eta*}^k(g) \in 
(\id -\lambda_{\eta*})(K_0(\A))
+ (\id -\lambda_{\rho*})(K_0(\A)) 
$
and
\begin{align*}
(\id & - \lambda_{\rho*})(g + g^{(1)} + \cdots + g^{(j-1)})
 \in  (\id -\lambda_{\rho*})(K_0(\A)),\\
(\id & - \lambda_{\eta*})(g^{(j)} + g^{(j,1)} + \cdots + g^{(j,k-1)}) 
\in (\id -\lambda_{\eta*})(K_0(\A)),
\end{align*}
we have 
 $$
 g \in 
(\id -\lambda_{\eta*})(K_0(\A))
+ (\id -\lambda_{\rho*})(K_0(\A)). 
$$
\end{pf}
Hence we obtain the following lemma for the cokernel.
\begin{lem}
The quotient group
$$
 K_0(\FRE) / 
((\id - \gamma_\eta^{-1})(K_0(\FRE)) + (\id - \gamma_\rho^{-1})(K_0(\FRE)))
$$ 
is isomorphic to the quotient group
$$
K_0(\A) /
((\id - \lambda_{\eta*})(K_0(\A)) + (\id - \lambda_{\rho*})(K_0(\A)))
$$
\end{lem}
\begin{pf}
Surjectivity of the quotient map
$$
q_{\A*}: K_0(\A) \longrightarrow 
K_0(\FRE) / 
((\id - \gamma_\eta^{-1})(K_0(\FRE))+(\id - \gamma_\rho^{-1})(K_0(\FRE)))
$$
comes from the preceding lemma.
As
$$
\Ker(q_{\A*}) = (\id - \lambda_{\eta*})(K_0(\A)) + (\id - \lambda_{\rho*})(K_0(\A))
$$
by the preceing lemma,
we have a desired assertion.
\end{pf}
For the kernel, we have 
\begin{lem}
The subgroup
$$
\Ker(\id - \gamma_\eta^{-1}) \cap \Ker(\id - \gamma_\rho^{-1})
\text{ in } 
K_0(\FRE)
$$
is isomorphic to the subgroup
$$
\Ker(\id - \lambda_{\eta*}) \cap \Ker(\id - \lambda_{\rho*})
\text{ in } 
K_0(\A)
$$
through 
$J_{\A*}$.
\end{lem}
\begin{pf}
For $g \in \Ker(\id - \gamma_\eta^{-1}) \cap \Ker(\id - \gamma_\rho^{-1})
$
in 
$K_0(\FRE)$,
one may assume that
$g \in K_0({\cal F}_{j,k})$ 
for some $j,k \in \Zp$
so that
$
(\gamma_\rho^{-1})^j\circ (\gamma_\eta^{-1})^k(g) \in K_0(\A).
$
Through the identification
between 
$J_{\A*}(K_0(\A))$
and
$K_0(\A)$
via $J_{\A*}$,
one may assume that
$\lambda_{\eta*} = \gamma_\eta^{-1}$
and
$
\lambda_{\rho*} = \gamma_\rho^{-1}
$
on
$K_0(\A)$.
As
$g \in \Ker(\id - \gamma_\eta^{-1}) \cap \Ker(\id - \gamma_\rho^{-1}),
$
one has
$
g = (\gamma_\rho^{-1})^j\circ (\gamma_\eta^{-1})^k(g) \in K_0(\A).
$
This implies that
$
g \in \Ker(\id - \lambda_{\eta*}) \cap \Ker(\id - \lambda_{\rho*})
$
in
$K_0(\A)$.
The converse inclusion relation
$
\Ker(\id - \lambda_{\eta*}) \cap \Ker(\id - \lambda_{\rho*})
\subset
\Ker(\id - \gamma_\eta^{-1}) \cap \Ker(\id - \gamma_\rho^{-1})
$
is clear through the above identification.
\end{pf}
Therefore we have
\begin{prop}
There exists a short  exact sequence: 
\begin{align*}
0 
& \longrightarrow 
K_0(\A) /
((\id - \lambda_{\eta*})(K_0(\A)) + (\id - \lambda_{\rho*})(K_0(\A))) \\
& \longrightarrow
K_0(\ORE) \\
& \longrightarrow
\Ker(\id - \lambda_{\eta*}) \cap \Ker(\id - \lambda_{\rho*})
\text{ in }
K_0(\A)
\longrightarrow 0.
\end{align*}
\end{prop}

Let
${\cal F}_\rho$ be the fixed point algebra
$({\cal O}_\rho)^{\hat{\rho}}$
of the $C^*$-symbolic dynamical system
$(\A,\rho,\Sigma^\rho)$.
The algebra ${\cal F}_\rho$
is isomorphic to
the subalgebra ${\cal F}_{\rho,0}$ of $\FRE$
in a natural way.
As in the proof of Lemma 7.14, the group
$K_0({\cal F}_{\rho,0})$
is regarded as a subgroup of 
$K_0(\FRE)$ and
the restriction of $\gamma_\eta^{-1}$
to $K_0({\cal F}_{\rho,0})$
satisfies 
$\gamma_\eta^{-1}(K_0({\cal F}_{\rho,0})) \subset K_0({\cal F}_{\rho,0})$
so that
$\gamma_\eta^{-1}$
yields an endomorphism on ${\cal F}_\rho$,
which we denote by $\gamma_\eta^{-1}$.

For the group $K_1(\ORE)$,
we provide several lemmas.
Their proofs are similarly to the above discussions.
\begin{lem}\hspace{5cm}
\begin{enumerate}
\renewcommand{\labelenumi}{(\roman{enumi})}
\item
An element in $K_0(\FRE)$ is equivalent 
to an element of
$K_0({\cal F}_{\rho,0})( = K_0({\cal F}_{\rho}))$
modulo the subgroup
$(\id - \gamma_\eta)(K_0(\FRE))$. 
\item 
If $g \in K_0({\cal F}_\rho)(=K_0({\cal F}_{\rho,0}))$
belongs to 
$(\id - \gamma_\eta)(K_0(\FRE))$,
then
$g $ belongs to
$(\id - \gamma_\eta)(K_0({\cal F}_{\rho}))$.
\end{enumerate}
\end{lem} 

\begin{lem}
The quotient group
$
K_0(\FRE) / (\id - \gamma_\eta^{-1})(K_0(\FRE))
$
is isomorphic to the quotient group
$
K_0({\cal F}_\rho) / (\id - \gamma_\eta^{-1})(K_0({\cal F}_\rho)),
$
that is also isomorphic to the quotient group
$
K_0(\A) / (1 - \lambda_\eta)K_0(\A)
$
such that the kernel of 
$\id - \gamma_\rho$ in
$
K_0(\FRE) / (\id - \gamma_\eta^{-1})(K_0(\FRE))
$
is isomorphic to
the kernel of 
$
\id - \lambda_\rho
$
 in
$
K_0(\A) / (1 - \lambda_\eta)K_0(\A).
$ 
That is
\begin{align*}
& \Ker(\id - \gamma_\rho) \text{ in } 
K_0(\FRE) / (\id - \gamma_\eta^{-1})(K_0(\FRE))  \\
\cong  &
\Ker(\id - \lambda_{\rho*})\text{ in } 
K_0(\A) / (\id - \lambda_{\eta*})(K_0(\A)).
\end{align*}
\end{lem}
 \begin{pf}
 The first assertion that
 the three quotient groups 
$$
K_0(\FRE) / (\id - \gamma_\eta^{-1})(K_0(\FRE)),\quad 
K_0({\cal F}_\rho) / (\id - \gamma_\eta^{-1})(K_0({\cal F}_\rho)),\quad
K_0(\A) / (1 - \lambda_\eta)K_0(\A)
$$
are naturally isomorphic 
is similarly proved to the previous discussions. 
 For the second assertion,
the kernel
$ \Ker(\id - \gamma_\rho)$
 in 
$ 
K_0(\FRE) / (\id - \gamma_\eta^{-1})(K_0(\FRE))
$
is isomorphic to
the kernel 
$ \Ker(\id - \gamma_\rho)$
 in 
$ 
K_0({\cal F}_\rho) / (\id - \gamma_\eta^{-1})(K_0({\cal F}_\rho))
$
which is isomorphic to
the kernel
$  
\Ker(\id - \lambda_{\rho*})
$
 in
$
K_0(\A) / (1 - \lambda_{\eta*})(K_0(\A)).
$ 
 \end{pf}

\begin{lem}
The kernel of
$\id - \gamma_\rho$ in
$
K_0(\FRE)
$
is isomorphic to
the kernel of
$\id - \gamma_\rho$ 
in
$
K_0({\cal F}_\rho)
$
that is also isomorphic to the kernel of
$
\id - \lambda_{\eta*}
$
in
$K_0(\A)$
such that 
the quotient group
$$
\Ker(\id-\gamma_\eta) 
\text{ in }
K_0(\FRE) / 
(\id - \gamma_\rho) 
(\Ker(\id-\gamma_\eta) 
\text{ in }
K_0(\FRE)) 
$$
 is isomorphic to
the quotient group
$$
\Ker(\id-\lambda_{\eta*}) 
\text{ in }
K_0(\A) / 
(\id  - \lambda_{\rho*}) 
(\Ker(\id-\lambda_{\eta*}) 
\text{ in }
K_0(\A). 
$$
That is
\begin{align*}
& \Ker(\id-\gamma_\eta) 
\text{ in }
K_0(\FRE) / 
(\id - \gamma_\rho) 
(\Ker(\id-\gamma_\eta) 
\text{ in }
K_0(\FRE)) \\
\cong
&
\Ker(\id-\lambda_{\eta*}) 
\text{ in }
K_0(\A) / 
(\id  - \lambda_{\rho*}) 
(\Ker(\id-\lambda_{\eta*}) 
\text{ in }
K_0(\A). 
\end{align*}
\end{lem}
\begin{pf}
The proofs are similar to the previous discussions.
\end{pf}
Therefore we have
\begin{prop}
There exists a short  exact sequence: 
\begin{align*}
0 
& \longrightarrow 
\Ker(\id-\lambda_{\eta*}) 
\text{ in }
K_0(\A) / 
(\id - \lambda_{\rho*}) 
(\Ker(\id-\lambda_{\eta*}) 
\text{ in }
K_0(\A)) \\
& \longrightarrow
K_1(\ORE) \\
& \longrightarrow
\Ker(\id - \lambda_{\rho*})\text{ in }
K_0(\A) / (\id - \lambda_{\eta*})(K_0(\A))
\longrightarrow 0.
\end{align*}
\end{prop}
Consequently we have
\begin{thm}
Suppose that a $C^*$-textile dynamical system 
$\CTDS$ forms square.
Then there exist short exact sequences for their K-theory groups
as in the following way:
\begin{align*}
0 
& \longrightarrow 
K_0(\A) /
(\id - \lambda_{\eta*})K_0(\A) + (\id - \lambda_{\rho*})K_0(\A) \\
& \longrightarrow
K_0(\ORE) \\
& \longrightarrow
\Ker(\id - \lambda_{\eta*}) \cap \Ker(\id - \lambda_{\rho*})
\text{ in }
K_0(\A)
\longrightarrow 0 \\
\intertext{and}
0 
& \longrightarrow 
\Ker(\id-\lambda_{\eta*}) 
\text{ in }
K_0(\A) / 
(\id - \lambda_{\rho*}) 
(\Ker(\id-\lambda_{\eta*}) 
\text{ in }
K_0(\A)) \\
& \longrightarrow
K_1(\ORE) \\
& \longrightarrow
\Ker(\id - \lambda_{\rho*})\text{ in }
K_0(\A) / (\id - \lambda_{\eta*})K_0(\A)
\longrightarrow 0
\end{align*}
where the endomorphisms
$
\lambda_\rho, \lambda_\eta: 
K_0(\A) \longrightarrow K_0(\A)
$
are defined by
\begin{align*}
\lambda_{\rho*}([p]) & = \sum_{\alpha \in \Sigma^\rho} [\rho_{\alpha*}(p)] \in K_0(\A)
\text{ for } [p] \in K_0(\A),\\
\lambda_{\eta*}([p]) & = \sum_{a \in \Sigma^\eta} [\eta_{a*}(p)] \in K_0(\A)
\text{ for } [p] \in K_0(\A).  
\end{align*}
\end{thm}

\section{Examples}

{\bf 1.  LR-textile $\lambda$-graph systems.}

A symbolic matrix $\M = [\M(i,j)]_{i,j=1}^N$
 is a matrix whose components consist of formal sums of elements of $\Sigma$,
such as 
\begin{equation*}
\M =
\begin{bmatrix}
a & a + c \\
c & 0
\end{bmatrix}
\qquad 
\text{ where }
\Sigma = \{ a,b,c\}.
\end{equation*}
$\M$ is said to be essential if
there is no zero column or zero row.
$\M$ is said to be 
left-resolving if for each column 
a symbol does not appear in two different rows.
For example,
$
\begin{bmatrix}
a & a + b \\
c & 0
\end{bmatrix}
$
is left-resolving,
but
$
\begin{bmatrix}
a & a + b \\
c & b
\end{bmatrix}
$
is not left-resolving
because of $b$ at the second column.
We henceforth asssume that
symbolic matrices are always essential and left-resolving.

Let
$\M = [\M(i,j)]_{i,j=1}^N$
and
$\M' = [\M'(i,j)]_{i,j=1}^N$
be symbolic matrices over
$\Sigma$ and $\Sigma'$ respectively.
Suppose that there is a bijection
$\kappa:\Sigma \longrightarrow \Sigma'$.
Following Nasu's terminology \cite{NaMemoir}
we say that 
$\M$ and $\M'$ are equivalent under specification $\kappa$,
or simply a specified equivalence 
if $\M'$ can be obtained from $\M$ by replacing every symbol
$\alpha \in \Sigma$ by
$\kappa(\alpha)$.
That is
if
$\M(i,j) = \alpha_1+\cdots+\alpha_n$,
then
$\M'(i,j) = \kappa(\alpha_1)+\cdots+ \kappa(\alpha_n)$.
We write this situation as
$\M \overset{\kappa}{\cong}\M'$ (see \cite{NaMemoir}).

For
a symbolic matrix $\M = [\M(i,j)]_{i,j=1}^N$
over $\Sigma^{\M}$,
we set 
for $\alpha \in \Sigma^{\M}, i,j=1,\dots,N$ 
\begin{equation}
A^\M(i,\alpha,j) =
\begin{cases}
1 & \text{ if } \alpha \text{ appears in } \M(i,j),\\
0 & \text{ otherwise}.
\end{cases}
\end{equation}
Put an $N\times N$ nonnegative matrix
$A^\M = [A_\M(i,j)]_{i,j=1}^N$
by setting
$
A^\M(i,j) = \sum_{\alpha \in\Sigma^{\M}} A^\M(i,\alpha,j).
$
Let
$\A$ be an $N$-dimensional commutative $C^*$-algebra ${\Bbb C}^N$
with minimal projections 
$E_1,\dots, E_n$
such that
$$
\A = {\Bbb C}E_1 \oplus \cdots \oplus {\Bbb C}E_n.
$$
We set for $\alpha \in \Sigma^{\M}$:
\begin{equation*}
\rho^{\M}_\alpha(E_i) = \sum_{j=1}^N A^\M(i,\alpha,j)E_j,
\qquad i=1,\dots,n.
\end{equation*}
Then we have a  $C^*$-symbolic dynamical system
$(\A,\rho^{\M},\Sigma^{\M})$.

Let
$\M = [\M(i,j)]_{i,j=1}^N$
and
$\N = [\N(i,j)]_{i,j=1}^N$
be symbolic matrices over
$\Sigma^{\M}$ and $\Sigma^{\N}$ 
respectively.
Suppose that there is a bijection
$\kappa:\Sigma^{\M} \longrightarrow \Sigma^{\N}$
such that 
$\kappa$ yields a specified equivalence
\begin{equation}
\M \N \overset{\kappa}{\cong} \N \M. \label{eqn:LR}
\end{equation}
Then we have two $C^*$-symbolic dynamical systems
$(\A,\rho^{\M},\Sigma^{\M})$
and
$(\A,\rho^{\N},\Sigma^{\N})$.  
Put 
\begin{align*}
\Sigma^{\M\N}
& = \{ (\alpha,b) \in \Sigma^{\M}\times\Sigma^{\N}
    \mid \rho^{\N}_b \circ \rho^{\M}_\alpha \ne 0 \}, \\ 
\Sigma^{\N\M}
& = \{ (a,\beta)\in \Sigma^{\N}\times\Sigma^{\M}
    \mid \rho^{\M}_\beta \circ \rho^{\N}_a \ne 0 \}.
\end{align*}
\begin{prop}
Keep the above situations.
$\kappa$ induces a specification 
$\kappa:\Sigma^{\M\N} \longrightarrow \Sigma^{\N\M}$
such that
\begin{equation}
\rho^{\N}_b \circ \rho^{\M}_\alpha
=
\rho^{\M}_\beta \circ \rho^{\N}_a
\quad
\text{ if }
\quad
\kappa(\alpha,b) = (a,\beta).
\end{equation}
Hence 
$(\A, \rho^{\M}, \rho^{\N},\Sigma^\M, \Sigma^\N,\kappa)$
yields a $C^*$-textile dynamical system.
\end{prop}
\begin{pf}
Since
$\M \N \overset{\kappa}{\cong} \N \M,
$
one sees that for $i,j=1,2,\dots,N$,
$\kappa(\M\N(i,j)) = \N\M(i,j).$
For
$(\alpha,b) \in \Sigma^{\M\N}$,
there exists $i,k=1,2,\dots,N$ such that
$\rho^\N_b \circ \rho^\M_\alpha(E_i) \ge E_k$.
As
$\kappa(\alpha,b)$ 
appears in 
$\N\M(i,k)$,
by putting
$(a,\beta) = \kappa(\alpha,b)$,
we have
$\rho^\M_\beta \circ \rho^\N_a(E_i) \ge E_k$.
Hence
$\kappa(\alpha,b) \in \Sigma^{\N\M}$.
One indeed sees that
$
\rho^{\N}_b \circ \rho^{\M}_\alpha
=
\rho^{\M}_\beta \circ \rho^{\N}_a
$
by the relation
$\M \N \overset{\kappa}{\cong} \N \M.
$
\end{pf}

We remark that symbolic matrices are presentations of labeled directed graphs.
Hence we may consider our discussions above in terms of labeled directed graphs.

Two symbolic matrices satisfying the relations \eqref{eqn:LR}
gives rise to 
LR textile systems that have been introduced by Nasu
(see \cite{NaMemoir}).
Textile systems introduced by Nasu play a important tool to analyze automorphisms
and endomorphisms of topological Markov shifts,

The author has generalized the LR-textile systems to LR-textile $\lambda$-graph systems
which consists of two commuting symbolic matrix systems (\cite{MaYMJ2008}). 
Let
$\T$
be an LR-textile $\lambda$-graph system defined by a specified equivalence:
\begin{equation}
\M_{l,l+1} \N_{l+1,l+2} \overset{\kappa}{\cong}
\N_{l,l+1} \M_{l+1,l+2},
\qquad
l \in \Zp
\end{equation}
through specification $\kappa$.
There exist two symbolic matrix systems
$(\M,I^\M)$ and
$(\N,I^\N)$.
Denote by
${\frak L}^\M$ 
and
${\frak L}^\N$
the associated $\lambda$-graph systems respectively. 
Since
${\frak L}^\M$ 
and
${\frak L}^\N$
form square in the sense of \cite[p.170]{MaYMJ2008},
they have common sequences 
$V_l^\M = V_l^\N, l \in \Zp$
of vertices and inclusion matrices
$I^\M_{l,l+1} =I^\N_{l,l+1}, l \in \Zp$.
We denote 
$V_l^\M = V_l^\N$
and
$I^\M_{l,l+1} =I^\N_{l,l+1}$
by  $V_l$ and $I_{l,l+1}$
respectively.

Let
$(\A_\M, \rho^\M,\Sigma^\M)$
and
$(\A_\N, \rho^\N,\Sigma^\N)$
be the associated $C^*$-symbolic dynamical systems with the $\lambda$-graph systems
${\frak L}^\M$ 
and
${\frak L}^\N$
respectively.
Hence one sees that
$\A_\M = \A_\N$ 
which is denoted by
$\A$.
It is easy to seee that
the relation (9.1) implies 
\begin{equation}
\rho^\M_\alpha \circ \rho^\N_b
=
\rho^\N_a \circ \rho^\M_\beta
\qquad
\text{ if }
\quad
\kappa(\alpha,b) = (a,\beta).
\end{equation}
\begin{prop}
An LR-textile $\lambda$-graph system $\T$
yields a $C^*$-textile dynamical system
$(\A, \rho^\M, \rho^\N, \Sigma^\M, \Sigma^\N,\kappa)$
which forms square.

Conversely, a  $C^*$-textile dynamical system
$\CTDS$ which forms square 
yields an LR-textile $\lambda$-graph system
${\cal T}_{{\K}^{\M^\rho}_{\M^\eta}}
$
such that the associated 
 $C^*$-textile dynamical system
$(\A, \rho^{\M^\rho}, \rho^{\M^\eta}, \Sigma^{\M^\rho}, \Sigma^{\M^\rho},\kappa)$
is a subsystem of
$\CTDS$ in the sense which satisfis the following relations:
$$
\A_{\frak L} \subset \A, \qquad \rho |_{\A^{\frak L}} = \rho^{\M^\rho},
\qquad
\eta |_{\A^{\frak L}} = \rho^{\M^\eta}.
$$
\end{prop}
\begin{pf}
Let $\T$ be an LR-textile $\lambda$-graph system.
As in the above discussions,
we have
a $C^*$-textile dynamical system
$(\A, \rho^\M, \rho^\N, \Sigma^\M, \Sigma^\N,\kappa)$.

Conversely,
let
$\CTDS$ be a  $C^*$-textile dynamical system
 which forms square.
 Put for $l \in {\Bbb N}$
\begin{equation*}
\A_l^\rho = C^*(\rho_\mu(1): \mu \in B_l(\Lambda_\rho)),
\qquad
\A_l^\eta = C^*(\eta_\xi(1): \xi \in B_l(\Lambda_\eta)).
\end{equation*}
Since 
$\A_l^\rho =\A_l^\eta$
and they are commutative and of finite dimensional,
the algebra
$$
\A_{\frak L} = \overline{\cup_{l \in \Zp}\A_l^\rho}
             = \overline{\cup_{l \in \Zp}\A_l^\eta}
$$
is a commutative AF-subalgebra of $\A$.
It is easy to see that both
$(\A,\rho,\Sigma^\rho)$
and
$(\A,\eta,\Sigma^\eta)$
are $C^*$-symbolic dynamical systems
such that 
\begin{equation}
\eta_b \circ \rho_\alpha = 
\rho_\beta \circ \eta_a
\qquad
\text{ if }
\quad
\kappa(\alpha,b) = (a,\beta) \label{sharp:eqn}
\end{equation}
By construction,
there exist $\lambda$-graph systems
${\frak L}^\rho$
and
${\frak L}^\eta$
whose $C^*$-symbolic dynamical systems 
are
$(\A,\rho,\Sigma^\rho)$
and
$(\A,\eta,\Sigma^\eta)$
respectively.
Let
$(\M^\rho, I^\rho)$
and
$(\M^\eta, I^\eta)$
be the associated symbolic dynamical systems.
It is easy to see that 
the relation \eqref{sharp}
implies
\begin{equation}
\M^\rho_{l,l+1} \M^\eta_{l+1,l+2} \overset{\kappa}{\cong}
\M^\eta_{l,l+1} \M^\rho_{l+1,l+2},
\qquad
l \in \Zp.
\end{equation}
Hence we have an LR-textile $\lambda$-graph system
$
{\cal T}_{{\cal K}^{\M^\rho}_{\M^\eta}}.
$
It is direct to see that the associated $C^*$-textile dynamical system
is
$$
(\A_{\frak L}, \rho|_{\A_{\frak L}}, \eta|_{\A_{\frak L}}, \Sigma^\rho, \Sigma^\eta,\kappa).
$$
\end{pf}
\medskip

Let $A$ be an $N \times N$ matrix with entries in nonnegative integers.
We may consider a directed graph
$G_A = (V_A,E_A)$ with vertex set $V_A$ and edge set $E_A$.
The vertex set 
$V_A$ consists of $N$ vertices which we denote by
$\{v_1,\dots,v_N \}$.
We equip $A(i,j)$ edges from the vertex
$v_i$ to the vertex $v_j$.
Denote by
$E_A$ such edges.
Let
$\Sigma^A = E_A$
and the labeling map
$\lambda_A: E_A \longrightarrow \Sigma^A$
be defined as the identity map.
Then we have a labeled directed graph denoted by $G_A$
as well as a symbolic matrix
$\M_A =[\M_A(i,j)]_{i,j=1}^N$ by setting
\begin{equation*}
\M_A(i,j) =
\begin{cases}
e_1+ \cdots + e_n & \text{ if } e_1,\cdots,e_n 
\text{ are edges from } v_i \text{ to } v_j, \\
0 & \text{ if there is no edge from } v_i \text{ to } v_j.
\end{cases}
\end{equation*}
Let
$B$ be an $N\times N$ matrix with entries in nonnegative integers such that
\begin{equation*}
AB = BA.
\end{equation*} 
Hence the numbers of pairs of directed edges
\begin{align*}
\{(e,f)\in E_A \times E_B 
& 
\mid s(e) = v_i, t(e) = s(f), t(f) = v_k \}\\ 
\{(f,e)\in E_B \times E_A 
& 
\mid s(f) = v_i, t(f) = s(e), t(e) = v_k \}.
\end{align*}
coincide with each other for each
$v_i$ and $v_k$, 
so that one may take a bijection
$\kappa: \Sigma^{AB} \longrightarrow \Sigma^{BA}$
which gives rise to
a specified equivalence
$\M_A \M_B \overset{\kappa}{\cong}\M_A \M_B$.
Therefore we have a $C^*$-textile dynamical system
\begin{equation*}
(\A, \rho^{\M_A}, \rho^{\M_B},\Sigma^A, \Sigma^B, \kappa)
\end{equation*}
which we denote by
\begin{equation*}
(\A, \rho^{A}, \rho^{B},\Sigma^A, \Sigma^B, \kappa).
\end{equation*}
The associated $C^*$-algebra
is denoted by
${\cal O}_{A,B}^\kappa$.
We remark that
the algebra ${\cal O}_{A,B}^\kappa$
is dependent on the choice of specification
$\kappa: \Sigma^{AB} \longrightarrow \Sigma^{BA}$.
The algebras are $2$-graph algebras of Kumjian and Pask \cite{KP}.
They are $C^*$-algebras associated to textile systems studied by V. Deaconu \cite{Dea}.

\begin{prop}
Keep the above situations.
There exist short exact sequences: 
\begin{align*}
0 
& \longrightarrow 
{\Bbb Z}^N /
((1 - A){\Bbb Z}^N + (1 - B) {\Bbb Z}^N) \\
& \longrightarrow
K_0({\cal O}_{A,B}^\kappa) \\
& \longrightarrow
\Ker(1 - A) \cap \Ker(1 - B)
\text{ in }
{\Bbb Z}^N
\longrightarrow 0 \\
\intertext{and}
0 
& \longrightarrow 
\Ker(1-B) \text{ in } {\Bbb Z}^N /
(1 - A)(\Ker(1 - B) \text{ in } {\Bbb Z}^N )\\
& \longrightarrow
K_1({\cal O}_{A,B}^\kappa) \\
& \longrightarrow
\Ker(1 - A) 
\text{ in }\quad
{\Bbb Z}^N /
(1 - B) {\Bbb Z}^N 
\longrightarrow 0.
\end{align*}
\end{prop}
We consider $1 \times 1$ matrices
$[N]$ and $[M]$ with its entries $N$ and $M$  
respectively for
$1< N,M \in {\Bbb N}$.
Let
$G_N$ be a directed labeled graph with 
one vertex and $N$-self directed loops.
Similarly 
we consider a
directed labeled graph $G_M$
with 
 $M$-self loops at the vertex.
Denote by 
$\Sigma^N = \{ f_1,\dots,f_N\}$ 
the set of directed $N$-self loops of $G_N$
and
$\Sigma^M = \{ e_1,\dots,e_M\}$ 
the set of directed $M$-self loops of $G_M$.
The correspondence
$(e,f)\in \Sigma^M \times \Sigma^N \longrightarrow (f,e) \in \Sigma^N \times \Sigma^M$
yields a specification $\kappa$,
which we will fix.
Put
$$
\rho^M_{e_i}(1) = 1, \qquad
\rho^N_{f_j}(1) = 1.
$$
Then we have a $C^*$-textile dynamical system
\begin{equation*}
({\Bbb C}, \rho^M, \rho^N, \Sigma^M, \Sigma^N, \kappa).
\end{equation*}
The associated $C^*$-algebra is denoted by 
${\cal O}_{M,N}$.



\begin{lem}
${\cal O}_{N,M} = {\cal O}_N \otimes {\cal O}_M$.
\end{lem}
\begin{pf}
Let
$s_i, i=1,\dots,N$ and
$t_j, i=1,\dots,M$ and
be the generating isometries of ${\cal O}_N$ 
and of ${\cal O}_M$ satisfying
$$
\sum_{i=1}^N s_i s_i^* =1, \qquad \sum_{j=1}^M t_j t_j^* =1
$$
Let
$S_i, i=1,\dots,N$ and
$T_j, i=1,\dots,M$ and
be the generating isometries of ${\cal O}_{N,M}$ 
satisfying
$$
\sum_{i=1}^N S_i S_i^* =1, \qquad \sum_{j=1}^M T_j T_j^* =1
$$
such that
$$ 
S_i T_j = T_j S_i, \qquad i=1,\dots,N, \quad j=1,\dots,M.
$$ 
Since
$
(s_i \otimes 1)(1\otimes t_j) = (1\otimes t_j)  (s_i \otimes 1), 
i=1,\dots,N, \quad j=1,\dots,M,
$
By the universality of 
${\cal O}_{N,M}$
subject to the relations,
one has a surjective homomorphism
$
\Phi: {\cal O}_{N,M} \longrightarrow {\cal O}_N \otimes {\cal O}_M
$
such that
$
\Phi(S_i) = s_i \otimes 1,
\ 
\Phi(T_j) = 1 \otimes t_j.
$
And also 
by the universality of the tensor product 
$
{\cal O}_N \otimes {\cal O}_M,
$
there exists a homomorphism
$
\Psi: {\cal O}_N \otimes {\cal O}_M \longrightarrow {\cal O}_{N,M}
$
such that
$
\Psi(s_i \otimes 1) = S_i,
\ 
\Psi(1 \otimes t_j) = T_j.
$
Since
$
\Phi \circ \Psi = \id,
\Psi \circ \Phi = \id,
$
one concludes that
$\Phi$ and $\Psi$ are inverses to each other so that
$
{\cal O}_{N,M}\cong {\cal O}_N \otimes {\cal O}_M. 
$
As both 
$
{\cal O}_N
$
and
$
 {\cal O}_M 
$
are simple, purely infinite,
we have 
$
{\cal O}_{N,M}
$
is simple, purely infinite.
\end{pf}

\begin{lem}
Put
$n = N-1, m= M-1$.
Then we have
the subgroup
$\{[k] \in {\Bbb Z}/m{\Bbb Z} \mid n k \in m{\Bbb Z} \}$
of
${\Bbb Z}/m{\Bbb Z}$
is isomorphic to
${\Bbb Z}/d{\Bbb Z}$.
\end{lem}
\begin{pf}
As $d = g.c.d(n,m)$,
there exist $n_0, m_0 \in {\Bbb N}$ 
such that
$m = m_0 d, n= n_0 d
$
and
$(n_0,m_0) =1$.
For $k \in {\Bbb Z}$
with
$n k \in m{\Bbb Z}$,
the condition
$(n_0,m_0) =1$
implies
$ k = m_0 k'$
for some 
$k' \in {\Bbb Z}$.
Hence $k \in m_0 {\Bbb Z}$
so that we see that
the subgroup
$\{[k] \in {\Bbb Z}/m{\Bbb Z} \mid n k \in m{\Bbb Z} \}$
of
${\Bbb Z}/m{\Bbb Z}$
is isomorphic to
$m_0{\Bbb Z}/ m_0 d{\Bbb Z}$,
which is isomorphic to
${\Bbb Z}/d{\Bbb Z}$.
\end{pf}
\begin{lem}
For $ 1 <N,M\in {\Bbb N}$
 with $d = g.c.d(N-1,M-1)$,
\begin{enumerate}
\renewcommand{\labelenumi}{(\roman{enumi})}
\item
${\Bbb Z} / ((N-1){\Bbb Z} + (N-1){\Bbb Z}) \cong  {\Bbb Z}/d{\Bbb Z}.$
\item
$\Ker(N-1) = \Ker(M-1) = 0$ \text{ in } ${\Bbb Z}$.
\end{enumerate}
\end{lem}
\begin{pf}
It is easy  to show 
that
the subgroup 
$(N-1){\Bbb Z} + (N-1){\Bbb Z}$ of ${\Bbb Z}$
coincides with $d{\Bbb Z}$.
(ii) is trivial.
\end{pf}
Therefore we have
\begin{prop}
For $ 1 <N,M\in {\Bbb N}$,
the $C^*$-algebra
$
{\cal O}_{N,M}
$
is simple, purely infinite,
such that
\begin{equation*}
K_0({\cal O}_{N,M}) 
\cong 
K_1({\cal O}_{N,M})
\cong {\Bbb Z}/ d{\Bbb Z}
\end{equation*}
where $d = g.c.d(N-1,M-1) $
the greatest common diviser of $N-1,M-1$.
\end{prop}
It is easy to see that
 the K-groups $K_i( {\cal O}_N \otimes {\cal O}_M)$
are  ${\Bbb Z}/ d{\Bbb Z}$ for $i=0,1$
by using the K{\"u}nneth formula  proved in \cite{RS}.
 
We will generalize the above examles from the view point of tensor products. 

\medskip

{\bf 2. Tensor products.}

Let
$(\A^\rho,\rho,\Sigma^\rho)$
and 
$(\A^\eta,\eta,\Sigma^\eta)$
be $C^*$-symbolic dynamical systems.
We will construct a $C^*$-textile dynamical system as follows:
Put
\begin{equation*}
\bar{\A} = \A^\rho \otimes \A^\eta, \qquad
\bar{\rho}_\alpha = \rho_\alpha \otimes \id, \qquad
\bar{\eta}_a = \id \otimes \eta_a, \qquad
\Sigma^{\bar{\rho}}=\Sigma^\rho ,\qquad
\Sigma^{\bar{\eta}} =\Sigma^\eta
\end{equation*}
for
$\alpha \in \Sigma^\rho, a \in \Sigma^\eta$,
where
$\otimes$ means the minimal $C^*$-tensor product $\otimes_{\min}$.
For 
$(\alpha,a) \in \Sigma^\rho\times\Sigma^\eta$,
we see
$
\eta_b \circ \rho_\alpha (1) \ne 0
$ 
if and only if
$
\eta_b(1) \ne 0,  \rho_\alpha (1) \ne 0,
$ 
so that
$$
\Sigma_{\bar{\rho}\bar{\eta}} = \Sigma^\rho \times\Sigma^\eta
\quad
\text{ and similarly }
\quad
\Sigma_{\bar{\eta}\bar{\rho}} = \Sigma^\eta \times\Sigma^\rho.
$$
\begin{lem}
Define 
$\bar{\kappa}:\Sigma_{\bar{\rho}\bar{\eta}}\longrightarrow \Sigma_{\bar{\eta}\bar{\rho}}$
by setting
$\bar{\kappa}(\alpha,b) = (b,\alpha)$.
We then have
$(\bar{\A}, \bar{\rho}, \bar{\eta},\Sigma^{\bar{\rho}}, \Sigma^{\bar{\eta}}, \bar{\kappa})$
is a $C^*$-textile dynamical system.
\end{lem}
\begin{pf}
By \cite{Arch},
we have 
$Z_{\bar{\A}} = Z_{\A^\rho} \otimes Z_{\A^\eta}$
so that 
$$
\bar{\rho}_\alpha(Z_{\bar{\A}}) \subset Z_{\bar{\A}}, 
\quad
\alpha \in \Sigma^{\bar{\rho}}
\quad
\text{ and }
\quad
\bar{\rho}_a(Z_{\bar{\A}}) \subset Z_{\bar{\A}}, 
\quad
a \in \Sigma^{\bar{\eta}}.
$$
We also have
$
\sum_{\alpha \in\Sigma^{\bar{\rho}}} \bar{\rho}_\alpha(1) = 
\sum_{\alpha \in\Sigma^{\rho}} \rho_\alpha(1)\otimes 1 \ge 1,
$
and similarly
$
\sum_{a \in\Sigma^{\bar{\eta}}} \bar{\eta}_(1)  \ge 1
$
so that
both families  
$\{ \bar{\rho}_\alpha \}_{\alpha \in\Sigma^{\bar{\rho}}} \}$
and
$\{ \bar{\eta}_a \}_{a \in\Sigma^{\bar{\eta}}} \}$
of endomorphisms
are essential.
Since
$\{ \rho_\alpha \}_{\alpha \in \Sigma^\rho}$
is faithful on $\A^\rho$,
the homomorphism
$$
x \in \A^\rho \longrightarrow  
\oplus_{\alpha \in \Sigma^\rho} \rho_\alpha(x) \in \oplus_{\alpha \in \Sigma^\rho} \A^\rho
$$
is injective
so that the homomorphism
$$
x\otimes y  \in \A^\rho \otimes \A^\eta \longrightarrow  
\oplus_{\alpha \in \Sigma^\rho} \rho_\alpha(x)\otimes y  \in \oplus_{\alpha \in \Sigma^\rho} \A^\rho\otimes \A^\eta
$$
is injective.
This implies that
$\{ \bar{\rho}_\alpha \}_{\alpha \in\Sigma^{\bar{\rho}}}$
is faithful and similary so is
$\{ \bar{\eta}_a \}_{a \in\Sigma^{\bar{\eta}}}$.
Hence 
$
(\bar{\A}, \bar{\rho}, \Sigma^{\bar{\rho}})
$
and
$(\bar{\A},  \bar{\eta}, \Sigma^{\bar{\eta}})$
are $C^*$-symbolic dynamical systems.
It is direct to see that
$\bar{\eta}_b \circ \bar{\rho}_\alpha = \bar{\rho}_\alpha\circ \bar{\eta}_b
$
for
$(\alpha,b) \in \Sigma_{\bar{\rho}\bar{\eta}}$.
Therefore
$(\bar{\A}, \bar{\rho}, \bar{\eta},\Sigma^{\bar{\rho}}, \Sigma^{\bar{\eta}}, \bar{\kappa})$
is a $C^*$-textile dynamical system.
\end{pf}
We call 
$(\bar{\A}, \bar{\rho}, \bar{\eta},\Sigma^{\bar{\rho}}, \Sigma^{\bar{\eta}}, \bar{\kappa})$
the tensor product between
$(\A^\rho,\rho,\Sigma^\rho)$
and
$(\A^\eta,\eta,\Sigma^\eta)$.
Denote by 
$S_\alpha, \alpha \in \Sigma^{\bar{\rho}}$,
$T_a, a \in \Sigma^{\bar{\eta}}$
the generating partial isometries 
of the $C^*$-algebra
${\cal O}_{\bar{\rho},\bar{\eta}}^{\bar{\kappa}}$
for the $C^*$-textile dynamical system
$(\bar{\A}, \bar{\rho}, \bar{\eta},\Sigma^{\bar{\rho}}, \Sigma^{\bar{\eta}}, \bar{\kappa})$.
By the universality for the algebra
${\cal O}_{\bar{\rho},\bar{\eta}}^{\bar{\kappa}}$
subject to the relations
$(\bar{\rho},\bar{\eta},\bar{\kappa})$,
one sees that the algebra
${\cal D}_{\bar{\rho},\bar{\eta}}$
is isomorphic to
the tensor product
${\cal D}_\rho \otimes {\cal D}_\eta$
through the correspondence
$$
S_\mu T_\xi(x \otimes y) T_\xi^* S_\mu^* \longleftrightarrow S_\mu x S_\mu^* \otimes T_\xi y T_\xi^*
$$
for
$\mu \in B_*(\Lambda_\rho), \xi\in B_*(\Lambda_\eta), \ $
$x \in A^\rho, y \in \A^\eta$.

\begin{lem}
Suppose that 
$(\A^\rho,\rho,\Sigma^\rho)$
and
$(\A^\eta,\eta,\Sigma^\eta)$
are both free (resp. AF-free).
Then the tensor product
$(\bar{\A}, \bar{\rho}, \bar{\eta},\Sigma^{\bar{\rho}}, \Sigma^{\bar{\eta}}, \bar{\kappa})$
is free (resp. AF-free).
\end{lem}
\begin{pf}
Suppose that 
$(\A^\rho,\rho,\Sigma^\rho)$
and
$(\A^\eta,\eta,\Sigma^\eta)$
are both free.
There exist 
increasing sequences 
$\A^\rho_l, l \in \Zp$ 
and 
$\A^\eta_l, l \in \Zp$
of $C^*$-subalgebras
of $\A^\rho$ and $\A^\eta$
satisfying the conditions of the freeness respectively.
Put
$\bar{\A}_l = \A^\rho_l \otimes \A^\eta_l, l\in \Zp $
It is clear that

(1) $\bar{\rho}_\alpha(\bar{A}_l) \subset \bar{\A}_{l+1}, \alpha \in \Sigma^{\bar{\rho}}$
    and 
    $ \bar{\eta}_a(\bar{A}_l) \subset \bar{\A}_{l+1}, a \in \Sigma^{\bar{\eta}}$ for $l \in \Zp$.
 
(2) $\cup_{l \in \Zp}\bar{\A}_l $ is dense in $\bar{\A}$. 

We will show that the condition (3) in the definition of freeness holds.
Take and fix arbitrary  $j,k,l\in {\Bbb N}$ with $j+k \le l$.
For 
$j \le l$, by the freeness of 
$(\A^\rho,\rho,\Sigma^\rho)$
there exists a projection 
$q_\rho \in {\cal D}_{\rho} \cap {\A^\rho_l}^\prime 
$
such that
\begin{enumerate}
\renewcommand{\labelenumi}{(\roman{enumi})}
\item
 $ q_\rho x \ne 0$ for $0\ne x \in \A^\rho_l$,
\item
 $
 \phi_{\rho}^n(q_\rho) q_\rho = 0
$ 
for all $n=1,2,\dots,j.$
\end{enumerate} 
Similarly
for 
$k \le l$, by the freeness of 
$(\A^\eta,\eta,\Sigma^\eta)$
there exists a projection 
$q_\eta \in {\cal D}_{\eta} \cap {\A^\eta_l}^\prime 
$
such that
\begin{enumerate}
\renewcommand{\labelenumi}{(\roman{enumi})}
\item
 $ q_\eta y \ne 0$ for $0\ne y \in \A^\eta_l$,
\item
 $
 \phi_\eta^m(q_\eta) q_\eta = 0
$ 
for all $m=1,2,\dots,k.$
\end{enumerate} 
Put
$q = q_\rho \otimes q_\eta \in {\cal D}_\rho \otimes {\cal D}_\eta (= {\cal D}_{\bar{\rho},\bar{\eta}})$
so that
$q \in {\cal D}_{\bar{\rho},\bar{\eta}}\cap \bar{\A}_l^\prime$.
As the maps
$\Phi_l^\rho: x \in \A_l^\rho \longrightarrow q_\rho x \in q_\rho \A_l^\rho$
and 
$\Phi_l^\eta: y \in \A_l^\eta \longrightarrow q_\eta x \in q_\eta \A_l^\eta$
are isomorphisms
so that the tensor product
$$
\Phi_l^\rho\otimes \Phi_l^\eta : 
x\otimes y  \in \A_l^\rho \otimes \A_l^\eta 
\longrightarrow 
(q_\rho\otimes q_\eta)(x\otimes y)  \in 
(q_\rho\otimes q_\eta) \in \A_l^\rho\otimes \A_l^\eta
$$
is 
isomorphic.
Hence
$qa \ne 0$ for $0\ne a \in \bar{\A}_l$.
It is straightforward to see that 
 $
 \phi_{\rho}^n(q) \phi_{\eta}^m(q) 
=\phi_{\rho}^n(( \phi_{\eta}^m(q)))q 
=\phi_{\rho}^n(q) q 
=\phi_{\eta}^m(q) q = 0
$ 
for all $n=1,2,\dots,j$, $m= 1,2,\dots, k.$
Therefore 
the tensor product
$(\bar{\A}, \bar{\rho}, \bar{\eta},\Sigma^{\bar{\rho}}, \Sigma^{\bar{\eta}}, \bar{\kappa})$
is free.
It is obvious to see that 
if 
both
$(\A^\rho,\rho,\Sigma^\rho)$
and
$(\A^\eta,\eta,\Sigma^\eta)$
are AF-free,
then 
$(\bar{\A}, \bar{\rho}, \bar{\eta},\Sigma^{\bar{\rho}}, \Sigma^{\bar{\eta}}, \bar{\kappa})$
is AF-free.
\end{pf}
\begin{prop}
Suppose that
$(\A^\rho,\rho,\Sigma^\rho)$
and
$(\A^\eta,\eta,\Sigma^\eta)$
are both free.
Then the $C^*$-algebra 
${\cal O}_{\bar{\rho},\bar{\eta}}^{\bar{\kappa}}$
for the tensor product $C^*$-textile dynamical system
$(\bar{\A}, \bar{\rho}, \bar{\eta},\Sigma^{\bar{\rho}}, \Sigma^{\bar{\eta}}, \bar{\kappa})$
is isomorphic to the tensor product 
${\cal O}_\rho \otimes {\cal O}_\eta$
of the $C^*$-algebras 
between
${\cal O}_\rho$
and
${\cal O}_\eta$.
If in particular,
$(\A^\rho,\rho,\Sigma^\rho)$
and
$(\A^\eta,\eta,\Sigma^\eta)$
are both irreducible,
the $C^*$-algebra 
${\cal O}_{\bar{\rho},\bar{\eta}}^{\bar{\kappa}}$
is simple.
\end{prop}
\begin{pf}
Suppose that
$(\A^\rho,\rho,\Sigma^\rho)$
and
$(\A^\eta,\eta,\Sigma^\eta)$
are both free.
By the preceding lemma,
the tensor product 
$(\bar{\A}, \bar{\rho}, \bar{\eta},\Sigma^{\bar{\rho}}, \Sigma^{\bar{\eta}}, \bar{\kappa})$
is free and hence satisfies condition (I).
Let
$s_\alpha, \alpha \in \Sigma^\rho$
and
$t_a, a \in \Sigma^\eta$
be the generating partial isometries of the $C^*$-algebras
${\cal O}_\rho$
and
${\cal O}_\eta$
respectively.
Let
$S_\alpha, \alpha \in \Sigma^{\bar{\rho}}$
and
$ T_a, a \in \Sigma^{\bar{\eta}}$
be the generating partial isometries 
of the $C^*$-algebra
${\cal O}_{\bar{\rho},\bar{\eta}}^{\bar{\kappa}}$.
By the uniqueness of the algebra
${\cal O}_{\bar{\rho},\bar{\eta}}^{\bar{\kappa}}$
with respect to the relations $(\bar{\rho}, \bar{\eta}, \bar{\kappa})$,
the correspondence
$$
S_\alpha \longrightarrow s_\alpha \otimes 1 \in {\cal O}_\rho \otimes {\cal O}_\eta,
\qquad
T_a \longrightarrow  1 \otimes t_a  \in {\cal O}_\rho \otimes {\cal O}_\eta
$$
naturally gives rise to an isomorphism from
${\cal O}_{\bar{\rho},\bar{\eta}}^{\bar{\kappa}}$
onto the tensor product 
${\cal O}_\rho \otimes {\cal O}_\eta$.

If in particular,
$(\A^\rho,\rho,\Sigma^\rho)$
and
$(\A^\eta,\eta,\Sigma^\eta)$
are both irreducible,
the $C^*$-algebras
${\cal O}_\rho$
and
${\cal O}_\eta$
are both simple
so that
${\cal O}_{\bar{\rho},\bar{\eta}}^{\bar{\kappa}}$
is simple.
\end{pf}

We remark that 
the tensor product
$(\bar{\A}, \bar{\rho}, \bar{\eta},\Sigma^{\bar{\rho}}, \Sigma^{\bar{\eta}}, \bar{\kappa})$
does not necessarily form square.
The  K-theory groups 
$K_*({\cal O}_{\bar{\rho},\bar{\eta}}^{\bar{\kappa}})$
are computed from the  
 K{\"u}nneth formulae for 
$K_*({\cal O}_\rho \otimes {\cal O}_\eta)$
\cite{RS}.

In \cite{Ma2011},
higher dimensional analogue
$(\A, \rho^1,\dots,\rho^N, \Sigma^1,\dots,\Sigma^N,\kappa)$
will be studied.

\end{document}